\documentclass[final,onefignum,onetabnum]{siamart220329}

\usepackage[american]{babel}
\usepackage{lipsum}
\usepackage{amsfonts}
\usepackage{graphicx}
\usepackage{epstopdf}
\usepackage{algorithm}
\usepackage{amsmath,bm}
\usepackage{amssymb}
\usepackage{faktor} 
\usepackage{algorithmic}
\usepackage{graphicx}
\usepackage{subcaption}
\usepackage{xr-hyper}
\usepackage{hyperref}
\hypersetup{
	colorlinks=true,
	linkcolor=black,
	filecolor=black,
	citecolor=black,      
	urlcolor=black,
}
\usepackage[aboveskip=3pt]{caption} 
\Crefname{ALC@unique}{Line}{Lines}
\ifpdf
\DeclareGraphicsExtensions{.eps,.pdf,.png,.jpg}
\else
\DeclareGraphicsExtensions{.eps}
\fi
\DeclareMathOperator*{\argmax}{argmax}
\DeclareMathOperator*{\argmin}{argmin}

\newsiamremark{remark}{Remark}
\newsiamremark{hypothesis}{Hypothesis}
\crefname{hypothesis}{Hypothesis}{Hypotheses}
\newsiamthm{claim}{Claim}

\headers{Entropic trust region}
{M. Torda, J.Y. Goulermas, R. P\'{u}\v{c}ek, and V. Kurlin}

\title{Entropic trust region for densest crystallographic symmetry group packings \thanks{\funding{This work was funded by the Leverhulme Research Centre for Functional Materials Design and supported by EPSRC grants EP/R018472/1, EP/X018474/1, and RAEng fellowship IF2122$\backslash$186.}
	}
	}

\author{Miloslav Torda
	 \thanks{
		Leverhulme Research Centre for Functional Materials Design, Department of Computer Science, University of Liverpool, Liverpool, UK 
		(\email{miloslav.torda@liverpool.ac.uk})}
	\and John Y. Goulermas \thanks{
		Department of Computer Science, University of Liverpool, Liverpool, UK (\email{j.y.goulermas@liverpool.ac.uk}, \email{vitaliy.kurlin@liverpool.ac.uk})}
	\and Roland P\'{u}\v{c}ek \thanks{Department of Mathematics, Friedrich Schiller University Jena, Jena, DE (\email{roland.pucek@uni-jena.de})}
	\and Vitaliy Kurlin \footnotemark[3]
}

\usepackage{amsopn}

\begin{document}

	\maketitle
	
	\maketitle
	
	\begin{abstract}	 	
		Molecular crystal structure prediction (CSP) seeks the most stable periodic structure given the chemical composition of a molecule and pressure-temperature conditions. Modern CSP solvers use global optimization methods to search for structures with minimal free energy within a complex energy landscape induced by intermolecular potentials. A major caveat of these methods is that initial configurations are random, making thus the search susceptible to convergence at local minima. Providing initial configurations that are densely packed with respect to the geometric representation of a molecule can significantly accelerate CSP. Motivated by these observations, we define a class of periodic packings restricted to crystallographic symmetry groups (CSG) and design a search method for the densest CSG packings in an information-geometric framework. Since the CSG induce a toroidal topology on the configuration space, a non-Euclidean trust region method is performed on a statistical manifold consisting of probability distributions defined on an $n$-dimensional flat unit torus by extending the multivariate von Mises distribution. Introducing an adaptive quantile reformulation of the fitness function into the optimization schedule provides the algorithm with a geometric characterization through local dual geodesic flows. Moreover, we examine the geometry of the adaptive selection-quantile defined trust region and show that the algorithm performs a maximization of stochastic dependence among elements of the extended multivariate von Mises distributed random vector. 
		We experimentally evaluate the behavior and performance of the method on various densest packings of convex polygons in $2$-dimensional CSGs for which optimal solutions are known, and demonstrate its application in the pentacene thin-film CSP.
	\end{abstract}
	
	

	\section{Introduction}
	\label{sec:introduction}	
\begin{figure}[t]
		\centering
		\begin{subfigure}{0.16\textheight}
			\includegraphics[trim={0 -92 0 0},clip,width=\columnwidth]{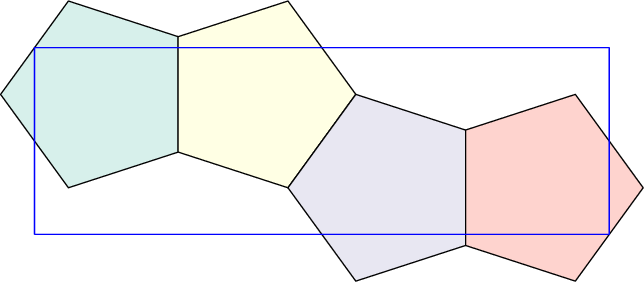}
		\end{subfigure}
		\begin{subfigure}{0.45\textheight}
			\includegraphics[trim={0 0 0 0},clip,width=\columnwidth]{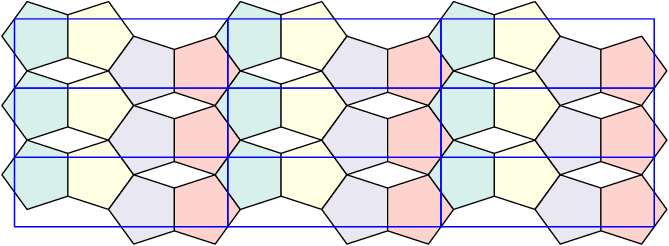}
		\end{subfigure}
		\caption{The $2$D periodic structure with the $p2mg$ plane group symmetry (i.e., consisting of a $2$-fold rotational symmetry operation, glide reflections and mirror reflections along two mirror planes. Colors represent these symmetry operations modulo lattice translations.) where $K$ is a regular convex pentagon with the packing density of approximately 0.8541019. (Left) A single unit cell. (Right) $15$ unit cells.}
		\label{fig:P2G3Gon_example}
	\end{figure}

	The work presented here is motivated by the problem of Crystal Structure Prediction (CSP), in which, given some molecular shape $K$, the goal is to predict a synthesizable periodic structure. Such a periodic structure may consist of several copies of $K$ within a unit cell formation (parallelepiped) that is periodically repeated along the three directions. \cref{fig:P2G3Gon_example} exemplifies such a formation with a 2D pentagonal crystal. CSP traditionally starts from an almost random configuration of molecules in a random unit cell and attempts to optimize a complex energy function depending on the given molecular structure and numerous problem parameters such as pressure-temperature conditions.
	
	Current CSP approaches have two main computational caveats. The first is energy computation, where either one of many empirical potentials needs to be chosen or computationally expensive, but precise density-functional theory calculations are used. The second is that the energy functions induce complicated energy landscapes \cite{wales2018}, increasing the likelihood of the global search methods converging to a local minimum basin and leading to over-prediction \cite{price2018}. The standard output of CSP computations is thousands of theoretical polymorphic structures, each representing some local optimum of the energy landscape \cite{woodley2020}. Afterwards, data analytic tools are employed to identify metastable structures. 
	
	Since crystals are solid materials, they are almost always very dense, and therefore CSP can be substantially accelerated when initial configurations are sufficiently dense and not random.
	Maximizing atomic densities has already been considered in current CSP software \cite{de2012, glass2006}, where various free energy approximations can be minimized, such as Lennard-Jones potentials, Buckingham potentials, and others. 
	
	In this work, we propose an approach based on discrete geometry to facilitate the CSP workflow by providing reasonable initial periodic configurations. Specifically, a polytope representation is assigned to a molecule based on its intrinsic properties \cite{Santolini, wang2021}, and the highest density packings of hypothetical structures are generated. These hypothetical configurations are then used as starting positions for the usual CSP energy minimizations, thus reducing the computational burden only to local explorations of the energy landscape.
	
	Geometric packings are well-studied objects in discrete and computational geometry \cite{toth} and are fundamental in solid-state physics modelling \cite{torquato2018}. Although molecular crystals are considered to be embeddings in the 3D Euclidean space, following the currently high interest in 2D materials \cite{das2021}, the same approach can also be employed here. That is, given a representation of a molecule by a polygon, the aim is to acquire the configuration that maximizes packing density and subsequently use this as a starting configuration in classical CSP workflows.  	
	Moreover, not only generally densest packings but also lower density but higher symmetry crystal structures that maximize packing density among a particular isomorphism class of periodic structures can point to possible, stable crystal phases. For example, \Cref{fig:P2G3Gon_example} illustrates a crystallization pattern of pentagonal proteins on lipid mono-layers \cite{wang1999} by densest packing of pentagon when the configuration space is restricted to the $p2mg$ plane group.
	
	This approach is well justified experimentally \cite{cui2018, ecija2012} and by previous molecular dynamics simulations using force-field methods \cite{mukhopadhyay, zhao}.
	Furthermore, the crystallization conjecture states that in the Euclidean space of dimensions two and three, the ground state energy of systems of interacting particles forms periodic configurations in the thermodynamic limit \cite{lewin2015}. \cite{theil2006} proved the equivalence between the crystallization conjecture for mono-atomic systems in two dimensions and the densest disc packing for a class of Lennard-Jones-like energy potentials. Later, \cite{flatley2015} proved the face-centred-cubic sphere packing model's optimality in terms of energy minimization using an additional three-body potential. Even though there are no such results for molecular systems, the usual correlation between packing density maximization and energy minimization suggests an equivalence between the densest packings of polytopes and the crystallization conjecture, at least for some molecular crystals.
	
	Despite having attracted the interest of various scientific communities for centuries now, constructing the densest packings of a set of given geometric entities is a notoriously hard problem, and only a few optimal solutions are known. In the $3$D Euclidean space, these include the general packing of the sphere \cite{hales2005} and the truncated rhombic dodecahedron \cite{bezdek1994} and the existence of an algorithm to construct the densest lattice packings of convex polygons \cite{betke2000}. In the 2D Euclidean plane, known general packings include that of the disk \cite{toth2013} and the pentagon \cite{hales}, algorithms to construct the densest packing of centrally symmetric convex polygons \cite{toth2013}, and algorithms for the densest lattice packings and double lattice packings of convex polygons \cite{mount1990, mount1991}.
	However, all aforementioned construction methods are tailored to a specific geometric shape. We aim to construct the densest packings for a large class of objects and symmetry groups using a robust and generic method without the restriction to given shapes. To this extent, we have developed an optimization framework based on the natural gradient method \cite{amari} used in evolution strategies \cite{hansen2001, wierstra} which are instances of a general information-geometric optimization framework \cite{ollivier}.
	
	Although our proposed optimization system was specifically developed to search for the densest crystallographic packings, it can also be applied in classical CSP computations by interfacing with force field methods \cite{gale2003} or density functional theory calculations \cite{hafner2008}. This can be achieved by replacing the maximizing packing density objective with minimizing free energy. In essence, due to the generic design of the entropic trust region based optimization, it can be used to solve any bounded and constrained black-box optimization problem.
	
	Since the currently established molecular crystal model is that of a crystallographic symmetry group (CSG), we restrict our search to the $230$ space groups and $17$ plane groups \cite{brock2016} and define a new subproblem of the general packing problem \cite{rogers}, the CSG packing. The periodic boundary conditions inherent to CSGs, innately induce a toroidal topology on the packing configuration space. We exploit this property by performing natural gradient ascent on a statistical manifold composed of probability distributions on an $n$D torus using an extension of the multivariate von Mises model for directional data \cite{mardia2008}. This effectively removes the optimization problem boundaries, which pose a considerable difficulty for many optimization methods when the solution lies on the configuration space boundaries. Moreover, by starting from the uniform distribution on the torus, we remove the algorithm's dependence on the initial configuration.
	
	The manuscript is divided into four parts. First, \Cref{sec:spacegrouppacking} introduces CSG packings and the related densest CSG packing problem. \Cref{sec:entTrust} introduces the entropic trust region method and the extended multivariate von Mises distribution and presents the exponential family reformulation for the extended multivariate von Mises distribution. We also introduce the adaptive selection quantile to facilitate the natural gradient ascent. Furthermore, using its connection with the proximal entropic method \cite{teboulle}, we examine the geometry of the entropic trust region based exponential family adaptive selection quantile and establish that the algorithm defines, in fact, parallel Riemannian gradient flows between two statistical submanifolds of an ambient manifold; namely, one maximizing entropy \cite{shore} and another maximizing multi-information \cite{Studeny1998}. Finally, the section introduces a method for refining solutions based on the localization of the search in a subspace of the initial configuration space. \Cref{sec:proofofconcept} contains various experiments that investigate the behaviour and performance of the entropic trust region method on the densest $p2$ packing of a regular octagon. Results demonstrate higher accuracy than in previous Monte--Carlo methods \cite{atkinson, schilling2005}.
	Lastly, in \Cref{sec:pentacene}, we demonstrate the application of our algorithm on the densest packings of a polygonal representation of pentacene in a few plane groups to determine crystal structures of pentacene thin films found in the literature. Many additional experiments and mathematical, algorithmic and implementation details are included in the supplementary material for this manuscript.

	\section{Densest CSG packings}
	\label{sec:spacegrouppacking}
	Even though regarding CSP applications, we are interested in dimensions $2$ and $3$, we introduce CSGs  for general $n$ dimensions following \cite{miller}.
	
	Let $E_n$ denote the group of all isometries in $\mathbb{R}^n$ that is all length preserving transformations of $\mathbb{R}^n$, and $O_n$ and $T_n$ denote orthogonal and translational subgroups of $E_n$, respectively. Specifically
	\begin{math}
		O_n=\{\textbf{R} \ | \ \textbf{R}^{ \intercal} \textbf{R} = \textbf{R} \textbf{R}^\intercal = \textbf{I};\textbf{R} \in \mathbb{R}^{n\times n} \}
	\end{math}
	is the set of all rotations and rotoinversions that preserve a fixed point $\mathbf{p} \in \mathbb{R}^n$, and
	\begin{math}
		T_n=\{T_{\mathbf{a}} \ | \ T_{\mathbf{a}}\mathbf{p}= \mathbf{p} + \mathbf{a} \ ; \mathbf{a} \in \mathbb{R}^n  \}
	\end{math}
	is the set of all translations.
	
	The $n$D lattice group
	\begin{equation}\label{eq:latticeG}
		L=\left\{ T_{\mathbf{a}} \ | \ \mathbf{a}= u_{1}\textbf{b}_{1}+\ldots + u_{n}\textbf{b}_{n} \ ; \ u_{1},\ldots , u_{n}  \in \mathbb{Z}\right\},
	\end{equation}
	is a discrete subgroup of $T_n$ where
	$\textbf{b}_{1},\ldots , \textbf{b}_{n} \in \mathbb{R}^n$ are the basic vectors. We denote the orbit of a point $\mathbf{p}$ under the action of $L$ by $\Lambda_{\mathbf{p}}=\left\{T_{\mathbf{a}}\mathbf{p} \ | \ T_{\mathbf{a}} \in L \right\}$ (the geometric lattice) and the convex hull of basic vectors by $\bar{\Lambda}=\text{conv}\left\{\mathbf{0}.\textbf{b}_{1},\ldots,\textbf{b}_{n}  \right\}$ (the primitive cell). The crystallographic point group is a finite subgroup of $O_n$ that maps lattice $\Lambda_{\mathbf{0}}$ to itself. The CSG $G$ is a discrete subgroup of $E_n$, such that $L = G \ \cap \ T_n$ is a lattice group of the form \cref{eq:latticeG} or, in other words, a discrete subgroup of the group of isometries of the $n$-dimensional Euclidean space containing a lattice subgroup. In the following, we denote the lattice $L$ associated with a CSG $G$ as $L_{G}$ and the primitive cell corresponding to $L_{G}$ as $\bar{\Lambda}_{G}$. An asymmetric unit is a subset of the primitive cell such that the whole $\mathbb{R}^n$ is filled
	when the CSG symmetry operations are applied.	
	We refer to the equivalence classes of $2$D CSGs as plane-group types
	and of $3$D CSGs as space-group types. The $2$D CSGs are classified into $17$ plane-group
	types which are assigned to $10$ point group conjugacy classes (called 
	geometric crystal classes). The $230$ space-group types are classified into $32$ geometric crystal
	classes.
	
	
	Given an $n$D CSG $G \in \mathcal{G}$, where $\mathcal{G}$ is a CSG equivalence class associated with $G$, and $K$ a compact subset of $\mathbb{R}^n$, by a CSG packing $\mathcal{K}_{G}$, we mean a collection of non-overlapping copies of $K$ generated as an orbit under the action of $G$-action on $\mathbb{R}^n$.
	Formally, $\mathcal{K}_{G}$ is a $G$-set defined as
	\begin{subequations}\label{eq:spacepacking}
		\begin{gather}
			\mathcal{K}_{G}=\bigcup_{g \in G} g K, \label{eq:spacepackingConfig} \\ 
			\text{int}\left(g_i  K\right) \ \cap \text{int}\left(g_j  K\right)=\emptyset, \quad \forall \ g_i, \ g_j \in G, \; \ g_i\neq \ g_j. \label{eq:spacepackingConst}
		\end{gather}
	\end{subequations}
	
	Every element $g$ of the CSG $G$ acting on some point $\mathbf{p}$ can be expressed as
	\begin{equation}\label{eq:gact}
		g\mathbf{p} = \textbf{R}\mathbf{p} + \mathbf{ a} + \mathbf{ l},
	\end{equation}
	where $\textbf{R} \in O_n$, $\mathbf{ a} \in T_{\mathbf{ a}}$ such that
	\begin{equation}\label{eq:alat}
		\mathbf{a}=\alpha_1\mathbf{ b}_1+\ldots + \alpha_n\mathbf{ b}_n
	\end{equation}
	where $\alpha_i \in \mathbb{R}, \ 0\leq\alpha_i < 1 $ and $\mathbf{ b}_1,\ldots  , \ \mathbf{ b}_n$ lattice basis vectors, and $\mathbf{ l} \in L_{G}$.	
	Since $\mathcal{K}_{G}$ is a periodic system of sets it can be expressed in terms of \cref{eq:gact} as
	\begin{math}
		\textbf{R}_{i} K + \mathbf{a}_{i} + \mathbf{l}_{j},
	\end{math}
	with a finite number of translation vectors  $\mathbf{a}_{i}$ of form \cref{eq:alat}, rotations and rotoinversions $R_{i} \in O_n$ for $i=1,\ldots, N$, and lattice translation vectors $\mathbf{l}_{j} \in L_{G}$ for $j=1,2,\ldots$. Following from the formula for the packing density of the periodic system \cite{rogers}, the plane group packing density has a simple closed form expression
	\begin{equation}\label{eq:density}
		\rho \left(\mathcal{K}_G \right) = \frac{N m(K)}{m(\overline{\Lambda}_{G})},
	\end{equation}
	where $N$ is the number of symmetry operations modulo lattice translations in a CSG $G$ given by the pair $(\textbf{R},\mathbf{ a})$ in \cref{eq:gact}, $\overline{\Lambda}_{G}$ is the primitive cell, and $m(\cdot)$ denotes the $n$D Jordan measure.
	
	Finally, we can state the crystallographic packing problem. Given a CSG isomorphism class $\mathcal{G}$ and $K$, a compact subset of $\mathbb{R}^n$, the goal is to find the CSG packing $\mathcal{K}_{G_{\max}}$ with maximum density over the whole $\mathcal{G}$. Formally expressed, we aim at finding a $\mathcal{G}$-packing, such that  
	\begin{equation}\label{eq:problemPack}
		\mathcal{K}_{G_{\max}}=\argmax_{\mathcal{K}_{G \in \mathcal{G}}}\rho \left(\mathcal{K}_G \right).
	\end{equation}
	Since the Jordan measure of $K$ and the number of symmetry operations $N$ for a given $n$D CSG in \cref{eq:density} are constant, maximizing density is equivalent to minimizing the Jordan measure of the primitive cell associated with the geometric crystal class of $\mathcal{G}$. The crystallographic packing problem can be restated as finding a $\mathcal{G}$-packing with minimal primitive cell volume
	\begin{equation*} 
		\mathcal{K}_{G_{\max}}=\argmin_{\mathcal{K}_{G \in \mathcal{G}}}m(\overline{\Lambda}_{G}).
	\end{equation*}
	We refer to the solution of \cref{eq:problemPack} as the densest $\mathcal{G}$-packing. Here we search not only over the whole $\mathcal{G}$ but also over all rotations and translations of $K$, whose
	centroid lies in the asymmetric unit of $G$ such that the resulting configuration is a CSG packing.
	
	We consider the densest $\mathcal{G}$-packing as a nonlinear bounded constrained optimization problem. The Jordan measure of the primitive cell is computed as the determinant of lattice generators \cref{eq:latticeG}, a polynomial of degree $n$. The bounds are given by the space group's asymmetric unit, the range of rotational freedom of the set $K$, and bounds of the size and shape of the primitive cell given by the geometric crystal class associated with $\mathcal{G}$. This means that the configuration space is compact, which implies that problem \cref{eq:problemPack} has a solution. Linear and nonlinear constraints are given by the CSG's asymmetric unit and non-overlap condition \cref{eq:spacepackingConst}, respectively.
	
	\section{Entropic trust region}
	\label{sec:entTrust}
	In the classic black box optimization setting, no requirements are imposed on the objective except that given any design parameters or query $\textbf{x}$ from the set of possible configurations $\mathcal{X}$, the function value can be calculated. This kind of objective, which is here denoted by $\textbf{F}$, is frequently referred to as a zero-order oracle due to the lack of accessibility to gradient information. In all of the following, we assume that $\textbf{F}$ is a function between Borel measurable spaces $\mathcal{X}$ and $\mathcal{Y}$ with measures $\mu$ and $\nu$, respectively, such that $\mu(\textbf{F}^{-1}(E))=0$ whenever $\nu(E)=0$ for $E \subset \mathcal{Y}$.
	
	An approach for solving such optimization problems is stochastic relaxation \cite{geman1984}, and it follows from the observation that optimizing $\textbf{F}$ is equivalent to solving
	\begin{equation}\label{eq:Jmax}
		\bm{\tilde{ \theta}} = \argmax_{\bm{\theta} \in \bm{\Theta}} J(\bm{\theta}),
	\end{equation}
	with 
	\begin{equation}\label{eq:efit}
		J(\bm{\theta}):=E[\textbf{F}|\bm\theta]=\int_{\mathcal{X}} \textbf{F(x)}dP(\bm{\theta})
	\end{equation}
	being the expected value of $\textbf{F}$ under some probability measure $dP(\bm{\theta})$ from a parametric family of probability measures
	\begin{equation}\label{eq:Smodel}
		S=\left\{dP(\bm{\theta}) \ | \ \bm{\theta} \in \bm{\Theta} \subseteq \mathbb{R}^n \right\}
	\end{equation}
	defined on some configuration space $\mathcal{X}$.
	
	Natural evolution strategies \cite{wierstra} solve \cref{eq:Jmax} algorithmically using a non-Euclidean trust region strategy where the candidate step $\delta\bm{\theta}$ is found by maximizing the first-order Taylor approximation of $J(\bm{\theta}+\delta\bm{\theta})$. The trust region radius is given by the square root of twice the second-order approximation of the Kullback-Leibler divergence from $P_{\bm{\theta}}$ to $P_{\bm{\theta} + \delta\bm{\theta}}$ in $S$. Specifically
	\begin{subequations}\label{eq:trust}
		\begin{gather}
			\max_{\delta\bm{\theta}} J(\bm{\theta}^{t}) +  \delta\bm{\theta}^\intercal \nabla_{\bm{\theta}} J(\bm{\theta}^{t}) \\ \label{eq:trustcon}
			\text{s.t.} \ \sqrt{2 D_{KL}\left(P_{\bm{\theta}} \ || \ P_{\bm{\theta} +\delta\bm{\theta}} \right)}\approx\sqrt{\delta\bm{\theta}^\intercal\mathcal{I}_{\bm{\theta}} \delta\bm{\theta}}\leq \Delta^{t},
		\end{gather}
	\end{subequations}
	where $\nabla_{\bm{\theta}}$ denotes the Euclidean gradient operator with respect to $\bm{\theta}$ coordinates,  $\Delta^{t}$ is the trust region radius, and $\mathcal{I}_{\bm{\theta}}$ is the Fisher information matrix with elements 
	\begin{equation} \label{eq:FIM}
		{\mathcal{I}_{\bm{\theta}}}_{ij}=\int_{\mathcal{X}}  \frac{\partial\ln\left(p(\bm{\theta})\right) }{\partial \theta_i} \frac{\partial\ln\left(p(\bm{\theta})\right)}{\partial \theta_j}dP(\bm{\theta});
	\end{equation}
	$p(\bm{\theta})=\frac{dP(\bm{\theta})}{d\bm{\nu}}$ is the Radon-Nikodym derivative of $P(\bm{\theta})$ with respect to some reference measure $\bm{\nu}$ defined on $\mathcal{X}$, and the Kullback-Leibler divergence (KLD) from $\bm{\theta}$ to $\bm{\theta} + \delta\bm{\theta} $ is given by
	\begin{equation}\label{eq:KLdiv}
		D_{KL}\left(P_{\bm{\theta}} \ || \ P_{\bm{\theta} +\delta\bm{\theta}} \right)=\int_{\mathcal{X}} \ln \left(\frac{dP(\bm{\theta})}{dP(\bm{\theta}+ \delta\bm{\theta})}\right)dP(\bm{\theta})
	\end{equation}	
	where $P_{\bm{\theta}} \in S$ is parametrized by $\bm{\theta}$.
	
	Using Lagrange multipliers to solve \cref{eq:trust} results in a trust region step size
	\begin{equation}\label{eq:truststep}
		\delta\bm{\theta}^t=\frac{\Delta^{t} \mathcal{I}_{\bm{\theta}}^{-1} \nabla_{\bm{\theta}} J(\bm{\theta}^{t})}{\sqrt{{\nabla_{\bm{\theta}} J(\bm{\theta}^{t})}^\intercal \mathcal{I}_{\bm{\theta}}^{-1} \nabla_{\bm{\theta}} J(\bm{\theta}^{t})}}.
	\end{equation}	
	
	Expression \cref{eq:truststep} also has a differential geometric interpretation.  By considering the statistical model $S$ of \cref{eq:Smodel}, then
	\begin{equation}\label{eq:unitGrad}
		\text{grad} J(\bm{\theta}) =\frac{ \mathcal{I}_{\bm{\theta}}^{-1} \nabla_{\bm{\theta}} J(\bm{\theta})}{\sqrt{{\nabla_{\bm{\theta}} J(\bm{\theta})}^\intercal \mathcal{I}_{\bm{\theta}}^{-1} \nabla_{\bm{\theta}} J(\bm{\theta})}}
	\end{equation}
	constitutes a geodesic vector field on some neighbourhood of $\bm{\theta}^{t}$ in the Riemannian manifold $\left(S,\mathcal{I}_{\bm{\theta}}\right)$, with $S$ being a statistical manifold \cite{lauritzen1987} and with the Fisher information matrix $\mathcal{I}_{\bm{\theta}}$ used as the associated metric tensor. \cite{amari} refers to $\mathcal{I}_{\bm{\theta}}^{-1} \nabla_{\bm{\theta}} J(\bm{\theta})$ as the natural gradient, denoted by $\widetilde{\nabla} J(\bm{\theta})$.
	
	Kullback-Leibler divergence induces a local distance between $P_{\bm{\theta}}$ and $P_{\bm{\theta} +\delta\bm{\theta}}$ on $S$ when $\delta\bm{\theta}$ is sufficiently small, given by $D_{KL}\left(P_{\bm{\theta}} \ || \ P_{\bm{\theta} +\delta\bm{\theta}} \right)=\frac{1}{2}\delta\bm{\theta}^\intercal \mathcal{I}_{\bm{\theta}}\delta\bm{\theta} + \mathcal{O}\left({\parallel \delta\bm{\theta} \parallel}^2 \right)$, where $\mathcal{O}({\parallel \delta\bm{\theta} \parallel}^2 )$ vanishes at least as fast as ${\parallel \delta\bm{\theta} \parallel}^2$ as $\delta\bm{\theta}$ tends to zero. The square of this distance is denoted by ${ds}^2=\delta\bm{\theta}^\intercal \mathcal{I}_{\bm{\theta}}\delta\bm{\theta}$. Then, combining the natural gradient $\widetilde{\nabla}$ and the trust region radius constraints \cref{eq:trustcon} given by $ds\leq \Delta^{t}$, the update equations to solve \cref{eq:Jmax} take the following form
	\begin{equation}\label{eq:updateIGO}
		\bm{\theta}^{t+1}=\bm{\theta}^{t}+\Delta^{t}\frac{ \widetilde{\nabla} J(\bm{\theta}^{t})}{\parallel\widetilde{\nabla} J(\bm{\theta}^{t})\parallel}_{\mathcal{I}_{\bm{\theta}}},
	\end{equation}
	where ${\parallel\cdot\parallel}_{\mathcal{I}_{\bm{\theta}}}$ is the norm associated with the inner product induced by the metric tensor $\mathcal{I}_{\bm{\theta}}$.
	
	Model \cref{eq:trust} can be considered a probabilistic equivalent of a Euclidean trust region model \cite{necedal}. Furthermore, as the search direction of the Euclidean first-order trust region model coincides with the search direction of the line search method, \cref{eq:updateIGO} can be similarly viewed as a geodesic search on the statistical manifold $S$, where the search moves along a geodesic given by \cref{eq:unitGrad} with the step length given by $\Delta^{t}$.
	
	\subsection{The extended multivariate von Mises distribution}
	\label{sec:directionalstatistics}
	Crystal lattice \cref{eq:latticeG} induces a quotient space $\mathbb{R}^n/ \text{L}$, where $\text{L}$ is a lattice group
	\cref{eq:latticeG}. Since $\mathbb{R}^n/L_{G}$ is homeomorphic to the $n$D torus denoted by $T^n$, a natural choice for the statistical model $S$ of \cref{eq:Smodel} is to restrict the search to a family of probability distributions defined on $T^n$.
	
	\cite{mardia1975} defined a probability distribution with the support on a $2$D unit flat torus and, using the sine submodel of the general bivariate von Mises model, introduced a probability distribution on an $n$D torus \cite{mardia2008}.  Following the general bivariate von Mises model, we extended the multivariate von Mises model to the family of distributions with the probability density function
	\begin{equation}\label{eq:fmvm}
		f(\bm{\theta}|\bm{\mu}, \bm{\kappa} , \textbf{D}) = \frac{1}{ Z(\bm{\mu}, \bm{\kappa} , \textbf{D})} \exp \left \{ \bm{\kappa}^\intercal c(\bm{\theta}-\bm{\mu})+ \frac{1}{2} \left[\begin{matrix}
			c(\bm{\theta}-\bm{\mu}) \\ 
			s(\bm{\theta}-\bm{\mu})
		\end{matrix}\right]^\intercal\textbf{D} \left[\begin{matrix}
			c(\bm{\theta}-\bm{\mu}) \\ 
			s(\bm{\theta}-\bm{\mu})
		\end{matrix}\right] \right \}
	\end{equation}
	where
	\begin{gather*}
		c(\bm{\theta}-\bm{\mu})=\left[cos(\theta_1-\mu_1),\ldots,cos(\theta_n-\mu_n)\right]^\intercal,\\ s(\bm{\theta}-\bm{\mu})=\left[sin(\theta_1-\mu_1),\ldots,sin(\theta_n-\mu_n)\right]^\intercal, \\
		0\leq\theta_i, \mu_i\leq2\pi, \quad  0\leq\kappa_i,
	\end{gather*}
	and the $2n \times 2n$ real valued symmetric matrix that controls the cosine-sine interactions
	\begin{equation*}
		\textbf{D}=\left [\begin{matrix}
			\textbf{D}^{cc} & \textbf{D}^{cs} \\ 
			\left.\textbf{D}^{cs}\right.^\intercal & \textbf{D}^{ss}
		\end{matrix} \right] \\
	\end{equation*}
	has all the diagonal elements $d_{ii}^{cc}, \ d_{ii}^{ss} , \ d_{ii}^{cs}$ of its corresponding submatrices set to zero.
	
	The explicit form of the normalizer $Z\left(\bm{\mu}, \bm{\kappa} , \textbf{D}\right)$ is known only in a few instances of the bivariate case \cite{jupp1980, singh2002}. When an explicit evaluation is impossible, the standard approach to compute integrals is to use Monte--Carlo methods. For the exponential family statistical model used here, the Monte--Carlo estimates are discussed in \cref{subsec:montecarlo}.
	
	It has to be noted that the extended multivariate von Mises model \cref{eq:fmvm} is not identifiable and thus problematic to parameterize in the form of an exponential family of distributions since for $\bm{\kappa}=\textbf{0}$ and any fixed $\textbf{D}$, density functions \cref{eq:fmvm} for any $\bm{\mu}_1 \neq \bm{\mu}_2$ are equal due to the vanishing of the term containing $\bm{\kappa}$. Nevertheless, when the concentration parameters $\bm\kappa$ are restricted to being strictly positive, the full multivariate von Mises model \cref{eq:fmvm} can be rewritten using trigonometric identities in the following form
	\begin{equation}\label{eq:eMVM}
		f(\bm{\theta}|\bm{\eta}, \textbf{E})= \exp\left \{ \left[\begin{matrix}
			c(\bm{\theta}) \\ 
			s(\bm{\theta})
		\end{matrix}\right]^\intercal
		\bm{\eta}
		+ \text{vec}\left( \left[\begin{matrix}
			c(\bm{\theta}) \\ 
			s(\bm{\theta})
		\end{matrix}\right] \left[\begin{matrix}
			c(\bm{\theta}) \\ 
			s(\bm{\theta})
		\end{matrix}\right]^\intercal\right)^\intercal \text{vec} \left(  \textbf{E}  \right) - \psi\left(\bm{\eta} , \textbf{E} \right) \right \}
	\end{equation}
	where $\text{vec}(\cdot)$ denotes vectorization, $\psi\left(\cdot,\cdot \right)$ is the logarithm of the normalizing constant or $\log$-partition function, and the canonical exponential family parameters $(\bm{\eta},\textbf{E})$ are given as
	\begin{equation}  \label{eq:canonicalPar}
		\bm{\eta}=\left[\begin{matrix}
			\bm{\kappa}\odot c(\bm{\mu}) \\
			\bm{\kappa}\odot s(\bm{\mu})
		\end{matrix}\right], \quad
		\textbf{E}=\left [\begin{matrix}
			\textbf{E}^{cc} & \textbf{E}^{cs} \\ 
			\left.\textbf{E}^{cs}\right.^\intercal & \textbf{E}^{ss}
		\end{matrix} \right].
	\end{equation}
	where $\odot$ denotes Hadamard product. The submatrices can be expressed as follows
	\begin{subequations}
		\begin{equation*}
			\textbf{E}^{cc}=\tfrac{1}{2}\textbf{D}^{cc}\odot c(\bm{\mu})c(\bm{\mu})^\intercal -\tfrac{1}{2}\textbf{D}^{cs}\odot c(\bm{\mu})s(\bm{\mu})^\intercal + \tfrac{1}{2}\textbf{D}^{ss}\odot s(\bm{\mu})s(\bm{\mu})^\intercal - \tfrac{1}{2}\left.\textbf{D}^{cs}\right.^\intercal\odot s(\bm{\mu})c(\bm{\mu})^\intercal
		\end{equation*}
		\begin{equation*}
			\textbf{E}^{cs}=\tfrac{1}{2}\textbf{D}^{cs}\odot c(\bm{\mu})c(\bm{\mu})^\intercal+\tfrac{1}{2}\textbf{D}^{cc}\odot c(\bm{\mu})s(\bm{\mu})^\intercal - \tfrac{1}{2}\textbf{D}^{ss}\odot s(\bm{\mu})c(\bm{\mu})^\intercal - \tfrac{1}{2}\left.\textbf{D}^{cs}\right.^\intercal\odot s(\bm{\mu})s(\bm{\mu})^\intercal
		\end{equation*}
		\begin{equation*}
			\textbf{E}^{ss}=\tfrac{1}{2}\textbf{D}^{ss}\odot c(\bm{\mu})c(\bm{\mu})^\intercal+\tfrac{1}{2}\textbf{D}^{cc}\odot s(\bm{\mu})s(\bm{\mu})^\intercal +\tfrac{1}{2}\textbf{D}^{cs}\odot s(\bm{\mu})c(\bm{\mu})^\intercal + \tfrac{1}{2}\left.\textbf{D}^{cs}\right.^\intercal\odot c(\bm{\mu})s(\bm{\mu})^\intercal .
		\end{equation*}
	\end{subequations}
	
	The property that the reparameterisation of the extended multivariate von Mises model \cref{eq:eMVM} is an exponential family probability distribution grants us a few beneficial properties besides the ability to use Monte--Carlo methods to estimate the natural gradients of \cref{subsec:montecarlo}. For example, from the maximum entropy property of the exponential family, we have that \cref{eq:fmvm} is the maximum entropy or minimum discrimination distribution when $(\bm{\mu}, \bm{\kappa} , \textbf{D})$ parameters are specified. Moreover, exponential families provide the statistical model $S$ \cref{eq:Smodel} with a dually flat structure \cite{amari2000}, a concept that is utilized in  \Cref{sec:subgradient} and \Cref{sec:geometryEntTrust}.
	
	Additional details on the extended multivariate von Mises distribution, parameter transformations as well as implementation details for the Gibbs sampler used in the Monte--Carlo estimates are presented in \Cref{sec:toroidalDist} and subsections therein.
	
	\subsection{Adaptive selection quantile}
	\label{sec:adaptiveselectionquantile}
	
	Using the raw expected fitness \cref{eq:efit} results in poor performance due to the sensitivity of the sample mean estimates of the expected fitness to extreme values \cite{beyer}. The usual way to address this in evolutionary computation methods \cite{ollivier, wierstra} is to employ a rank-preserving transformation of the fitness function, referred to as the selection quantile \cite{beyer2001}. The main advantage of using quantiles instead of raw fitnesses is that quantiles are highly robust measures of position. This section proposes a selection quantile fitness transformation inspired by the simulated annealing control parameter \cite{laarhoven}.
	
	One particular selection quantile is of the form
	\begin{equation}\label{eq:qES}
		q \ \mathbf{ 1}_{\textbf{F}_{1- \frac{1}{q}}^{\tilde{\bm{\theta}}}}(\textbf{x}):=\left \{\begin{matrix}
			q & \ \text{if} \ \textbf{F}(\textbf{x}) \geq \textbf{F}_{1- \frac{1}{q}}^{\tilde{\bm{\theta}}},  \\
			0& \text{otherwise},
		\end{matrix} \right.
	\end{equation}
	where $\textbf{F}_{1- \frac{1}{q}}^{\tilde{\bm{\theta}}}$ is the $P(\tilde{\bm{\theta}})$  $ (q - 1)$-th $q$-quantile of the fitness \textbf{F}, using the indicator function $\mathbf{ 1}_{\bm\cdot}\!(\cdot)$.   Since the expected value of the indicator function \cref{eq:qES} is equal to the probability of observing $\textbf{F(x)}$ being greater than $\textbf{F}_{1- \frac{1}{q}}$ as
	\begin{equation}\label{eq:qESJtheta}
		J(\bm{\theta})=\int_{\mathcal{X}} q\mathbf{ 1}_{\textbf{F}_{1- \frac{1}{q}}^{\tilde{\bm{\theta}}}}(\textbf{x})dP(\bm{\theta})=q\int_{\mathcal{V}} dP(\bm{\theta})=qP\left(\textbf{F}(\textbf{X}) \geq \textbf{F}_{1- \frac{1}{q}}^{\tilde{\bm{\theta}}}|\bm{\theta}\right)
	\end{equation}
	where
	\begin{math}
		\mathcal{V}=\{\textbf{x} | \textbf{F}(\textbf{x}) \geq \textbf{F}_{1- \frac{1}{q}}^{\tilde{\bm{\theta}}}\}
	\end{math}
	if we set $\tilde{\bm{\theta}}=\bm{\theta}$, then we have $qP\left(\textbf{F}(\textbf{X}) \geq \textbf{F}_{1- \frac{1}{q}}^{\tilde{\bm{\theta}}}|\tilde{\bm{\theta}}\right)=1$, and 
	\begin{equation}\label{eq:truncDist}
		J(y  |  \tilde{\bm{\theta}})=\int_{\textbf{F}^{\tilde{\bm{\theta}}}_{1- \frac{1}{q}}}^{y}q dP\circ \textbf{F}^{-1}(z)(\tilde{\bm{\theta}})
	\end{equation}
	can be regarded as a probability distribution function of the random variable $Y=\textbf{F}(\textbf{X})$ truncated at $\textbf{F}_{1- \frac{1}{q}}^{\tilde{\bm{\theta}}}$ with the density function $q p(\textbf{F}^{-1}(y)|\tilde{\bm{\theta}})$ where $p(\textbf{x}|\tilde{\bm{\theta}})$ is the density function of the random variable $\textbf{X}$.
	
	Specifically, given a random vector $\mathbf{ X}$ with the exponential family density function
	\begin{equation}\label{eq:expDensity}
		p(\mathbf{x}|\tilde{\bm{\theta}})=\exp\left\{ \tilde{\bm{\theta}}^\intercal \textbf{t}\left(\textbf{x}\right) - \psi\left( \tilde{\bm{\theta}} \right) \right\},
	\end{equation}
	the density function of the random variable $Y=\textbf{F}(\textbf{X})$ truncated at $\textbf{F}_{1- \frac{1}{q}}^{\tilde{\bm{\theta}}}$ can be expressed as
	\begin{equation}\label{eq:truncExpDensity}
		p(y|\tilde{\bm{\theta}})=\exp\left\{ \tilde{\bm{\theta}}^\intercal \textbf{t}\left(\textbf{F}^{-1}(y)\right) - \psi \left( \tilde{\bm{\theta}} \right)+\ln\left(q\right) \right\}w(y)
	\end{equation}
	for some Borel measurable function $w$ on $\mathcal{Y}$.
	
	The $P(\tilde{\bm{\theta}})$  $ (q - 1)$-th $q$-quantile of the fitness \textbf{F} is generally unknown and has to be estimated. Given $N$ observations of $\textbf{x}_i \sim P (\tilde{\bm{\theta}})$, the empirical quantile function is constructed by assigning the ranks
	\begin{equation*}
		r_i^{\tilde{\bm{\theta}}}:=\left\{ i \ | \ \textbf{F}(\textbf{x}_1)\leq\textbf{F}(\textbf{x}_2)\leq\ldots\leq\textbf{F}(\textbf{x}_i)\leq\ldots\leq\textbf{F}(\textbf{x}_{N-1})\leq\textbf{F}(\textbf{x}_N) \ | \ \textbf{x}_i \sim P (\tilde{\bm{\theta}})  \right\}
	\end{equation*}
	to each fitness value. Then the $\textbf{F}_{1- \frac{1}{q}}^{\tilde{\bm{\theta}}}$ estimate $\hat{\textbf{F}}_{1- \frac{1}{q}}^{\tilde{\bm{\theta}}}$ is given by $\hat{\textbf{F}}_{1- \frac{1}{q}}^{\tilde{\bm{\theta}}}=\textbf{F} (\textbf{x}_{r_{\left\lfloor N- \frac{N}{q}\right\rfloor}^{\tilde{\bm{\theta}}}})$, where $\lfloor \cdot \rfloor$ denotes rounding to the nearest lower integer.
	
	This process can be repeated multiple times to estimate $\textbf{F}_{1- \frac{1}{q}}^{\tilde{\bm{\theta}}}$ by taking the maximum of the empirical quantile $\textbf{F} (\textbf{x}_{r_{\left\lfloor N- \frac{N}{q}\right\rfloor}^{\tilde{\bm{\theta}}}})$ from all batch sampling iterations. The estimate then takes the form
	\begin{equation*}
		\hat{\textbf{F}}_{1- \frac{1}{q}}^{\tilde{\bm{\theta}}}=\max_{i=1\ldots}\left[ \textbf{F} (\textbf{x}^{i}_{r_{\left\lfloor N- \frac{N}{q}\right\rfloor}^{\tilde{\bm{\theta}}}}) \right].
	\end{equation*}
	
	In practice, with the above, we simulate realizations from a truncated probability distribution of the distribution in \cref{eq:truncDist}, similarly to the simulated annealing method where the homogeneous algorithm \cite{laarhoven} is a Metropolis-Hastings one for generating realizations from the Boltzmann distribution. Continuing with this analogy, the $q$-quantile can be considered equivalent to the temperature control parameter.
	
	At the beginning of the search, it is beneficial to have a smaller $q$-quantile to maintain a more extensive profile of the overall optimization landscape. As the algorithm progresses and we need the distribution to concentrate on samples with higher fitness values, higher values of $q$ and a more localized search is preferred. We implement this using a time varying $q_t$-quantile by setting
	\begin{equation}\label{eq:cDet}
		q_{t+1}=q_{t}\exp\left\{\beta t\right\}.
	\end{equation}
	Consequently, for fixed $\tilde{\bm{\theta}}$, when $t$ approaches infinity, the probability distribution \cref{eq:truncDist} converges to a distribution where all the probability mass is concentrated on the extrema of the fitness $F$.
	
	Following the fitness transformation \cref{eq:qES}, with $P(\bm{\theta})$ from the exponential family of distributions with the $q$-quantile fixed, the gradient of \cref{eq:qESJtheta} at $\tilde{\bm{\theta}}$ takes the form
	\begin{equation}\label{eq:adaThetaGrad}
		\nabla_{\bm{\theta}} J(\bm{\theta})\mid_{\bm{\theta}=\tilde{\bm{\theta}}} \ = \ \bm{\mu}_{\textbf{F}_{1- \frac{1}{q_{t}}}}-\bm{\mu},
	\end{equation}
	where $\bm{\mu}_{\textbf{F}_{1- \frac{1}{q_{t}}}}$ is the expectation parametrization of the truncated exponential probability distribution \cref{eq:truncExpDensity}, given by
	\begin{equation}\label{eq:qTrunctE}
		\bm{\mu}_{\textbf{F}_{1- \frac{1}{q_{t}}}} = \int\limits_{\textbf{F}^{\tilde{\bm{\theta}}}_{1- \frac{1}{q_{t}}}}^{\infty} \textbf{F}^{-1}(y) \exp\left\{ \tilde{\bm{\theta}}^\intercal \textbf{t}\left(\textbf{F}^{-1}(y)\right) - \psi( \tilde{\bm{\theta}})+\ln\left(q_{t}\right) \right\} d\textbf{F}^{-1}(y),
	\end{equation}
	and $\bm{\mu}$ is the expectation parametrization of the exponential distribution with density function \cref{eq:expDensity}, equal to
	\begin{math}
		\bm{\mu} = \int \textbf{x} \exp\{ \tilde{\bm{\theta}}^\intercal \textbf{t}\left(\textbf{x}\right) - \psi( \tilde{\bm{\theta}}) \}dx.
	\end{math}
	Finally, the estimate of \cref{eq:adaThetaGrad} is given by
	\begin{equation}\label{eq:adaSelEst}
		\nabla_{\bm{\theta}} \hat{J}(\bm{\theta})\mid_{\bm{\theta}=\tilde{\bm{\theta}}}=\hat{E}\left[ \left. \textbf{t}\left(\textbf{x}\right) \right| \textbf{F}(\textbf{x}) \geq \hat{\textbf{F}}_{1- \frac{1}{q_{t}}}^{\tilde{\bm{\theta}}} \right]-\hat{E}\left[\textbf{t}\left(\textbf{x}\right)\right],
	\end{equation}
	with $\hat{E}\left[\cdot\right]$ denoting the sample mean.
	
	A caveat of simulated annealing with a fixed annealing schedule is the risk of the algorithm being trapped in a local maximum if the control parameter converges too fast. Therefore, the scheduling constant $\beta$ in \cref{eq:cDet} has to be set optimally, which is usually done experimentally to counter this issue.
	
	Another way to mitigate premature convergence caused by the control parameter \cref{eq:cDet} is by introducing self-adaptation into the control schedule where the selection quantile \cref{eq:qES} is used in conjunction with the trust region method \cref{eq:updateIGO}. For a fixed $q_t$, the algorithm performs a local search with the neighbourhood given by the $q_t$-quantile. By increasing $q_t$ the neighbourhood becomes more localized, resulting in a decreased chance of escaping the attraction of local maxima. This behaviour is desirable in later stages when optima basins have already been selected. On the other hand, in instances when the trust region path reverses direction frequently, indicating a tendency toward moving to multiple local optima, it is beneficial to increase the search neighbourhood by decreasing $q_t$ and provide the algorithm with less localized information for inference of various optima basins of the optimization landscape. 
	
	To assess the trajectory's current state, we compare the directions of three consecutive updates
	$\bm{\theta}^{t-2}, \bm{\theta}^{t-1},  \bm{\theta}^{t}$ by expressing their difference
	\begin{gather*}
		\Delta \bm{\theta}^{t} = \bm{\theta}^{t} - \bm{\theta}^{t-1}, \\
		\Delta \bm{\theta}^{t-1} = \bm{\theta}^{t-1} - \bm{\theta}^{t-2},
	\end{gather*}
	as elements of the tangent space $T_{\bm{\theta}^{t-1}}S$ of a statistical manifold $S$ at point $\bm{\theta}^{t-1}$. The angle $\alpha^t$ between vectors $\Delta \bm{\theta}^{t}$ and $\Delta \bm{\theta}^{t-1}$ induced by the scalar product $<\cdot,\cdot>_{\mathcal{I}_{\bm{\theta}^{t-1}}}$ with the Fisher metric tensor $\mathcal{I}_{\bm{\theta}^{t-1}}$ is given by
	\begin{equation}\label{eq:diffangle}
		\cos(\alpha^t)=\frac{<\Delta\bm{\theta}^t, \Delta\bm{\theta^}{t-1}>_{\mathcal{I}_{\bm{\theta}^{t-1}}}}{||\Delta\bm{\theta}^t||_{\mathcal{I}_{\bm{\theta}^{t-1}}}||\Delta\bm{\theta}^{t-1}||_{ \mathcal{I}_{\bm{\theta}^{t-1}}}}.
	\end{equation}
	By combining the quantile control parameter $q_t$ \cref{eq:cDet} and the direction of the current state of algorithm \cref{eq:diffangle}, we can use an adaptive selection quantile scheme as
	\begin{equation}\label{eq:adaSQ}
		q_{t+1}=q_{t}\exp\left\{\beta \cos(\alpha^t)\right\}=q_{0}\exp\left\{\beta \sum_{i=3}^{t} \cos(\alpha^i)\right\},
	\end{equation}
	where $q_{0}$ is the initial value of the quantile control parameter. 
	
	The expression \cref{eq:adaThetaGrad} has a clear geometric interpretation. $\nabla_{\bm{\theta}} J(\bm{\theta})\mid_{\bm{\theta} = \tilde{\bm{\theta}}}$ can be regarded as a vector in the direction of $\bm{\mu}_{\textbf{F}_{1- \frac{1}{q_{t}}}}$, that is, in the direction given by the $ (q_{t} \!-\! 1) $th $q_{t}$-quantile of the fitness. When the trust region steps \cref{eq:updateIGO} follow the same general direction, the adaptive selection quantile \cref{eq:adaSQ} increases, forcing the trust region to move toward higher fitness values. On the other hand, if $\cos\left(\alpha^t\right)<0$, the algorithm backtracks and $q_t$ is decreased, resulting in a broader search. In terms of evolutionary computation, the adaptive selection quantile models a type of time-varying selective pressure.
	
	The proposed selection quantile can be used independently from the entropic trust region method, and its implementation details are provided in the form of pseudocode in \ref{sec:quntileHillClimb}. The value of parameter $\beta$ used in our experiments is discussed in \ref{sec:algparameters}.
	
	\subsection{Entropic proximal maximization for exponential families}
	\label{sec:subgradient} 
	
	Proximal mappings are used in convex optimization to construct approximations of the objective function that have a smoothing effect and preserve the set of minimizers \cite{rockafellar2009}.
	\cite{beck} showed equivalence between the mirror descent algorithm \cite{nemirovskij} and a modified proximal minimization algorithm by taking into account only first-order approximation of the objective function and considering a more general class of proximal operators \cite{censor1992}. 
	Subsequently, \cite{raskutti} showed equivalence between the natural gradient descent method \cite{amari} and mirror descent. Following these results, we show that the entropic trust region method \cref{eq:updateIGO} for exponential family distributions is a special case of a proximal maximization where the projection term in the proximal operator is given by the KLD \cref{eq:KLdiv}, resulting in a gradient ascent in dual coordinates, maximizing the likelihood and minimizing cross entropy of parameter update estimates.
	
	Taking the first-order Taylor expansion of $J\left(\bm{\theta}\right)$ in the neighbourhood of $\bm{\theta}^{t}$ with respect to the exponential family statistical model $S$ with canonical parameterization $\bm{\theta}\in \bm{\Theta}$ and expectation parameterization $\bm{\mu}\in \mathcal{M}$, the proximal method for problem \eqref{eq:Jmax} can be equivalently rewritten to
	\begin{equation}\label{eq:thetaProx}
		\bm{\theta}^{t+1}=\argmax_{\bm{\theta}\in \bm{\Theta}}\left\{ \bm{\theta}^\intercal\nabla_{\bm{\theta}} J(\bm{\theta}^{t}) - \frac{1}{\epsilon} D_{KL}\left( P_{\bm{\theta}^t} \ || \  P_{\bm{\theta}} \right) \right\},
	\end{equation}
	with parameter $\epsilon>0$  where the projection of $\bm{\theta}^{t}$ onto $\bm{\Theta}$ is given by the KLD. From the optimality condition of \cref{eq:thetaProx}, we directly get the mirror ascent in the form
	\begin{subequations} \label{eq:mirrorAscMu}
		\begin{gather}
			\bm{\mu}^{t+1}= \bm{\mu}^{t}  + \epsilon\nabla_{\bm{\theta}} J(\bm{\theta}^{t}) \label{eq:mirrorCondMu}, \\
			\bm{\theta}^{t+1}=\nabla_{\bm{\mu}} \phi ( \bm{\mu}^{t+1}) \label{eq:mirrorCondTheta},
		\end{gather}
	\end{subequations}
	where $\bm{\theta}$ and $\bm{\mu}$ are the dual canonical, and exponential parametrizations of the exponential family given by the Legendre transforms $\bm{\theta}=\nabla_{\bm{\mu}} \phi ( \bm{\mu})$ and $\bm{\mu}=\nabla_{\bm{\theta}} \psi ( \bm{\theta})$
	of the dual free energy functionals $\psi ( \bm{\theta})$ and $\phi ( \bm{\mu})$. 
	
	Differentiating $\bm{\mu}$ as a function of $\bm{\theta}$ with respect to time yields the time derivative
	\begin{equation*}
		\frac{d\bm{\mu}\left(\bm{\theta}\right)}{dt}=\nabla^2_{\bm{\theta}}\psi\left(\bm{\theta}\right)\frac{d\bm{\theta}}{dt}.
	\end{equation*}
	For a sufficiently small $dt$ and from the well-known fact that the Hessian of the free energy $\psi\left(\theta\right)$ is equal to the Fisher information matrix $\mathcal{I}_{\bm{\theta}}$ in $\bm{\theta}$ parameterization \cite{nielsen}, the relationship between update step sizes in $\bm{\theta}$ and $\bm{\mu}$ parametrizations is given by
	\begin{equation}\label{eq:tanTrans}
		\bm{\mu}^{t+1} - \bm{\mu}^{t} = \mathcal{I}_{\bm{\theta}^{t}}\left(\bm{\theta}^{t+1} - \bm{\theta}^{t}\right).
	\end{equation}
	Using \cref{eq:tanTrans}, the mirror ascent update \cref{eq:mirrorCondMu} can be restated in the form of
	\begin{equation}\label{eq:thetaIGO}
		\bm{\theta}^{t+1}= \bm{\theta}^{t} + \epsilon\mathcal{I}_{\bm{\theta}^{t}}^{-1}\nabla_{\bm{\theta}} J(\bm{\theta}^{t}),
	\end{equation}
	which is the natural gradient method \cite{amari} or information geometric optimization algorithm \cite{ollivier}. Furthermore, by setting 
	\begin{equation}
		\label{eq:epsylon}
		\epsilon=\frac{\Delta^{t}}{\sqrt{{\nabla_{\bm{\theta}} J(\bm{\theta}^{t})}^\intercal\left(\mathcal{I}_{\bm{\theta}^{t}}\right)^{-1} \nabla_{\bm{\theta}} J(\bm{\theta}^{t})}},
	\end{equation}
	we get the entropic trust region method of \cref{eq:updateIGO}.
	
	Statistically, the update $\bm{\mu}^{t+1}$ in \cref{eq:mirrorCondMu} is unknown, and we are working with the point estimate
	\begin{equation}
		\label{eq:updateEstimate}
		\hat{\bm{\mu}}^{t+1}=\frac{1}{N}\sum_{i=1}^{N}\bm{\mu}^{t} +\epsilon\textbf{F}(\textbf{X}_i) (\textbf{X}_i-\bm{\mu}^{t}),
	\end{equation}
	where $\textbf{X}_i$ are independent identical exponential distributed family random vectors. This means that the update estimate $\hat{\bm{\mu}}^{t+1}$ is by itself a random vector from some parametric family of probability distributions, although generally not exponential family distributed. Nevertheless, using \cref {eq:updateEstimate}, the proximal map \cref{eq:thetaProx} can be rewritten as
	\begin{equation*}\label{eq:thetaProx1}
		\bm{\theta}^{t+1}=\argmax_{\bm{\theta}\in \bm{\Theta}}\left\{ \bm{\theta}^\intercal\hat{\bm{\mu}}^{t+1} - \psi ( \bm{\theta}) \right\},
	\end{equation*}
	where $\ln p(\hat{\bm{\mu}}^{t+1},\bm{\theta})= \bm{\theta}^\intercal\hat{\bm{\mu}}^{t+1} - \psi ( \bm{\theta})$ is the log-likelihood of the exponential family statistical model $S$ given estimate \cref{eq:updateEstimate}. Thus, the mirror map projection step \cref{eq:mirrorCondTheta} is given by the maximum likelihood estimate of $\bm{\theta}^{t+1}$.
	
	Due to the duality of exponential families, we can rewrite the mirror ascent and thus the proximal map  \cref{eq:thetaProx} in the expectation parameters providing us with dual trust region updates.   Expressing the expected fitness \cref{eq:efit} in terms of expectation parameters $\bm{\mu}$, yields 
	\begin{multline*}
		J(\bm{\theta})=\int_{\textbf{X}} \textbf{F(x)}\exp\left\{ \bm{\theta}^\intercal \textbf{t}\left(\textbf{x}\right) - \psi\left( \bm{\theta} \right) \right\}dP_0\left(\textbf{x}\right)= \\
		=\int_{\textbf{X}} \textbf{F(x)}\exp\left\{ \nabla_{\bm{\mu}} \phi\left(\bm{\mu}\right)^\intercal\left(\textbf{t}\left(\textbf{x}\right)-\bm{\mu}\right) + \phi\left( \bm{\mu} \right) \right\}dP_0\left(\textbf{x}\right)=J(\bm{\mu}),
	\end{multline*}
	which, when differentiated with respect to $\bm{\mu}$, gives the $\bm{\theta}$ gradient of the expected fitness $J$ as
	\begin{equation}\label{eq:muGradExp}
		\nabla_{\bm{\mu}} J( \bm{\mu}) =  \nabla_{\bm{\mu}}^{2} \phi ( \bm{\mu}) \nabla_{\bm{\theta}} J(\bm{\theta}).
	\end{equation}
	
	Using \cref{eq:muGradExp} and the relationship between the Fisher information matrices in canonical and dual parametrizations and the Hessian of the Legendre dual $\phi$ given by
	\begin{equation*}
		\mathcal{I}_{\bm{\theta}}^{-1}=\mathcal{I}_{\bm{\mu}}=\nabla_{\bm{\mu}}^{2} \phi ( \bm{\mu}),
	\end{equation*}
	produces the mirror ascent in the expectation parameters, which can be written as
	\begin{subequations}\label{eq:mirrorAscTheta}
		\begin{gather}
			\bm{\theta}^{t+1}= \bm{\theta}^{t} + \epsilon \nabla_{\bm{\mu}} J(\bm{\mu}^{t}) \label{eq:mirrorIterTheta}, \\
			\bm{\mu}^{t+1}= \nabla_{\bm{\theta}} \psi ( \bm{\theta}^{t+1}) \label{eq:mirrorIterMu},
		\end{gather}
	\end{subequations}
	or equivalently re-expressed as the proximal map maximization
	\begin{equation}\label{eq:muProximal}
		\bm{\mu}^{t+1}=\argmax_{\bm{\mu}\in \mathcal{M}}\left\{ \bm{\mu}^\intercal\nabla_{\bm{\mu}} J(\bm{\mu}^{t}) - \frac{1}{\epsilon}D_{KL}\left( P_{\bm{\mu}} \ || \  P_{\bm{\mu}^{t}}  \right) \right\}.
	\end{equation}
	
	The preceding paragraphs make it clear that when the underlying statistical model is of the exponential family, the entropic trust region in \cref{eq:updateIGO} forms trajectories in the parameter spaces $\Theta$ and $\mathcal{M}$ related by the Bregman divergence derived from the free energy functional $\psi ( \bm{\theta})$.
	In fact, these trajectories are given by a gradient flow of the expected fitness in dual coordinates, which will be further explored in the following section.		
	
	As in the case of the proximal map expressed in the canonical parametrization \cref{eq:thetaProx}, the estimate of mirror ascent iteration \cref{eq:mirrorIterTheta} defines the random vector
	\begin{equation}\label{eq:thetaEstimate}
		\hat{\bm{\theta}}^{t+1}=
		\frac{1}{N}\sum_{i=1}^{N} \bm{\theta}^{t} + \epsilon {\nabla_{\bm{\theta}}^{2} \psi (\bm{\theta}^t)}^{-1} \textbf{F}(\textbf{X}_i) \left(\textbf{X}_i-\nabla_{\bm{\theta}} \psi (\bm{\theta}^t)\right).
	\end{equation}
	Then the proximal maximization \cref{eq:muProximal} can be rewritten as
	\begin{equation*}
		\bm{\mu}^{t+1}=\argmax_{\bm{\mu}\in \mathcal{M}}\left\{ \bm{\mu}^\intercal\hat{\bm{\theta}}^{t+1} - \phi ( \bm{\mu}) \right\},
	\end{equation*}	
	and consequently as
	\begin{equation*}
		\bm{\mu}^{t+1}=\argmin_{\bm{\mu} \in \mathcal{M}} D_{KL}\left( P_{\nabla_{\bm{\mu}} \psi ( \bm{\mu})} \ || \ Q_{\hat{\bm{\theta}}^{t+1}}  \right),
	\end{equation*}
	assuming that $\hat{\bm{\theta}}^{t+1}$ belongs to a domain of the log-partition function $\psi(\cdot)$ and using the convex conjugate $\psi(\bm{\theta})=\bm{\theta}^{\intercal}\bm{\mu}-\phi(\bm{\mu})$, where $Q_{\hat{\bm{\theta}}^{t+1}}$ is the probability distribution of \cref{eq:thetaEstimate}. 
	Thus, the mirror ascent projection \cref{eq:mirrorIterMu} is the minimum cross entropy distribution from $S$, given the prior distribution of \cref{eq:thetaEstimate} \cite{shore}.
	
	\subsection{Geometry of the exponential family quantile-based entropic trust region method}
	\label{sec:geometryEntTrust}
	\Cref{sec:subgradient} introduced a relationship between the entropic trust region method and entropic proximal mappings for the particular case when the search is performed on an exponential family parametric statistical model. 
	When the quantile re-expression of the fitness of \cref{sec:adaptiveselectionquantile} is added, the dual gradient flows associated with the mirror ascents \cref{eq:mirrorAscMu} or alternatively \cref{eq:mirrorAscTheta} result in a geometrical interpretation relating evolutionary strategies, simulated annealing method and recurrent neural computing as instances of more general graphical interaction models \cite{cowell2007}. This section presents the exponential family quantile-based entropic trust region as an iterative solution to a minimax problem of the KLD between two hypersurfaces of an ambient statistical manifold.
	
	We consider a general exponential family
	\begin{equation}\label{eq:Se}
		S^{e}:=\left\{dP^{e}(\bm{\theta})=\exp\{ \bm{\theta}^\intercal \textbf{t}(\textbf{x}) - \psi( \bm{\theta} ) \}dP_0\left(\textbf{x}\right) \right\},
	\end{equation}
	where $dP^{e}$ is absolutely continuous with respect to some reference measure $dP_0\left(\textbf{x}\right)$ and $\psi\left( \bm{\theta} \right)=\ln \int \exp\{ \bm{\theta}^\intercal \textbf{t}\left(\textbf{x}\right) \}dP_0\left(\textbf{x}\right)$. Given the random variable $Y=\textbf{F}\left(\textbf{X}\right)$, where $\textbf{F}$ is as in \cref{sec:entTrust}, $S^{e}$ can be rewritten to
	\begin{equation}\label{eq:expFamY}
		S^{e}=\left\{\left. dP^{e}(\bm{\theta}) \right| \bm{\theta} \in \Theta^{e} \subset \mathbb{R}^{n} \right\},
	\end{equation}
	where $dP^{e}(\bm{\theta})=\exp\left\{ \bm{\theta}^\intercal \textbf{t}\left(\textbf{F}^{-1}(y)\right) - \psi \left( \bm{\theta} \right) \right\}dP_0 \circ \textbf{F}^{-1}(y)$ and $P_0\circ \textbf{F}^{-1}$	denotes the image of measure $P_0$ under the mapping $\textbf{F}$, and a family of probability distributions derived from \cref{eq:expFamY},  introduced in \cref{sec:adaptiveselectionquantile}, by truncating $S^{e}$  at the $P(\tilde{\bm{\theta}})$  $ (q\!-\!1) $th $q$-quantile of fitness \textbf{F}, as
	\begin{equation}\label{eq:expTrunc}
		S^{q}=\left\{\left.dP^{q}(\bm{\theta}) \right| \bm{\theta} \in \Theta^{q} \subset \mathbb{R}^{n} \right\},
	\end{equation}
	where 
	\begin{equation*}
		dP^{q}=\left\{\begin{array}{cc}	
			\exp\left\{ \bm{\theta}^\intercal \textbf{t}\left(\textbf{F}^{-1}(y)\right) - \psi \left( \bm{\theta}\right) + \ln(q) \right\}dP_0 \circ \textbf{F}^{-1}(y) & y \geq \textbf{F}_{1- \frac{1}{q}}^{\bm{\theta} }
			\\
			0 & \text{otherwise}
		\end{array}	\right.
	\end{equation*}
	and
	$\frac{1}{q}= \int_{\textbf{F}_{1- \frac{1}{q}}^{\bm{\theta}}} ^{\infty} \exp\{ \bm{\theta}^\intercal \textbf{t}\left(\textbf{F}^{-1}(y)\right) -\psi \left( \bm{\theta}\right) \}dP_0 \circ \textbf{F}^{-1}(y)$ with truncation parameter $\gamma=\textbf{F}_{1- \frac{1}{q}}^{\bm{\theta}}$.

	For a fixed truncation parameter $\gamma$, the parametric family of probability distributions $S^{q}$ is that of a regular exponential family derived from the $S^{e}$. On the other hand, if we consider $q \in [1; \infty)$ as another parameter, we can think of $S^{q}$ as the leaves of a foliation \cite{lee2013} of some statistical manifold and set a new family of probability distributions as
	\begin{equation}\label{eq:Sfam}
		S= \bigcup_{q \in [1; \infty)} S^{q}
	\end{equation}
	with the parametric space $\Theta \times [1; \infty) \subset \mathbb{R}^{n+1}$.	
	Since the truncation parameter in \cref{eq:expTrunc} is exactly determined by $\bm{\theta}$ and $q$, $S$ retains the dual structure given by the Bregman divergence \cite{bregman1967} induced by the convex function
	\begin{equation}\label{eq:freeEnergyTrunc}
		\psi \left( \bm{\theta}^q\right)=\psi \left( \bm{\theta}\right) - \ln(q),
	\end{equation}
	with the parametrization $\bm{\theta}^q=\left[ \bm{\theta},q\right]$. Moreover, since there is a bijection between regular exponential families and regular Bregman divergences \cite{banerjee2005}, $S$ is a regular exponential family. Consequently, the dual expectation parametrization is given by $\bm{\mu}^q=\left[\nabla_{\bm{\theta}}\psi \left( \bm{\theta}^q\right),-\frac{1}{q}\right]$ and the Fisher metric tensor by $\mathcal{I}_{\bm{\theta}^q} =  \left [\begin{matrix}
		\nabla^{2}_{\bm{\theta}}\psi ( \bm{\theta}^q) & \bm{0} \\
		\bm{0} & q^{-2}
	\end{matrix} \right]$.
	
	In the following, we show that the entropic proximal maximization algorithm of \cref{sec:subgradient} with the selection quantile fitness transformation \cref{eq:qES} constitutes dual gradient flows in $S$, maximizing the KLD from $S^{q^*}$ \cref{eq:expTrunc} to $S^1$ \cref{eq:expFamY} and minimizing this divergence from $S^1$ \cref{eq:expTrunc} to $S^{q^*}$ \cref{eq:expFamY}. In mathematical terms, we wish to solve the following max-min problem,
	\begin{equation}\label{eq:manifoldMin}
		\max_{P_{\bm{\theta}^{q^*}} \in S^{q^*}} \min_{P_{\bm{\theta}^1} \in S^1} D_{KL}(P_{\bm{\theta}^{q^*}} \ || \ P_{\bm{\theta}^1} ).
	\end{equation}
	
	By fixing $P_{\bm{\theta}^{q^*}}$, we can minimize $D_{KL}\left(P_{\bm{\theta}^{q^*}} \ || \ P_{\bm{\theta}^1} \right)$ with respect to $P_{\bm{\theta}^1}$. Since the KLD from $P_{\bm{\theta}^{q^*}}$ to $P_{\bm{\theta}^e}$ equals the Bregman divergence derived from the convex function \cref{eq:freeEnergyTrunc}, we have
	\begin{multline}\label{eq:bregTheta}
		D_{KL}(P_{\bm{\theta}^{q^*}} \ || \ P_{\bm{\theta}^1} )=D_{\psi}(\bm{\theta}^1,\bm{\theta}^{q^*})= \\ = \psi ( \bm{\theta}^q)\left|\right. _{\bm{\theta}^1} - \psi \left( \bm{\theta}^q\right)\left|\right. _{\bm{\theta}^{q^*}} -\nabla_{\bm{\theta}^q}\psi \left( \bm{\theta}^q\right)^\intercal \left|\right. _{\bm{\theta}^{q^*}}(\bm{\theta}^1-\bm{\theta}^{q^*}).
	\end{multline}
	Differentiating with respect to $\bm{\theta}^1$, the Riemannian gradient of \cref{eq:bregTheta} takes the form
	\begin{equation}\label{eq:gradMin}
		\text{grad}D_{\psi}(\bm{\theta}^1,\bm{\theta}^{q^*})= \mathcal{I}_{\bm{\theta}^1}^{-1}\left(\nabla_{\bm{\theta}^q}\psi \left( \bm{\theta}^q\right)\left|\right. _{\bm{\theta}^1} - \nabla_{\bm{\theta}^q}\psi \left( \bm{\theta}^q\right)\left|\right. _{\bm{\theta}^{q^*}} \right)= \\
		=\mathcal{I}_{\bm{\theta}^1}^{-1}( \bm{\mu}^1 - \bm{\mu}^{q^*}  ),
	\end{equation}
	and the gradient descent update equations to solve the minimization part of \cref{eq:manifoldMin} then takes the following form,
	\begin{equation}\label{eq:geoGradTheta}
		{\bm{\theta}^1}^{t+1}={\bm{\theta}^1}^{t}+ \epsilon\mathcal{I}_{{\bm{\theta}^1}^t}^{-1}( {\bm{\mu}^{q^*}}^t - {\bm{\mu}^1}^t  ).
	\end{equation}
	
	On the other hand, for a fixed $P_{\bm{\theta}^1}$, we can minimize $D_{KL}\left(P_{\bm{\theta}^{q}} \ || \ P_{\bm{\theta}^1}\right)$ with respect to $\bm{\theta}^q$. First, we express the KLD from $P_{\bm{\theta}^{q^*}}$ to $P_{\bm{\theta}^1}$ in terms of the Bregman divergence derived from the Legendre dual of $\psi \left( \bm{\theta}^q\right)$, given by
	\begin{math}
		\phi \left( \bm{\mu}^q\right)={\bm{\theta}^q}^\intercal{\bm{\mu}^q} -\psi \left( \bm{\theta}^q \right),
	\end{math}
	to the following form
	\begin{multline}\label{eq:bregMu}
		D_{KL}(P_{\bm{\theta}^{q^*}} \ || \ P_{\bm{\theta}^1})=D_{\phi}(\bm{\mu}^{q^*},\bm{\mu}^1)= \\
		=\phi \left( \bm{\mu}^q\right)\left|\right. _{\bm{\mu}^{q^*}} - \phi \left( \bm{\mu}^q\right)\left|\right. _{\bm{\mu}^1} -\nabla_{\bm{\mu}^q}\phi \left( \bm{\mu}^q\right)^\intercal \left|\right. _{\bm{\mu}^{q^*}}(\bm{\mu}^{q^*}-\bm{\mu}^1)
	\end{multline}
	and then by differentiating \cref{eq:bregMu} with respect to $\bm{\mu}^{q^*}$, the Riemannian gradient is expressed as 
	\begin{equation}\label{eq:gradMax}
		\text{grad} D_{\phi}(\bm{\mu}^{q^*},\bm{\mu}^1)= 
		\mathcal{I}_{\bm{\mu}^{q^*}}^{-1}   (\nabla_{\bm{\mu}^q}\phi \left( \bm{\mu}^q\right)\left|\right. _{\bm{\mu}^{q^*}} - \nabla_{\bm{\mu}^q}\phi \left( \bm{\mu}^q\right)\left|\right. _{\bm{\mu}^1} )
		=\mathcal{I}_{\bm{\mu}^{q^*}}^{-1}( \bm{\theta}^{q^*} - \bm{\theta}^{1}  ).
	\end{equation}
	Then, the gradient ascent update equations to solve the maximization part of \cref{eq:manifoldMin} are
	\begin{equation}\label{eq:geoGradMu}
		{\bm{\mu}^{q^*}}^{t+1}=
		{\bm{\mu}^{q^*}}^{t}+ \epsilon\mathcal{I}_{\bm{\mu}^{{q^*}^t}}^{-1}( {\bm{\theta}^{q^*}}^t - {\bm{\theta}^{1}}^t  ).
	\end{equation}
	
	For a fixed time $t$, the trust region search is being performed on some neighbourhood of $P_{\bm{\theta}^t}$ independently from $q^t$ since $q^t$ is updated between trust region updates by \cref{eq:adaSQ}. 
	Considering $S^{{q^*}}$ and $S^{{1}}$ as embeddings of $S^{e}$ \cref{eq:Se} in $S$ \cref{eq:Sfam}, we recover $S^{e}$ from $S$ via the projection $[\theta_{i}^{1}]_{i=1,\ldots,n} \rightarrow \bm{\theta}$ in the $\bm{\theta}$-coordinates for $\bm{\theta}^1 \in S^{{1}}$. For $\bm{\mu}^{q^{*}} \in S^{{q^*}}$, the $\bm{\mu}$-coordinates of $S^{e}$ are such that the following expression holds,
	\begin{math}
		\bm{\mu} = \frac{1}{{q^{*}}}[\mu^{q^{*}}_{i}]_{i=1,\ldots,n} + \frac{q^{*}-1}{{q^{*}}}\bm{\mu}_{\textbf{F}^C_{1-\frac{1}{q^{*}}}},
	\end{math}
	where 
	\begin{equation*}
		\bm{\mu}_{\textbf{F}^C_{1-\frac{1}{q^{*}}}} = \int\limits_{-\infty}^{\textbf{F}_{1- \frac{1}{q^{*}}}} \textbf{F}^{-1}(y) \exp\left\{ \bm{\theta}^\intercal \textbf{t}\left(\textbf{F}^{-1}(y)\right) - \psi( \bm{\theta})+\ln\left(\frac{q^{*}}{q^{*}-1}\right) \right\} dP_0 \circ \textbf{F}^{-1}(y)
	\end{equation*}	
	is the expected value of the complementary distribution to the truncated exponential distribution parameterized by $\bm{\theta}^{q^{*}}$ in $\cref{eq:expTrunc}$. From the dual gradient updates \cref{eq:geoGradTheta} and \cref{eq:geoGradMu} we then get the proximal maximization map in dual coordinates \cref{eq:thetaProx} and \cref{eq:muProximal} for the expected fitness $J(\bm{\theta})$ given by the selection quantile \cref{eq:qES}, and consequently the $\bm{\theta}$-coordinate trust region update \cref{eq:updateIGO} by setting $\epsilon$ as in \cref{eq:epsylon}.
	
	\cite{Ay2002} analyzed the infomax principle proposed in \cite{Linsker1997} as a universal rule to train artificial neural networks based on observations in biological neural networks. The infomax principle states that network weights should be set in such a way that the mutual information between input and output is maximized. Later, \cite{Ay2006} extended this work to maximizers of multi-information, a generalization of mutual information, previously defined in \cite{Watanabe} as total correlation and later proposed by \cite{Studeny1998} as a measure of stochastic dependence in complex systems. Under this perspective, the dual gradient flow induced by \cref{eq:gradMin} and \cref{eq:gradMax} can be interpreted as an instance of the generalized infomax principle. 
	
	To obtain a better view of this, we can first observe that the maximum entropy estimate gives the solution to the minimizing part of \cref{eq:manifoldMin}.
	Let 
	\begin{equation*}
		\hat{\bm{\theta}}^{1} =\argmin_{P_{\bm{\theta}^{1}} \in S^{1}} D_{KL}(P_{\bm{\theta}^{q^*}} \ || \ P_{\bm{\theta}^1} )
	\end{equation*}
	denote this estimate for some $P_{\bm{\theta}^{q^*}} \in S^{q^*}$ and $P_{\hat{\bm{\theta}}_{s}^{1}} \in S^1$ a split model corresponding to the maximum entropy estimate $\hat{\bm{\theta}}^{1}$ by setting all interaction parameters to zero. In our case, this means setting all the elements of the $\textbf{E}$ matrix \cref{eq:canonicalPar} in the exponential reparametrization of the full multivariate von Mises model \cref{eq:eMVM} to zero.	
	Since the Bregman divergence associated with the free energy \cref{eq:freeEnergyTrunc} induces a dually affine flat structure on $S$ \cref{eq:Sfam}, the information projection theorems and the generalized Pythagorean relation hold \cite{amari2000}.  That is, $\hat{\bm{\theta}}^{1}$ is given by an orthogonal geodesic projection of $P_{\bm{\theta}^{q^*}}$ onto $S^1$, meaning that the dual geodesic connecting $P_{\bm{\theta}^{q^*}}$ and $P_{\hat{\bm{\theta}}^{1}}$ and the geodesic connecting $P_{\hat{\bm{\theta}}^{1}}$ and $P_{\hat{\bm{\theta}}_{s}^{1}}$ are orthogonal intersecting at $P_{\hat{\bm{\theta}}^{1}}$. Thus, the KLD from $P_{\bm{\theta}^{q^*}}$ to $P_{\hat{\bm{\theta}}_{s}^{1}}$ can be decomposed into
	\begin{math}
		D_{KL}(P_{\bm{\theta}^{q^*}} \ || \ P_{\hat{\bm{\theta}}_{s}^{1}} )= 
		D_{KL}(P_{\bm{\theta}^{q^*}} \ || \ P_{\hat{\bm{\theta}}^{1}} ) + D_{KL}(P_{\hat{\bm{\theta}}^{1}} \ || \ P_{\hat{\bm{\theta}}_{s}^{1}} )
	\end{math}
	or alternatively in the form
	\begin{math}
		D_{KL}(P_{\bm{\theta}^{q^*}} \ || \ P_{\hat{\bm{\theta}}^{1} }) = 
		D_{KL}(P_{\bm{\theta}^{q^*}} \ || \ P_{\hat{\bm{\theta}}_{s}^{1}} ) - D_{KL}(P_{\hat{\bm{\theta}}^{1}} \ || \ P_{\hat{\bm{\theta}}_{s}^{1}} ).
	\end{math}
	Subsequently, for a fixed $P_{\hat{\bm{\theta}}^{1}}$ solving
	\begin{math}
		\max_{P_{\bm{\theta}^{q^*}} \in S^{q^*}} D_{KL}(P_{\bm{\theta}^{q^*}} \ || \ P_{\hat{\bm{\theta}}^{1}} )
	\end{math}
	in \cref{eq:manifoldMin} is equivalent to solving
	\begin{equation}\label{eq:splitModel}
		\max_{P_{\bm{\theta}^{q^*}} \in S^{q^*}} D_{KL}(P_{\bm{\theta}^{q^*}} \ || \ P_{\hat{\bm{\theta}}_{s}^{1}} ).
	\end{equation}
	
	Multi-information is defined as the KLD from the joint probability distribution $p\left(x_1,\ldots,x_n\right)$ to the product of the marginals $p\left(x_1\right)*\ldots*p\left(x_n\right)$ of a random vector $\textbf{X}=\left(X_1,\ldots, X_n\right)$. It can state that the components of the random vector $\textbf{X}$ are independent; that is, $p\left(x_1,\ldots,x_n\right) = p\left(x_1\right)*\ldots*p\left(x_n\right)$ if and only if the KLD or the multi-information vanishes. 
	Given a split model, that is a model without any interactions between state variables, parametrized by $\hat{\bm{\theta}}_{S}^{1}$ of the form \cref{eq:expFamY}, the components of the random vector associated with the split model are independent and the joint probability measure $dP(\hat{\bm{\theta}}_{s}^{1})$ can be decomposed into a product of marginal measures
	\begin{equation*}
		dP(\hat{\bm{\theta}}_{s}^{1})=\left[\prod_{i=1}^{n}\exp\left\{\left(\hat{\theta_{s}}^{1}\right)_{i}t_i\left(\left(\textbf{F}^{-1}(y)\right)_{i}\right)\right\}\right] \exp\left\{- \psi \left( \hat{\bm{\theta}}_{s}^{1} \right) \right\}dP_0 \circ \textbf{F}^{-1}(y),
	\end{equation*}
	where $t_i\left(x_i\right)$ is the $i$th component of the sufficient statistic $\textbf{t}(\textbf{x})=(t_1(x_1),\ldots,t_n(x_n))$. 
	
	Furthermore, the KLD in \cref{eq:splitModel} can be rewritten as
	\begin{equation*}
		D_{KL}(P_{\bm{\theta}^{q^*}} \ || \ P_{\hat{\bm{\theta}}_{s}^{1}} )=
		\int \limits_{\left\{\textbf{F}_{1-\frac{1}{q^*}}\leq\textbf{F}\left(\textbf{x}\right)\right\}} \ln \frac{dP(\bm{\theta}^{q^*})}{ dP(\hat{\bm{\theta}}_{s}^{1})}dP(\bm{\theta}^{q^*}).
	\end{equation*}
	due to $dP(\bm{\theta}^{q^*})=0$ for $\textbf{F}_{1-\frac{1}{q^*}} > \textbf{F}\left(\textbf{x}\right)$.	
	Since $q^*$ is constant for a given iteration, the dual gradient ascent \cref{eq:geoGradTheta} and \cref{eq:geoGradMu} and the related entropic proximal ascent in dual coordinates \cref{eq:thetaProx} and \cref{eq:muProximal} can be interpreted as moving towards maximizing stochastic dependence. In other words, they are increasing the amount of information shared among elements of the random vector $\textbf{X}$ due to the maximization of multi-information from the truncated exponential family $S^{q^*}$ defined by the selection quantile \cref{eq:qES} to the exponential family $S^{1}$ with full support.

	\subsection{Refining solutions}\label{sec:refining}
	Stochastic optimization methods are generally well suited for problems where approximate solutions are sufficient. Overall, it is practically impossible to sample the exact solution, even when, in theory, an algorithm converges to a delta distribution concentrated on the solutions of the optimization problem. The entropic trust region \cref{eq:updateIGO} is no exception. To attain higher accuracy, we use the following procedure to refine the solution attained after the initial optimization by reducing the configuration space to a smaller region.
	
	An $\bm{\epsilon}^{r}$ neighbourhood is created around the current best solution $\textbf{x}^{\text{best}}$. The optimal so far solution at the $r$th run is recorded as
	\begin{math}
		\textbf{x}^{\text{best}}=\{\textbf{x}^s | \textbf{F}(\textbf{x}^s)\geq \textbf{F}(\textbf{x}_i), \ i=0,1\dots,r \},
	\end{math}
	with  $\bm{\epsilon}^{r}$ set to
	\begin{equation}\label{eq:epsilonC}
		\bm{\epsilon}^r = \left(\frac{1}{c_{\epsilon}}\right)^r\left(\textbf{u}-\textbf{l}\right),
	\end{equation}
	where $c_{\epsilon}>1$ is the rate of exponential decay, $r=1,2,\ldots$ denotes the current run of the refining process and $\textbf{u}, \ \textbf{l}$ are the initial upper and lower bounds, given by the optimization problem. 
	
	Based on $\bm{\epsilon}^{r}$, we reduce the two bounds $\textbf{u}_\epsilon, \ \textbf{l}_\epsilon$ in the next run by setting
	\begin{equation}\label{eq:epsilonBoundary}
		\begin{array}{l}
			\textbf{u}^\epsilon=\min\left({\textbf{x}^{\text{best}}+\bm{\epsilon}^r},\textbf{u}\right), \\
			\textbf{l}^\epsilon=\max\left({\textbf{x}^{\text{best}}-\bm{\epsilon}^r},\textbf{l}\right),
		\end{array}
	\end{equation}
	taken along respective optimization variables and proceed with the next iteration of the algorithm set to the new upper and lower bounds $\textbf{u}_\epsilon, \ \textbf{l}_\epsilon$.
	We iterate this process until the desired accuracy is achieved. For example, in the setting of CSG packings, accuracy can be measured by the minimal Euclidean distance between sets $gK$ of a $G$-packing $\mathcal{K}_G$ in \cref{eq:spacepacking}. The algorithm concentrates on the expected optimum with higher precision with each run. This refining process can be interpreted as a variant of a restart strategy.
	
	\begin{figure}
		\centering
		\includegraphics[trim={60 32 -10 10},clip,height=0.4\linewidth]{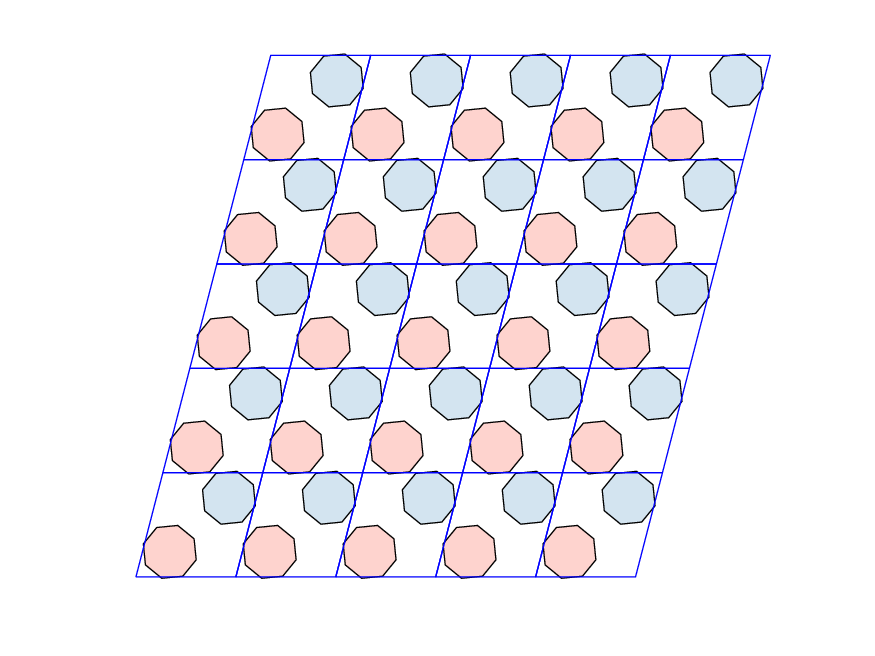}
		\includegraphics[trim={22 -90 21 -90},clip,height=0.4\linewidth]{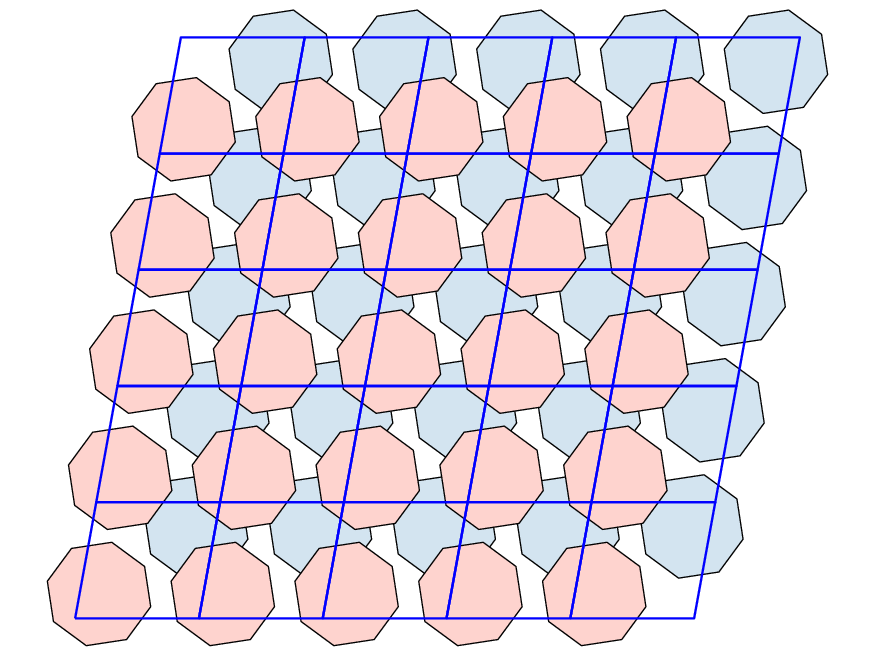}
		\caption{Visualization of $25$ cells of two $p2$ configurations of a regular octagon from the initial sampling in \cref{fig:samplehist}. (Left) Feasible solution of a $p2$-packing with density $\rho \left(\mathcal{K}_{p2} \right) =0.413705837593271$. (Right) Unfeasible solution due to overlapping.}
		\label{fig:initExamples}
	\end{figure}
	
	\begin{figure}
		\centering
		\begin{subfigure}[b]{0.49\textwidth}
			\includegraphics[trim={30 0 0 0},clip,width=1\linewidth]{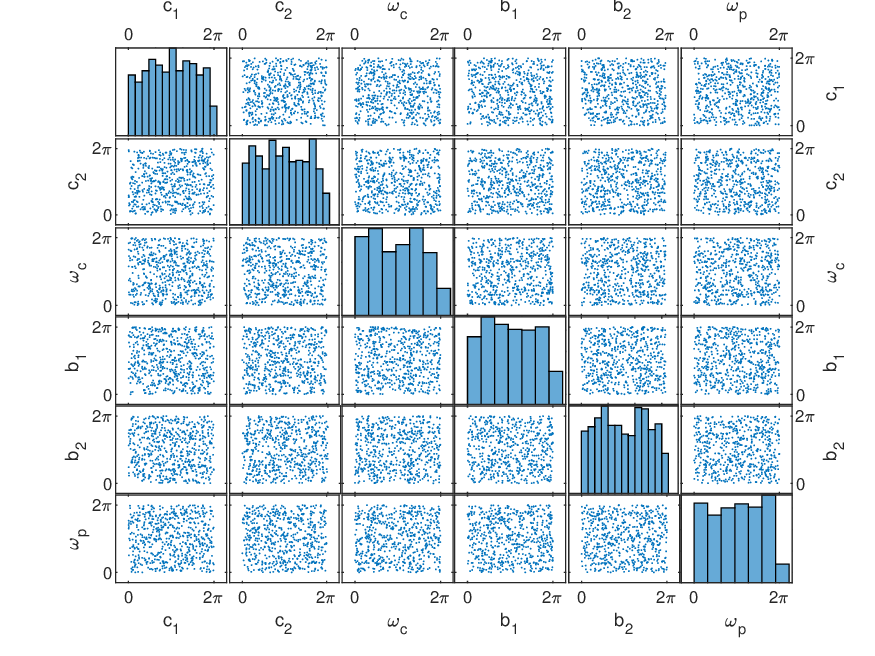}
		\end{subfigure}
		\begin{subfigure}[b]{0.49\textwidth}
			\includegraphics[trim={30 0 0 0},clip,width=1\linewidth]{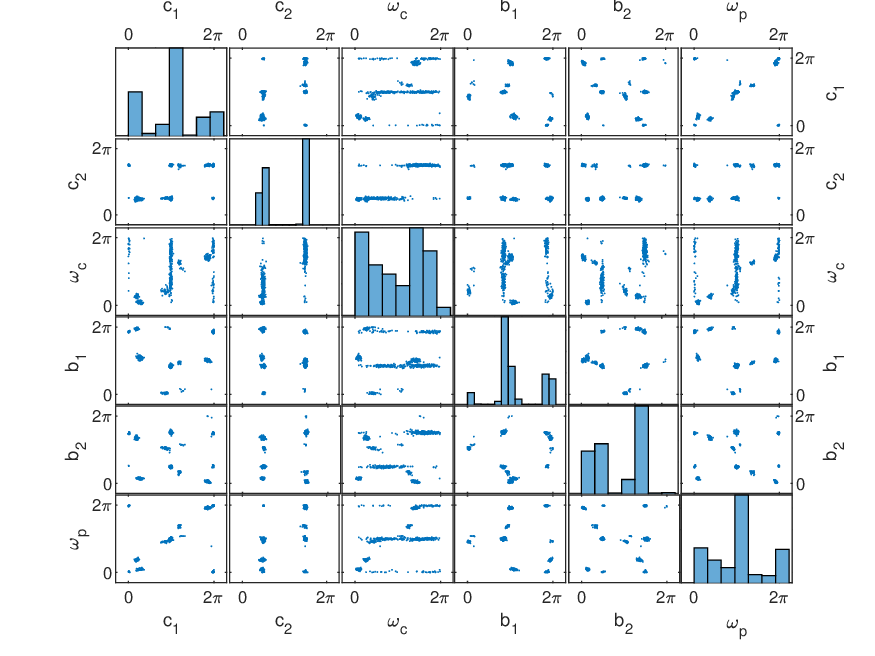}
		\end{subfigure}
		\caption{2D projections along coordinate axes and histograms of univariate marginals corresponding to the respective optimization variables of $600$ realizations of the exponential reparametrization of the full multivariate von Mises distribution. (Left) Initial distribution. (Right) Output distribution.}
		\label{fig:samplehist}
	\end{figure}
	
	\section{Experiments:  A case study for dense plane group packings of convex polygons}
	\label{sec:proofofconcept}
	
	In this section, we analyze the behaviour and performance of the introduced entropic trust region on the problem of packing convex polygons in 2D CSGs. In \cref{sec:spacegrouppacking}, we introduced these symmetry groups, also referred to as plane groups, and defined a subproblem of the general packing problem \cite{rogers}. Specifically, we present an experimental examination for the case study of the densest $p2$-packings of a regular octagon and results from implementations of our packing algorithm to various plane group packings of convex polygons in instances for which the optimal solutions are known \cref{tab:tabel}. More experimental results are presented in the supplementary material. That is additional details related to the behaviour of entropic trust region in $p2$-packing of regular octagons test case (examination of distribution parameter trajectories \cref{sec:trajectories} and dependence of stability on the sampling pool size  \cref{sec:samplingPool}) and a more detailed presentation of results from additional plane group packings (\cref{sec:additionalPackings}).		
	Technical details concerning the implementation of the entropic trust region method to the densest plane group packings are discussed in \cref{sec:implementation} \footnote{Matlab source code, instructions to use and examples are available at \url{https://milotorda.net/software/}.}.
	
	Here, we chose the problem of $p2$-packing of convex octagons to showcase the behaviour and performance of the algorithm for several reasons. First, the general optimal packing of octagons is known. Second, the problem, as we have defined it in section \ref{sec:spacegrouppacking}, has at least $64$ global maxima due to the $8$-fold rotational symmetry of the octagon. The way the plane groups are constructed and the initial boundary values make this case more challenging than the other test cases.
	
	\begin{figure}
		\centering
		\begin{subfigure}[!b]{0.49\textwidth}
			\includegraphics[width=1\linewidth]{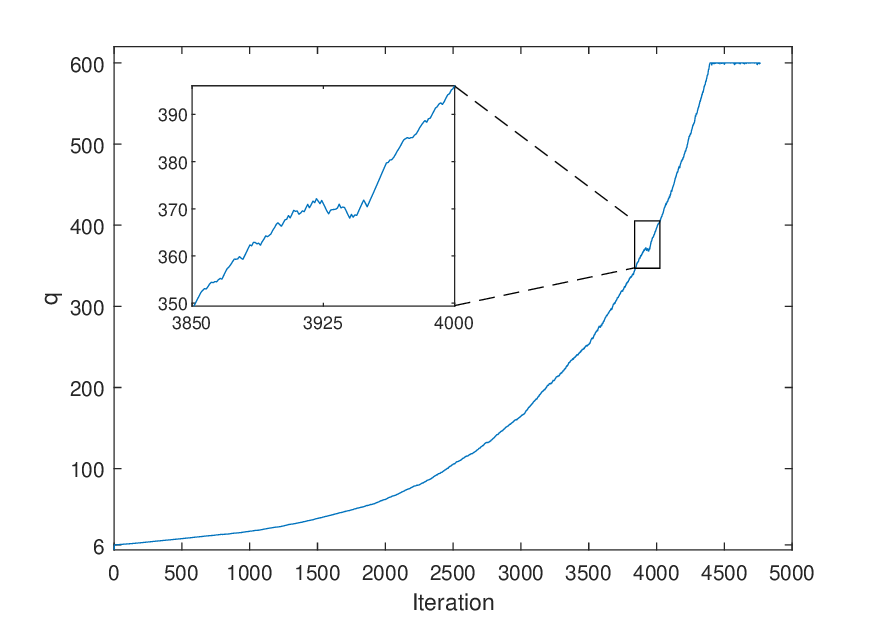}
		\end{subfigure}
		\begin{subfigure}[!b]{0.47\textwidth}
			\includegraphics[width=1\linewidth]{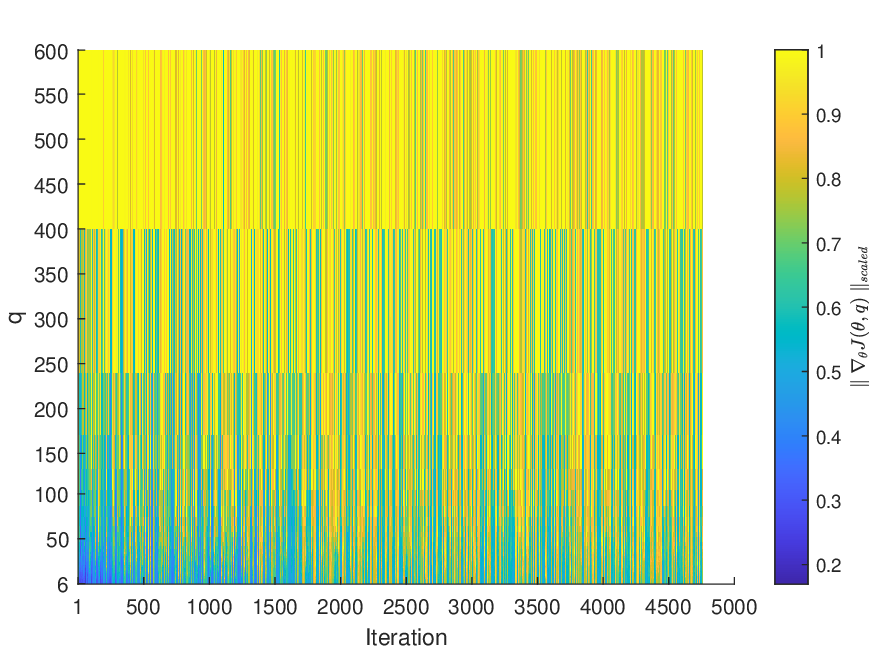}
		\end{subfigure}
		
		\begin{subfigure}[!b]{0.94\textwidth}
			\includegraphics[width=1\linewidth]{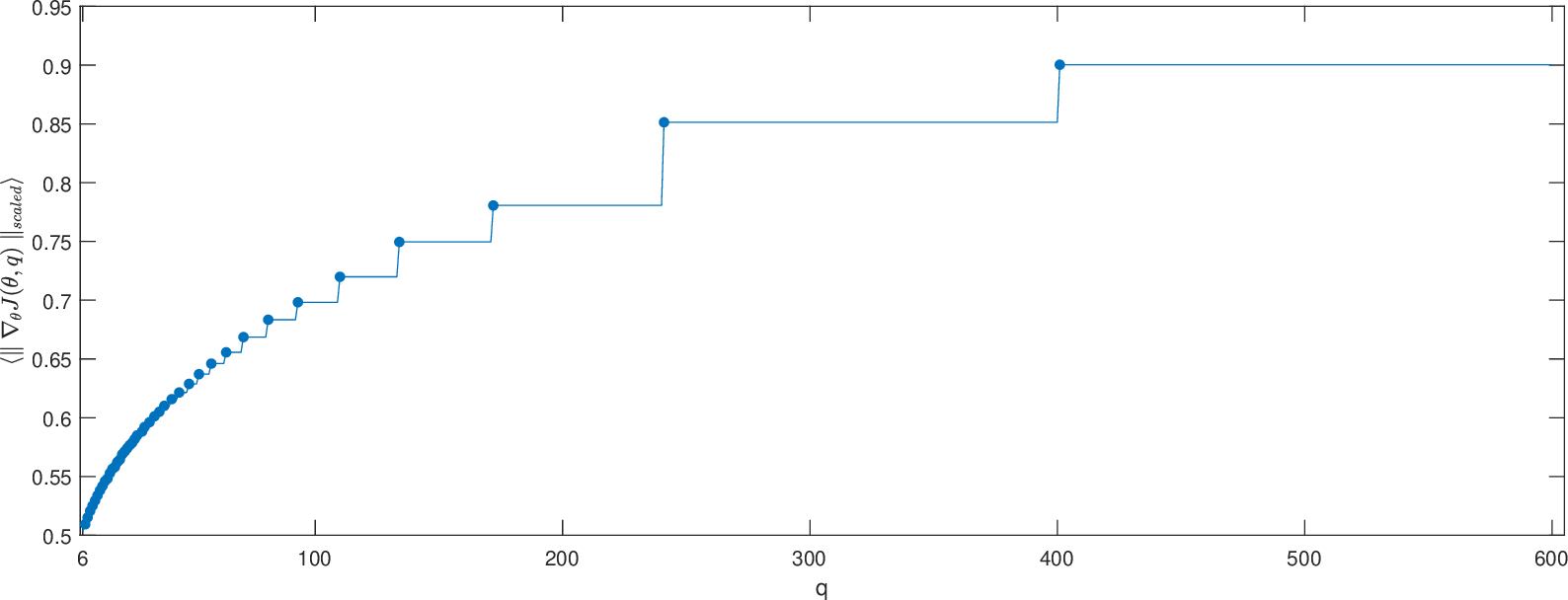}
		\end{subfigure}
		\caption{(Top left) Evolution of the adaptive selection quantile $q$. (Top right) Relationship between selection quantile $q$  and the scaled expected fitness gradient \cref{eq:gradScaled} at every iteration of the entropic trust region run, and (Bottom) the scaled expected fitness gradient averaged through iterations where markers denote changes in the value of $<\parallel\nabla_{\theta} J(\theta,q)\parallel_{scaled}>$.}
		\label{fig:qViz}
	\end{figure}
	
	The $p2$ plane group is a semi-direct product of a point group and a lattice group introduced in \cref{sec:spacegrouppacking}. The point group is given by two symmetry operations, that is the identity operation ($\textbf{R}=\textbf{I}$) and a rotation by $180^\circ$. 
	The lattice group, given by the oblique crystal system, has three degrees of freedom, namely the length of the lattice group generators or basic vectors $\textbf{b}_1$, $\textbf{b}_2$ and an angle $\omega_p$ between them, constituting part of the design variables to optimize.
	Additional degrees of freedom of the optimization configuration space are given by the fractional coordinates (a coordinate system where basis vectors are defined by the lattice basic vectors \eqref{eq:latticeG}) of the position of the centroid of the octagon in the asymmetric unit $\textbf{c}_1$, $\textbf{c}_2$ and an angle of rotation $\omega_c$ of the octagon. Altogether, we have six optimization variables for this specific test case.
	
	The optimal general packing of centrally symmetric convex polygons is lattice packing \cite{rogers1951}. Since the regular octagon is a polygon with central symmetry, the packing is known with the density of $\rho_{\text{opt}}=\frac{4+4\sqrt{2}}{5+4\sqrt{2}} \approx 0.90616367$.
	
	The optimization boundaries are set as follows. From the asymmetric unit restrictions of the $p2$ group, octagon centroid fractional coordinates $c_1 \in \left(0,1\right]$ and $c_2 \in \left(0,\frac{1}{2}\right]$, angle of rotation of the octagon $\omega_c \in \left(0,2\pi \right]$, lengths of the lattice generators $b_i \in \left[0,2d\right]$ for $i=1,2$ where $d$ denotes the diameter of the octagon's circumcircle, and angle between lattice generators $\omega_p \in \left[0,\frac{\pi}{2}\right]$. 
	
	\begin{figure}
		\centering
		\begin{subfigure}[!b]{0.49\textwidth}
			\includegraphics[clip,width=1\linewidth]{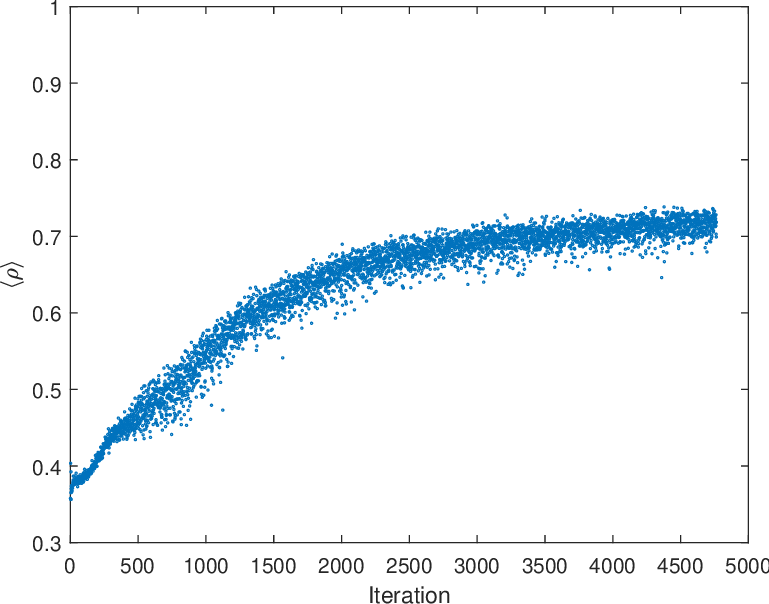}
		\end{subfigure}
		\begin{subfigure}[!b]{0.49\textwidth}
			\includegraphics[width=1\linewidth]{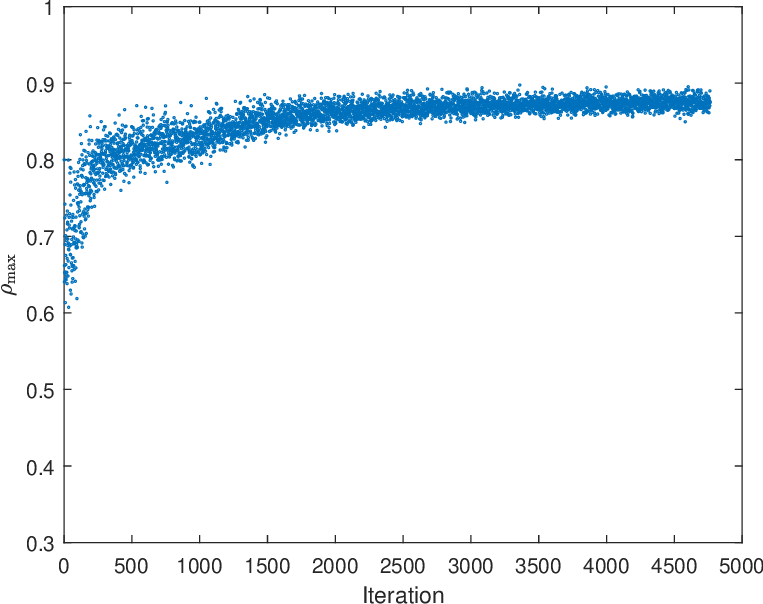}
		\end{subfigure}
		\caption{Evolution of (left) arithmetic mean and (right) maximum packing density of $600$ solutions generated in each iteration.}
		\label{fig:trajectories}
	\end{figure}
	
	\begin{figure}
		\centering
		\begin{subfigure}{0.54\textwidth}
			\includegraphics[width=1\linewidth]{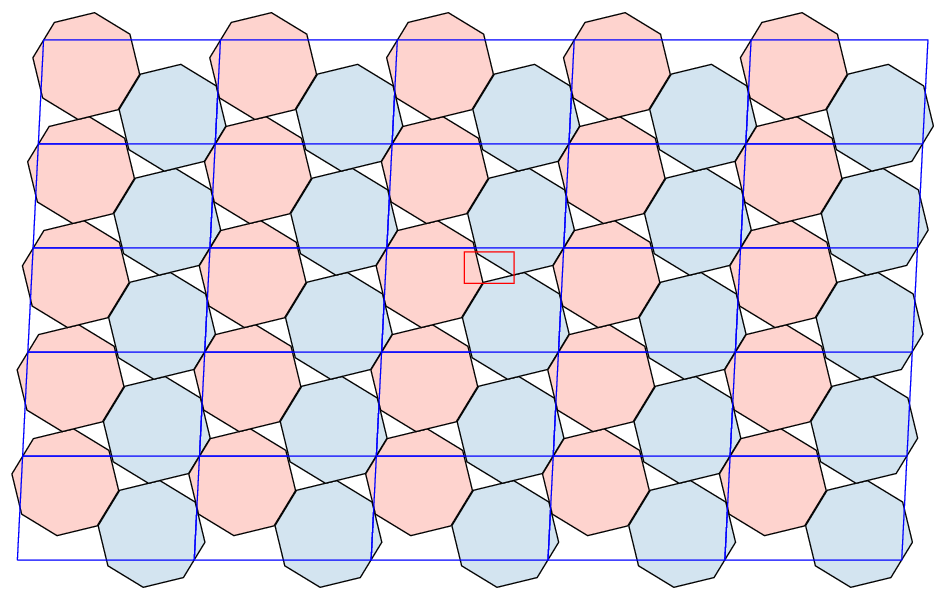}
		\end{subfigure}		
		\begin{subfigure}{0.45\textwidth}
			\includegraphics[width=1\linewidth]{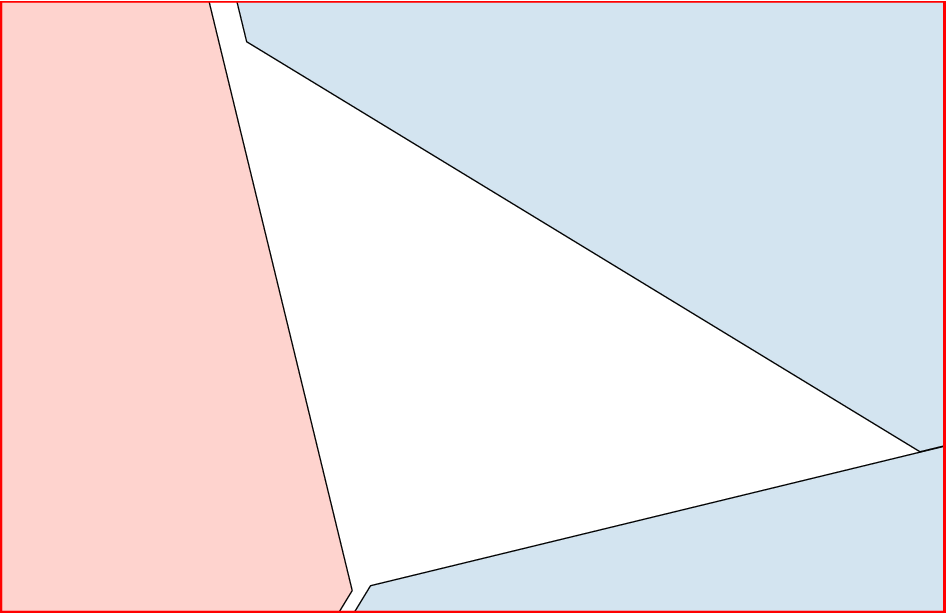}
		\end{subfigure}
		\caption{Best solution found during the initial optimization run. (Left) $25$ cells drawn from the whole packing configuration. (Right) An enlargement of the configuration segment marked by a red rectangle.}
		\label{fig:P2Gon8Run}
	\end{figure}
	
	At the beginning of each optimization run all of the parameters of the exponential reparametrization of the full multivariate von Mises distribution $\left(\bm{\eta},\textbf{E}\right)$ in \cref{eq:eMVM} are set to zero, and $600$ realizations from this distribution are generated. Effectively, at the start of each optimization run, we are sampling from a uniform distribution on a $6$D flat torus, providing a homogeneous and unbiased first exploration of the optimization landscape. \Cref{fig:samplehist} illustrates $2$D projections along the respective coordinate axes of $600$ samples generated in this way.	
	Two candidate solutions from these initial $600$ samples are shown in \cref{fig:initExamples}.
	
	Since there is no correspondence between the full multivariate von Mises distribution \cref{eq:fmvm} and its exponential form \cref{eq:eMVM} for $\left(\bm{\eta},\textbf{E}\right)=\bm{0}$, the initial trust region update is done using \cref{eq:updateIGO} for the exponential von Mises distribution \cref{eq:eMVM}. By doing so, we remove the dependency of the optimization on the initial position that burdens many optimization methods. After the initial step, all updates are performed in the $\left(\bm{\mu},\bm{\kappa},\textbf{D}\right)$ parametrization \cref{eq:fmvm}.
	
	During execution, the $q$ parameter of the adaptive selection quantile \cref{eq:adaSQ} gradually increases, as shown in \cref{fig:qViz} and intended by our design. The starting value of $q$ is $6$, meaning that at the beginning, $99\%$ of the $600$ samples are taken into account for the truncated exponential distribution expected value $\bm{\mu}_{\textbf{F}_{1- \frac{1}{q_{t}}}}$ \cref{eq:qTrunctE}. After about $4,000$ iterations, $q$ reaches $600$, which is the overall number of samples used in each iteration. Meaning, that only a single realization from the sampling distribution with the highest fitness is assigned to $\bm{\mu}_{\textbf{F}_{1- \frac{1}{q_{t}}}}$, and the algorithm moves in the direction of this point. An interesting observation from the evolution of $q$ is that the directions between the vectors in the tangent space $T_{\bm{\theta}^{t-1}}S_{\bm{\theta}}$ measured by \cref{eq:diffangle} does not change much, resulting in only small fluctuations in the $q$ evolution trajectory.
	
	\begin{figure}
		\centering
		\includegraphics[trim={0 20 0 0},clip,width=1\linewidth]{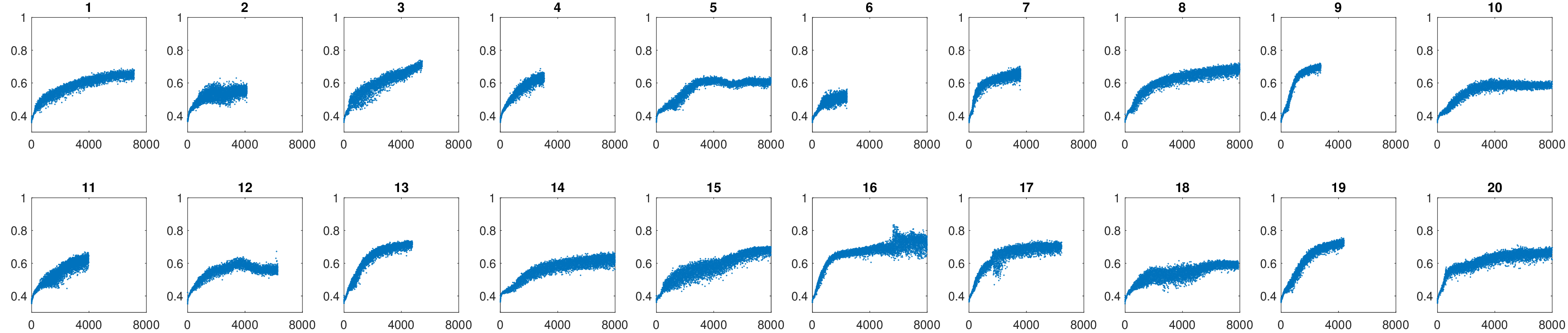}
		\caption{Evolution of the average density $\left<\rho\right>$ in $20$ optimization runs.}
		\label{fig:paperExpMean}
	\end{figure}
	\begin{figure}
		\centering
		\includegraphics[trim={0 0 0 20},clip,width=1\linewidth]{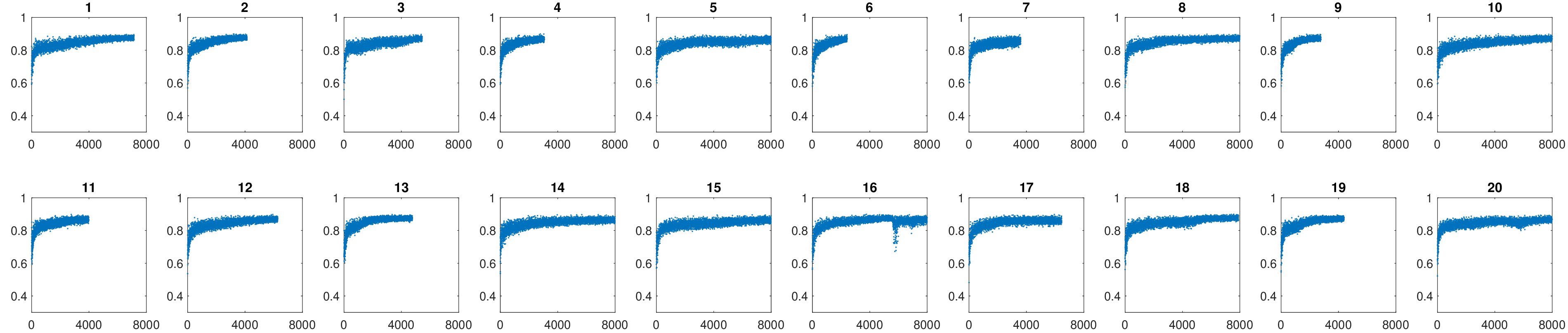}
		\caption{Evolution of maximum density $\rho_{\max}$ in $20$ optimization runs.}
		\label{fig:paperExpMax}
	\end{figure}
	
	The adaptive expected fitness gradient \cref{eq:adaThetaGrad} with increasing $q$ not only forces the trust region to move in the direction of the distribution that most likely represents areas with the highest packing density, but by doing so, it speeds up the rate of convergence. To illustrate this, we computed the $L_1$ norm of $\bm{\mu}_{\textbf{F}_{1- \frac{1}{q_{t}}}}-\bm{\mu}$ at every optimization iteration for the values of $q=6,\ldots,600$. We scale the results such that the maximum at every iteration is not greater than one, precisely
	\begin{equation}\label{eq:gradScaled}
		\parallel\nabla_{\bm{\theta}} J(\bm{\theta}^t,q)\parallel_{\text{scaled}} \ = \ \frac{\parallel\bm{\mu}^{t}_{\textbf{F}_{1- \frac{1}{q}}}-\bm{\mu}^{t}\parallel_{1}}{\parallel\bm{\mu}^{t}_{\textbf{F}_{1- \frac{1}{q^{\max}}}}-\bm{\mu}^{t}\parallel_{1}},
	\end{equation}
	where $q^{\max}=\argmax_{q}\parallel\bm{\mu}^{t}_{\textbf{F}_{1- \frac{1}{q}}}-\bm{\mu}^{t}\parallel_{1} \ , \ q=6,\ldots,600$.
	
	Although there are instances where
	\begin{math}
		\parallel\nabla_{\bm{\theta}} J(\bm{\theta}^t,q_1)\parallel_{\text{scaled}} \ \leq \ \parallel\nabla_{\bm{\theta}} J(\bm{\theta}^t,q_2)\parallel_{\text{scaled}}
	\end{math}
	for some $q_1<q_2$ does not hold, \cref{fig:qViz} shows that, on average, the higher the $q$, the larger the gradient size in $\bm{\theta}$-coordinates. 
	This kind of behaviour is beneficial since the adaptive selection quantile accelerates convergences at later stages of the execution when an attractor has already been singled out.
	
	The main principle of the algorithm is maximizing the expected fitness $J\left(\bm{\theta}\right)$ \cref{eq:efit}. In our case, it means maximizing the packing density \cref{eq:density} or minimizing primitive cell volume. \cref{fig:trajectories} shows the evolution of the average density at each iteration defined by
	\begin{math}
		\left<\rho\right>=\frac{1}{N}\sum_{i=1}^{N}\frac{2\text{AREA}\left(K\right)}{\textbf{F}\left(\mathbf{x}_i \right)},
	\end{math}
	where $K$ is the regular octagon, $N$ is the number of candidate solutions $\mathbf{x}_i$ sampled at each iteration and $\textbf{F}(\mathbf{x} )$ is the penalty function \cref{eq:penalty} with the objective function defined as the area of the unit cell $f(\mathbf{x}):= \text{AREA}(\overline{\Lambda}_{p2}(\mathbf{x}))$ and the constraint violation of the form $g\left(\mathbf{x}\right):=\text{dist}\left(\mathcal{K}_{p2}(\mathbf{x})\right)$ defined by \cref{eq:distPack}. The algorithm gradually increases the average density $\left<\rho\right>$ with the maximum attained at the $4,378$th iteration. 
	
	\begin{figure}
		\centering
		\begin{subfigure}{0.5\textwidth}
			\includegraphics[width=1\linewidth]{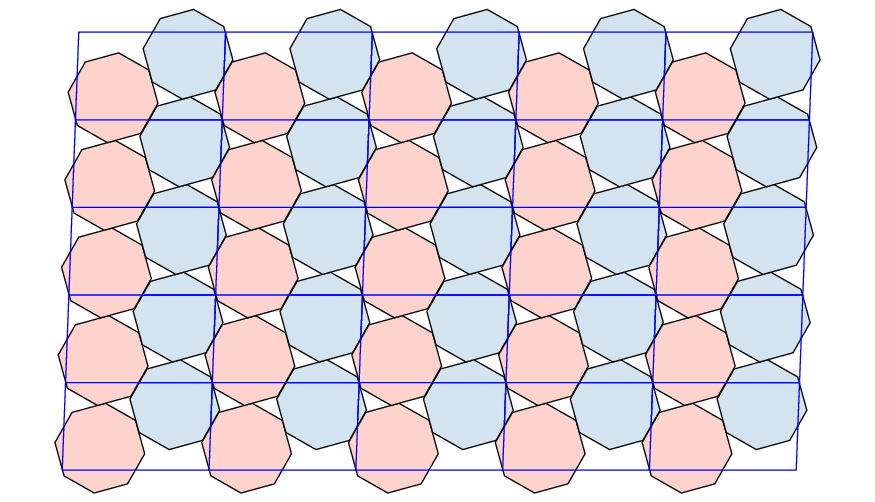}
		\end{subfigure}
		\begin{subfigure}{0.49\textwidth}
			\includegraphics[width=1\linewidth]{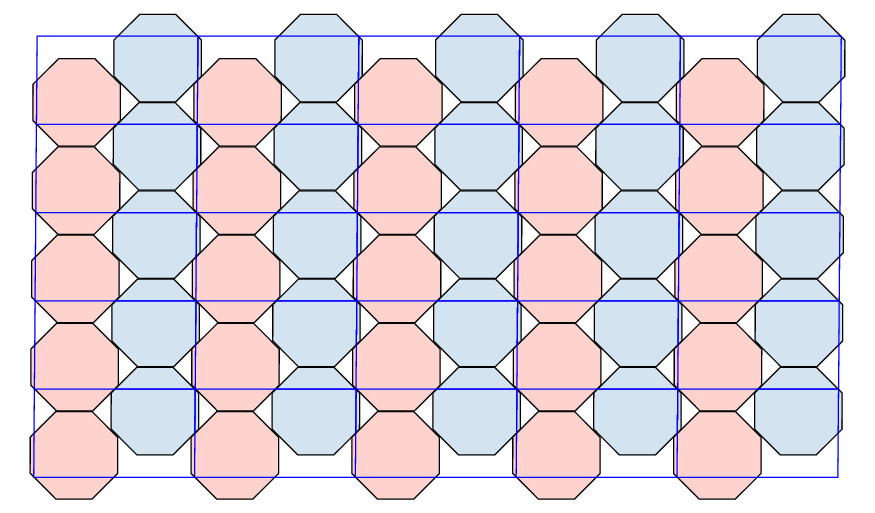}
		\end{subfigure}
		\caption{(Left) Best and (Right) worst densest packing configurations found in $20$ optimization runs.}
		\label{fig:p2gon8BestWorst}
	\end{figure}
	
	\begin{figure}
		\centering
		\includegraphics[width=0.5\linewidth]{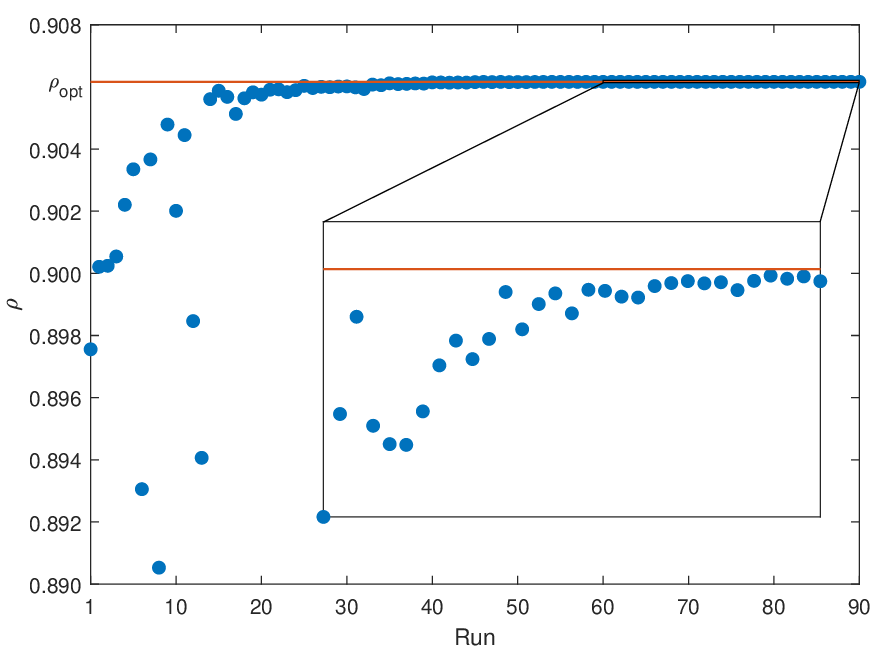}
		\caption{Densities of the best solutions in each run of the packing refining process. $\rho_{opt}$ denotes the theoretical optimal density.}
		\label{fig:paperplotruns}
	\end{figure}	
	
	Similar behaviour can be observed in the evolution of the maximum packing $\rho_{\max}$ found at each iteration, shown in \cref{fig:trajectories}. The best solution was found at the $3,360$th iteration with packing density $\rho \left(\mathcal{K}_{p2} \right)=0.897526117081202$ and minimal Euclidean distance between octagons in the configuration $\text{dist} \left(\mathcal{K}_{p2} \right)=0.000385651690559$. A visualization of $25$ cells from this configuration is presented in \cref{fig:P2Gon8Run}. The difference from the theoretical optimal packing density $\rho_{\text{opt}}$ defined as
	\begin{equation}\label{eq:packDiff}
		\Delta_{\mathcal{K}_{p2}}= \rho_{\text{opt}}-\rho (\mathcal{K}_{p2} ),
	\end{equation}
	is $\Delta_{\mathcal{K}_{p2}}=0.008637561562744$.	
	
	As the packing algorithm is stochastic in nature, it is unlikely to repeat the same result. Therefore, the maximum packing density attained during one execution can be considered a random variable by itself. To assess performance and robustness, we perform $20$ consecutive runs with the same hyperparameters (\cref{tab:1}) and the seed of the uniform distribution pseudorandom number generator used in the Gibbs sampler \cref{eMvMgibbs} initialized using system time. \cref{fig:paperExpMean} shows the evolutions of average density, and \cref{fig:paperExpMax} of maximum density. The Hodges-Lehmann estimator of the pseudomedian of the maximum packing density computed from the maximum densities attained in each of the $20$ runs is $\hat{m}=0.8970032$ with the $95\%$ confidence interval based on Wilcoxon's signed rank test equal to $\text{CI}=\left(0.8959246, 0.8980686\right)$. On average, the maximum packing density was attained at the $4,624$th iteration. The highest maximal packing density configuration was attained in run $10$ with packing density $\rho \left(\mathcal{K}_{p2} \right)=0.899848387789551$, and the lowest maximal packing density was attained in run $20$ with packing density $\rho \left(\mathcal{K}_{p2} \right)=0.893366540122926$. Both solutions are shown in \cref{fig:p2gon8BestWorst}. In a closer examination, it can be noticed that both configurations look similar, which is not surprising considering they both represent solutions from different global optima basins due to the multi-modality of the problem stemming from the symmetries of the regular octagon.
	
\begin{figure}
		\centering
		\begin{subfigure}[!b]{0.54\columnwidth}
			\includegraphics[width=1\linewidth]{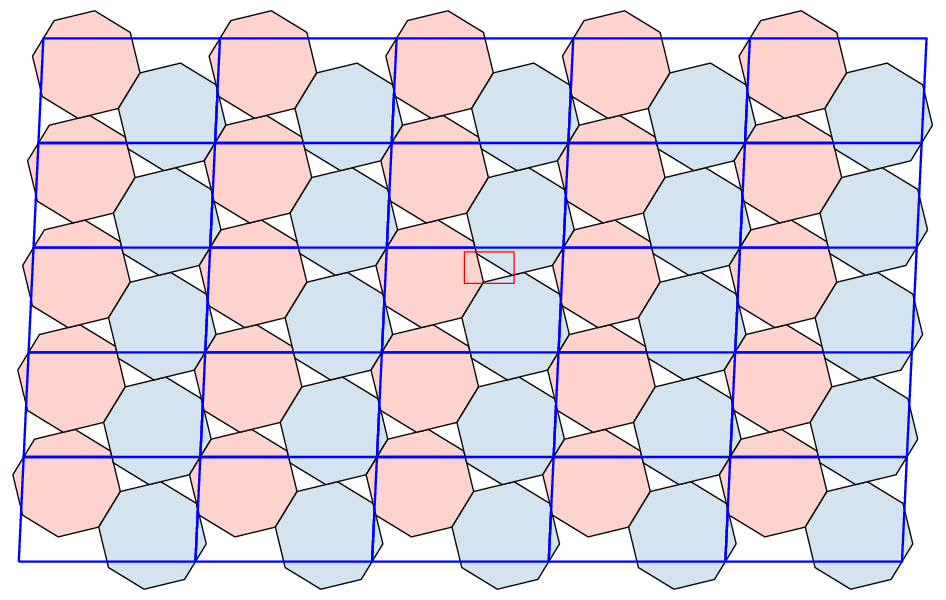}
		\end{subfigure}		
		\begin{subfigure}[!b]{0.45\columnwidth}
			\includegraphics[width=1\linewidth]{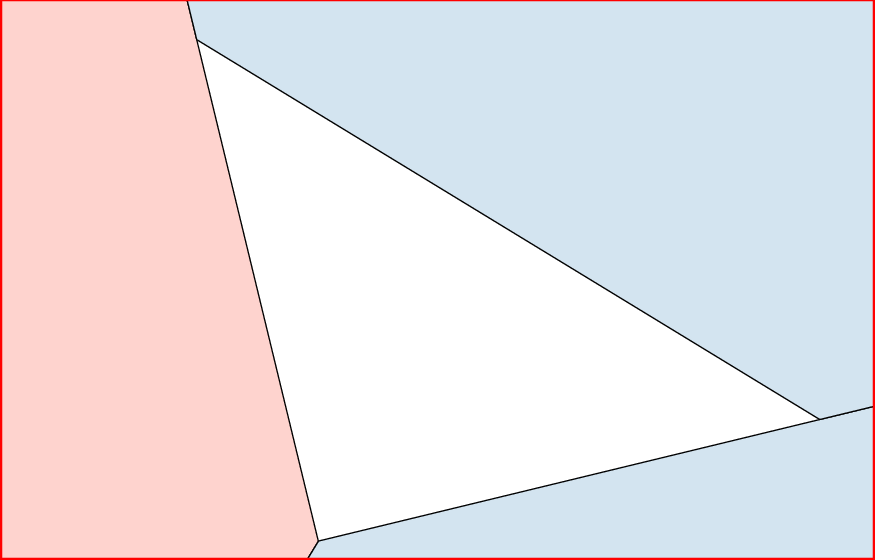}
		\end{subfigure}
		\caption{Refined solution. (Left) $25$ cells from the whole configuration. (Right) The red rectangle marks an enlargement of the configuration segment.}
		\label{fig:P2Gon8Last}
	\end{figure}

	After the initial run, although the maximal packing density configuration found during the run of \cref{fig:P2Gon8Run} visually resembles the theoretical optimal packing of regular octagons, the difference of packing densities $\Delta_{\mathcal{K}_{p2}}$ \cref{eq:packDiff} is in the $3$rd decimal place. This can also be observed either by noticing the minimal Euclidean distance between octagons in the packing $\text{dist} \left(\mathcal{K}_{p2} \right)$ defined by \cref{eq:distPack} or visually by enlarging a part of the packing, as is shown in \cref{fig:P2Gon8Run}. Therefore, as introduced in \cref{sec:refining}, we perform a refining process by taking the configuration with maximal packing density attained in the initial run (\cref{fig:P2Gon8Run}) and creating an $\bm{\epsilon}^r$ neighbourhood around this configuration's coordinates. In this way, we define a new configuration space with the boundaries given by \cref{eq:epsilonBoundary} and rerun the algorithm with these new boundaries. At this point the optimization variables $c_1$, $c_2$ and $\omega_c$ lose their inherent periodicity and have to be treated as aperiodic in the boundary mapping.   \cref{fig:paperplotruns} illustrates the convergence of the maximum density attained in each run during $90$ runs of the decaying boundary neighbourhood. The highest packing density was attained at run $87$ with the value $\rho \left(\mathcal{K}_{p2} \right)=0.90616363432568$ and the theoretical optimum difference $\Delta_{\mathcal{K}_{p2}}=3.1481 \ 10^{-8}$. A visualization of this configuration is presented in \cref{fig:P2Gon8Last}. For visual comparison, we include enlargement in \cref{fig:P2Gon8Last} of the same area of the output packing configuration as in the initial run (\cref{fig:P2Gon8Run}).
	
	\begin{figure}[!t]
		\centering
		\includegraphics[trim={120 90 110 80},clip,width=0.5\linewidth]{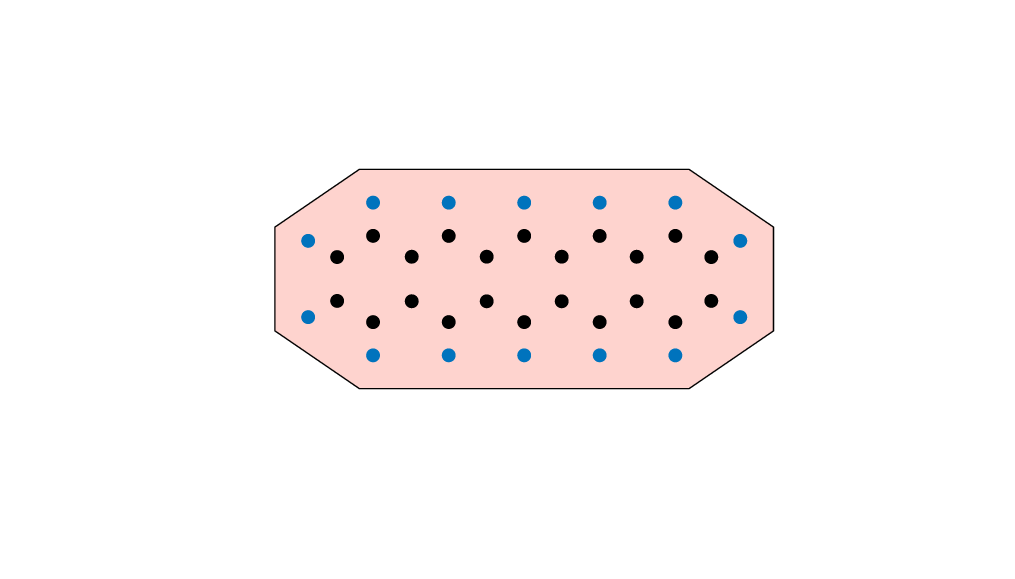}
		\caption{A geometric representation of pentacene. The dots symbolize atomic positions of (blue) hydrogen and (black) carbon.}
		\label{fig:pentacene2D}
	\end{figure} 
	
	\begin{figure}[!t]
		\centering
		
		\includegraphics[trim={0 0 0 0},width=\textwidth]{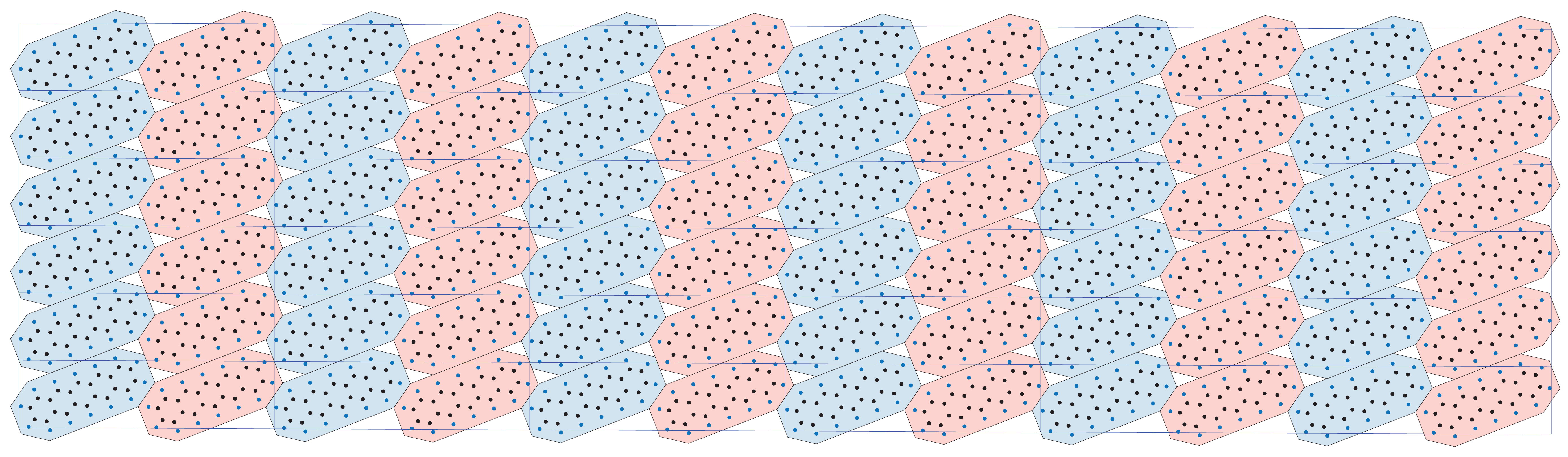}
		
		\includegraphics[trim={0 0 0 0},width=\textwidth]{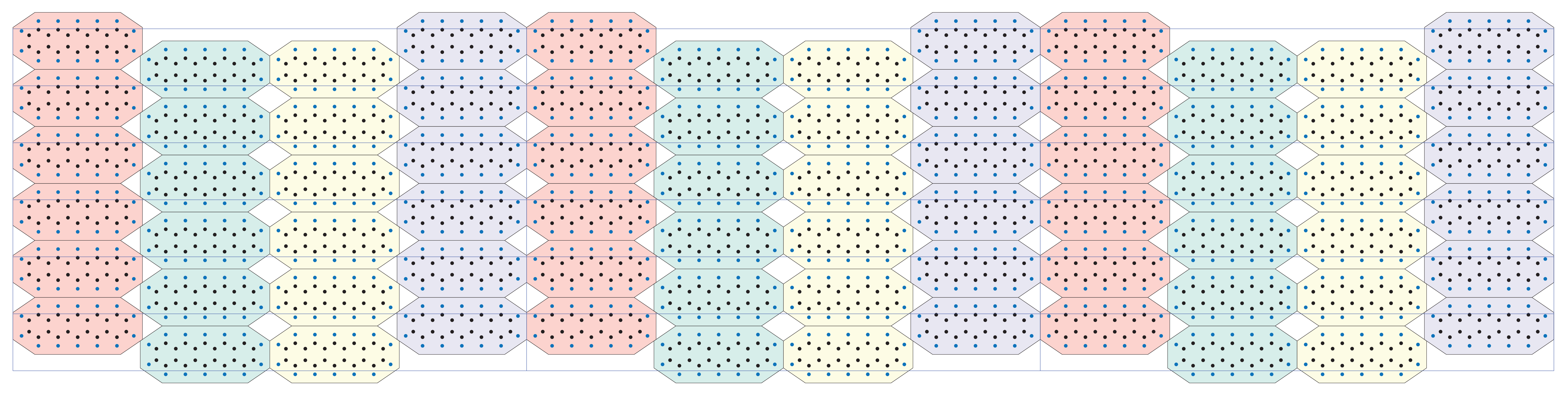}
		
		\includegraphics[trim={0 0 0 0},width=\textwidth]{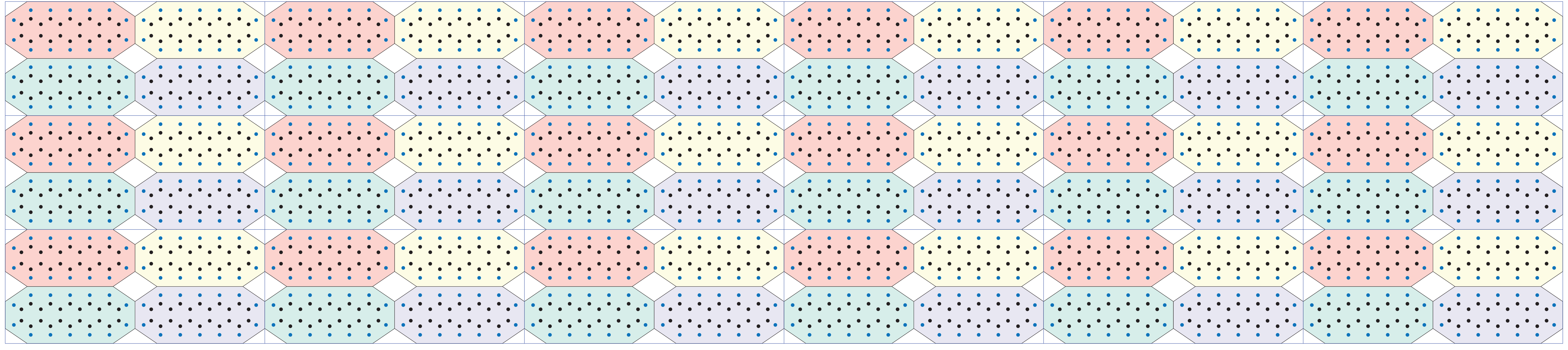}

		\caption{Visualization of the output configurations of the densest pentacene representation (top) $p2$-packing, (middle) $cm$-packing and (bottom) $p2mm$-packing. Colors represent CSG symmetry operations modulo lattice translations.}
		\label{fig:pentaceneP2}
	\end{figure}
	
	\section{Pentacene representation packings}
	\label{sec:pentacene}
	
	We demonstrate how the densest CSG packings are intended to be used in a molecular CSP workflow on pentacene thin-films. First, a geometric representation of a molecule by a polytope and an invertible map from the molecule's atomic positions to the polytope's interior is constructed. Afterwards, the densest packings of this representation are obtained in various CSGs. Lastly, the parameters of output configurations are used as input parameters for CSP computations, reducing the CSP search only to the neighbourhood of these configurations. 
	
	Pentacene is a planar molecule consisting of five serially connected benzene rings \cite{campbell1961}, explored as an organic thin-film semiconductor \cite{mei2013}.	
	Since the molecules in a crystal do not touch due to repulsive intermolecular forces, we built a $2$D representation of pentacene as the convex hull of the atomic positions of the molecule with an offset given by hydrogen's van der Waals radius of $1.09$\AA{} \cite{rowland1996}. The result is an irregular octagon illustrated in \cref{fig:pentacene2D}. 
	
	We employed the entropic trust region packing algorithm to search for this representation's densest plane group packings. \cref{fig:pentaceneP2} presents output configurations of the densest $p2$, $cm$ and $p2mm$ packings. 
	The approximate density of the $p2$-packing is $\rho \left(\mathcal{K}_{p2} \right)=0.953382110797399$ and resembles the configuration of single layer pentacene thin-film on graphite surface found in \cite{chen2008}. Moreover, this structure was found via simulation of self-assembly of a disordered system of pentacene molecules on graphene surface driven by the minimization of molecule-molecule interactions using the Lennard--Jones potential \cite{zhao}.	
	The densities of the output configurations in the $cm$ and $p2mm$ instances are $\rho \left(\mathcal{K}_{cm} \right)=0.918715231405704$ and $\rho \left(\mathcal{K}_{p2mm} \right)=0.910916580526171$. These configurations represent pentacene monolayer crystal phases on a Cu(110) surface found in \cite{sohnchen2004}.
	
	\section{Conclusions}
	\label{sec:conclusions}
	
	The problem of molecular CSP is to predict stable periodic structures from the knowledge of the chemical composition of a molecule. The most straightforward formulation means finding minima on a complicated energy landscape induced by one of the many free energy potentials. This presents a formidable task for optimization methods and search algorithms, leading in many cases to an over--prediction of polymorphic structures \cite{price2018}.  Providing current CSP solvers with densely packed initial
	configurations in terms of the geometric representation of a molecule
	can significantly accelerate CSP, as opposed to random starting
	structures, especially due to the recent complete isometry invariants
	of periodic structures
	\cite{anosova2021isometry,widdowson2022average,widdowson2022resolving}.
	
	The densest packing of geometric shapes is a notoriously hard problem in discrete and computational geometry \cite{toth} and is used in a large body of work in solid-state physics modelling \cite{torquato2018}. Since we are interested only in a particular class of periodic structures given by the crystallographic symmetry groups, in \Cref{sec:spacegrouppacking} we introduced a novel class of periodic packings, the CSG packings, by restricting possible packing configurations to a CSG and formulated the densest CSG packing problem as a nonlinear bounded and constrained optimization problem. 
	
	Our motivation was to develop a search method for the densest packing for 2D and 3D CSP that is robust to a given geometric representation of a molecule. Moreover, the method needed to be agnostic to the search configuration space and the objective function properties. For example, in our experimental setting of octagon $p2$-packings, the non-overlap constraint incorporated into the penalty function is a continuous but not differentiable function, which renders the new objective function also not differentiable. In this manner, \Cref{sec:entTrust} restated the densest packing of polytopes via stochastic relaxation \cite{geman1984} and formulated a non-Euclidean trust region method that solves this problem approximately. The resulting entropic trust region method performs updates along the geodesics on a statistical manifold where the trust region is given by KLD in a fashion similar to natural evolution strategies \cite{wierstra} and can be seen as an instance of the information geometry optimization framework \cite{ollivier}.
	
	The CSG restriction induces a toroidal topology on the packing configuration space. Therefore we perform the entropic trust region updates on a statistical manifold $S$ consisting of a parametric family of probability distributions on an $n$D unit flat torus by extending the parameter space of multivariate von Mises distributions \cite{mardia2008}, introduced in \Cref{sec:directionalstatistics}. Moreover, the exponential family reparametrization of the extended multivariate von Mises distribution equips $S$ with a dually flat structure \cite{amari2000}.
	
	Inspired by the simulated annealing control parameter \cite{laarhoven}, \Cref{sec:adaptiveselectionquantile} introduced an adaptive quantile rewriting of the fitness into the entropic trust region update schedule to facilitate the search strategy. Consequently, the natural gradient of the adaptive selection quantile-based expected fitness points in the direction of the $(q-1)$-th $q$ quantile of the fitness transformed random vector, serving as an adaptive step length method.
	
	The natural gradient descent \cite{amari} and the generalized proximal minimization algorithm \cite{censor1992} share a common ground due to the Bregman divergence characterized by the exponential family log-partition function discussed in \Cref{sec:subgradient}. In \Cref{sec:geometryEntTrust} we used this knowledge together with the dual structure given by the exponentially reparametrized extended multivariate von Mises statistical model and examined the geometry of the adaptive selection quantile equipped trust region. Embedding the statistical model $S$ into a statistical model consisting of probability distributions derived from $S$ by truncating $S$ at the $(q-1)$-th $q$-quantile of the fitness, where for every fixed $q$, $S$ becomes a submanifold of codimension $1$, we show that the resulting dual geodesic flow induced by the trust region search directions performs minimax of KLD between two hypersurfaces, one given by the statistical model $S$ and the other by $S$ truncated at the $(q-1)$-th $q$-quantile of the fitness. 
	Moreover, this minimax maximizes the stochastic dependence between the elements of the extended multivariate von Mises distributed random vector, measured by multi-information \cite{Studeny1998} or total correlation \cite{Watanabe}, providing the entropic trust region with even greater model interpretability.
	
	\begin{table}[!t]
		\centering
		\begin{tabular}{|c|c|c|c|}
			\hline
			& $\rho \left(\mathcal{K}_{G} \right)$  & $\text{dist} \left(\mathcal{K}_{G} \right)$  & $\Delta$  \\
			\hline
			regular octagon in $p2$	& 0.90616363432568 & $3.1481 \ 10^{-8}$ & $4.4318 \ 10^{-8}$
			\\
			\hline
			regular pentagon in $pg$ 	& 0.92131060131385 & $2.1811 \ 10^{-8}$ & $7.2852 \ 10^{-8}$ \\
			\hline
			regular heptagon in $p2gg$	& 0.89269066997639 & $8.2820 \  10^{-10}$ & $1.6150 \ 10^{-8}$  \\	
			\hline
			irregular pentagon in $p4$ & 0.99999999503997 & $1.5809 \ 10^{-9}$  & $4.9600 \ 10^{-9}$  \\
			\hline
			regular hexagon in $p3$	& 0.99999993380570  & $1.9423 \ 10^{-9}$ & $6.6194 \ 10^{-8}$ \\
			\hline
			$30-60-90$ triangle in $p6mm$ & 0.99999999871467 & $4.4161 \ 10^{-11}$  & $1.2853 \ 10^{-9}$  \\
			\hline
		\end{tabular}
		\caption{Packing density $\rho \left(\mathcal{K}_{G} \right)$, the minimal distance between a pair of polygons in a packing $\text{dist} \left(\mathcal{K}_{G} \right)$ and packing density difference from the theoretical optimum $\Delta$ for the presented test cases.}
		\label{tab:tabel}
	\end{table}
	
	Applying the proposed algorithm to the densest $p2$-packing of regular octagons, presented in \Cref{sec:proofofconcept}, showed competitive performance, even considering the relatively low number of samples used in the Monte--Carlo estimates, compared to the order of the exponential family rewriting of the extended multivariate von Mises distribution or dimensionality of the statistical model $S$, and the multi-modality of the optimization landscape. Furthermore, through the refining solution process, the algorithm achieved high accuracy measured by differences to the known theoretical optima. Moreover, the output configuration \cref{fig:P2Gon8Last} shows higher symmetry of the densest regular octagon packing than that of a lattice packing. The algorithm performed equally well when applied to the densest packings of regular and irregular convex polygons for which theoretical optimal solutions are known \cref{tab:tabel}. In all cases, the difference from the theoretical optimal solutions is lower than $10^{-7}$ and potentially could yield better approximations provided the refining process is allowed to run longer.
	
	Although we chose plane group packings for demonstrating the behavior and performance of the entropic trust region, the optimization algorithm is constructed to search CSGs of arbitrary dimensions. For example, in the setting of the densest space group packing of a convex polyhedron for the triclinic crystal system, the configuration space constitutes a $12$D torus given by three fractional coordinates of the polyhedron centroid, three angles of rotation of the polyhedron around the respective axes, three lengths of primitive cell edges and three angles between primitive cell edges. Moreover,  the algorithm is implemented modularly, with objective function and constraints user-specifiable as inputs which render the algorithm applicable to any nonlinear bounded constrained optimization.
	
	However, there are a few caveats regarding higher dimensional packing. The main computational bottleneck is the extended multivariate von Mises distribution Gibbs sampler (\Cref{sec:GibbsSampler}) which needs to scale better. By increasing the CSG dimension, the dimensionality of the configuration space rises polynomially, and more efficient sampling methods are necessary. Since the standard multivariate von Mises model is a stationary distribution of a Langevin diffusion stochastic differential equation \cite{garcia2019}, a natural approach is to construct a Metropolis-adjusted Langevin algorithm \cite{girolami2011}. Additionally, the stabilization of the Fisher metric tensor estimate is another caveat related to the increased dimensionality of the higher dimensional CSG packing problem. In the current implementation, a spectral radius scaling of the Fisher matrix (\Cref{sec:spectralScaling}) is used to improve the stability of the dynamical system underlying the Entropic trust region packing algorithm. A strategy to further improve the stability and reduce the number of samples necessary for accurate estimation is to regard the scaled Fisher metric tensor as a diffusion matrix \cite{fioresi2020} and derive conditions when the resulting Riemannian gradient defines a contraction mapping.
	The second most significant computational bottleneck is the overlap constraint violation computation (\Cref{sec:nonoverlap}) for a given CSG configuration. Here efficiency can be likewise improved by various heuristics. For instance, for two polytopes in a CSG configuration, it is not necessary to compute the degree of overlap if they do not intersect. Further, if the polytope circumspheres do not intersect, the polytopes do not intersect. Since the computation of sphere overlap is just one operation, the overlap constraint violation can be significantly improved for cases where polytope circumspheres do not intersect, compared to the full search for the separating hyperplane.
	
	Subsequent work is to incorporate the presented search method into existing CSP solvers to guide CSP tasks. This requires assigning a geometric representation to a molecule that can be done either manually by examining intrinsic topological properties given by the chemical composition of a molecule \cite{Santolini, wang2021} or automatically by, for example, taking the convex hull of the point set generated by the atomic coordinates of each atom as we demonstrated in \Cref{sec:pentacene}. The situation is more complicated in the $2$D CSP case since the molecule is usually defined as being embedded in $3$D Euclidean space. Thus to construct a polygon representation of the molecule, one needs to choose a projection onto the $2$D Euclidean space.

	\section*{Acknowledgements}
	
	The authors express their gratitude to Bernd Souvignier, Viktor Zamaraev, and two anonymous referees for their insightful comments and suggestions. Their contributions greatly enhanced the presentation of this work.
	
	\appendix
	
	\section{Estimating natural gradients}
	\label{subsec:montecarlo}
	Generally, the explicit computation of integrals for the natural gradient $\widetilde{\nabla} J(\bm{\theta})$ in \cref{eq:updateIGO} is not possible.  For example,  in our case the normalizer $Z(\bm{\mu}, \bm{\kappa} , \textbf{D})$ in \cref{eq:fmvm} is unknown.	The standard workaround in these situations is to use Monte--Carlo methods.
	
	For exponential families, the situation becomes simpler. The Fisher information matrix $\mathcal{I}$ \cref{eq:FIM} equals the variance of sufficient statistic \textbf{t}  and then the estimate takes the following form
	\begin{equation}\label{eq:expF}
		\widehat{\mathcal{I}_{\bm{\theta}}}=\widehat{\text{VAR}}\left[\textbf{t(\textbf{x})}|\bm{\theta}\right],
	\end{equation}
	where $\widehat{\text{VAR}}$ denotes the sample covariance matrix.
	
	The $\bm{\theta}$ gradient of $J(\bm{\theta})$ in the case of exponential families is obtained by differentiating \cref{eq:efit} as
	\begin{equation*}
		\nabla_{\bm{\theta}} J(\bm{\theta})=\int_{\textbf{X}} \textbf{F(x)}\left[ \textbf{t(\textbf{x})}-\nabla_{\bm{\theta}}\psi\left(\bm{\theta}\right)\right] dP(\bm{\theta}),
	\end{equation*}
	which is the expected value of $\textbf{F(x)}\left[ \textbf{t(\textbf{x})}-\nabla_{\bm{\theta}}\psi\left(\bm{\theta}\right)\right]$ with respect to the probability distribution $P(\bm{\theta})$. The standard gradient of the free energy $\psi\left(\bm{\theta}\right)$ equals the expected value of the sufficient statistic $\nabla_{\bm{\theta}} \psi\left(\bm{\theta}\right)=\int_{\textbf{X}}  \textbf{t(\textbf{x})}dP(\bm{\theta})$ and the estimate of the standard gradient of the expected fitness then takes the form
	\begin{equation}\label{eq:expgrad}
		\widehat{\nabla_{\bm{\theta}} J(\bm{\theta})}= \widehat{E}\left[\textbf{F(x)}\left(\textbf{t(\textbf{x})}-\widehat{E}\left[\textbf{t(\textbf{x})}\right]\right)|\; \bm{\theta}\right],
	\end{equation}
	where $\widehat{E}$ denotes the sample mean. The expression \cref{eq:expgrad} can be simplified further and receives a clear geometric interpretation using the selection quantile introduced in \cref{sec:adaptiveselectionquantile}.
	
	By combining \cref{eq:expF} and \cref{eq:expgrad}, the natural gradient estimate takes the form 
	\begin{equation}\label{eq:natGradEst}
		\widehat{\widetilde{\nabla_{\bm{\theta}}} J(\bm{\theta})}=\widehat{\mathcal{I}_{\bm{\theta}}}^{-1}\widehat{\nabla_{\bm{\theta}} J(\bm{\theta})}.
	\end{equation}	
	The above expression is also used to estimate the natural gradient for exponential families in \cite{malago} and \cite{ollivier}.
	
	\section{Toroidal distributions}
	\label{sec:toroidalDist}
	
	The general bivariate von Mises distribution density function has the following form \cite{supmardia1975}
	\begin{multline}\label{eq:gvm}
		f(\theta_1,\theta_2|\kappa_1 , \kappa_2 , \mu_1, \mu_2 ,\mathbf{A}) =  
		\frac{1}{ Z(\kappa_1 , \kappa_2 , \mu_1, \mu_2 ,\mathbf{A})}  \exp \left \{ \kappa_1 \cos(\theta_1-\mu_1) + \right.  \\ \left. + \kappa_2 \cos(\theta_2-\mu_2) +
		\left[\begin{matrix}
			\cos(\theta_1 - \mu_1) \\ 
			\sin(\theta_1 - \mu_1)
		\end{matrix}\right]^\mathtt{T}\mathbf{ A} \left[\begin{matrix}
			\cos(\theta_2 - \mu_2) \\ 
			\sin(\theta_2 - \mu_2)
		\end{matrix}\right] \right \},
	\end{multline}
	where $0\leq \theta_1,\theta_2 < 2\pi$ represent angles of corresponding unit circles $S^1$ of the product space $T^2=S^1 \times S^1$, $Z(\kappa_1 , \kappa_2 , \mu_1, \mu_2 ,\mathbf{A})$ is the normalizer, $\mu_1, \mu_2$ are the mean direction parameters, $\kappa_1 , \kappa_2$ are concentration parameters and $\mathbf{A}$ is a $2 \times 2$ matrix representing dependencies between angles $\theta_1,\theta_2$.
	
	Based on the sine submodel of the general bivariate von Mises model \cref{eq:gvm}, \cite{mardia2008} defined a probability distribution on an nD torus $T^n$ with probability density
	\begin{equation}\label{eq:mvm}
		f(\bm{\theta}|\bm{\mu}, \bm{\kappa} , \mathbf{\Lambda}) = \frac{1}{ Z(\bm{\mu}, \bm{\kappa} , \mathbf{\Lambda})} \exp \left \{ \bm{\kappa}^{\mathtt{T}}c(\bm{\theta}-\bm{\mu})+\frac{1}{2} s(\bm{\theta}-\bm{\mu})^{\mathtt{T}}\mathbf{\Lambda}s(\bm{\theta}-\bm{\mu}) \right \}
	\end{equation}
	where
	\begin{gather*}
		c(\bm{\theta}-\bm{\mu})=\left[cos(\theta_1-\mu_1),\ldots,cos(\theta_n-\mu_n)\right]^{\mathtt{T}},\\ s(\bm{\theta}-\bm{\mu})=\left[sin(\theta_1-\mu_1),\ldots,sin(\theta_n-\mu_n)\right]^{\mathtt{T}}, \\
		-\pi\leq\theta_i\leq\pi, \  -\pi\leq\mu_i\leq\pi,\  0\leq\kappa_i,\  -\infty\leq\lambda_{ij}\leq\infty, \\
		\mathbf{\Lambda}_{ij} = \lambda_{ij}=\lambda_{ji}, \ i\neq j, \ , \lambda{ii}=0. 
	\end{gather*}
	for $i,j=1,\ldots,n$.
	
	Following the full bivariate von Mises model \cref{eq:gvm}, we extend the multivariate von Mises model \cref{eq:mvm} to the full model with probability density
	\begin{equation}\label{eq:fmvm1}
		f(\bm{\theta}|\bm{\mu}, \bm{\kappa} , \textbf{D}) = \frac{1}{ Z(\bm{\mu}, \bm{\kappa} , \textbf{D})} \exp \left \{ \bm{\kappa}^{\mathtt{T}}c(\bm{\theta}-\bm{\mu})+ \left[\begin{matrix}
			c(\bm{\theta}-\bm{\mu}) \\ 
			s(\bm{\theta}-\bm{\mu})
		\end{matrix}\right]^\mathtt{T}\textbf{D} \left[\begin{matrix}
			c(\bm{\theta}-\bm{\mu}) \\ 
			s(\bm{\theta}-\bm{\mu})
		\end{matrix}\right] \right \}
	\end{equation}
	where \textbf{D} is $2n \times 2n$ matrix with no restrictions whatsoever. 
	
	Although there are no restrictions on $\textbf{A}$ in  \cref{eq:gvm} only a few specific submodels are considered in actual applications \cite{jammalamadaka2001} due to the difficulty in the statistical interpretation of the distribution parameters and a more direct relationship with the bivariate normal distribution. The same applies to the multivariate von Mises distribution \cref{eq:mvm} where the matrix $\text{diag}\left(\bm{\kappa}\right)-\bm{\Lambda}$ can be interpreted as the precision matrix of a multivariate normal distribution. However, our application of the extended multivariate von Mises model \cref{eq:fmvm1} is not dependent upon statistical interpretations of the parameters but is related to the interpretation of the probabilistic trust region method \cref{eq:Jmax}. The multivariate von Mises model \cref{eq:mvm} is a submodel of the extended multivariate model \cref{eq:fmvm1}, and the additional degrees of freedom enable the distribution to better approximate the optimization landscape induced by the expected fitness \cref{eq:efit}.

	\subsection{Exponential reformulation of the extended multivariate von Mises distribution}
	\label{subsec:exponentialfamilies}
	By restricting $\bm{\kappa}>\textbf{0}$ and using trigonometric identities we expand and rewrite \cref{eq:fmvm1} to	
	\begin{multline*} 
		f(\bm{\theta}|\bm{\mu}, \bm{\kappa} , \textbf{D}) 
		= \frac{1}{ Z(\bm{\mu}, \bm{\kappa} , \textbf{D})} \exp\left \{\sum_{i=1}^{n}\kappa_{j} \cos(\theta_i-\mu_i) +  \right.  \\
		+  \sum_{i=1}^{n} \sum_{j=1}^{n} d^{cc}_{ij}\cos(\theta_i-\mu_i)\cos(\theta_j-\mu_j)
		+   \sum_{i=1}^{n} \sum_{j=1}^{n} d^{cs}_{ij}\cos(\theta_i-\mu_i)\sin(\theta_j-\mu_j)+ \\
		+  \sum_{i=1}^{n} \sum_{j=1}^{n} d^{sc}_{ij}\sin(\theta_i-\mu_i)\cos(\theta_j-\mu_j)+  \sum_{i=1}^{n} \sum_{j=1}^{n} d^{ss}_{ij}\sin(\theta_i-\mu_i)\sin(\theta_j-\mu_j) \Biggl\} = 
	\end{multline*}	
	\begin{multline*}
		=\frac{1}{ Z(\bm{\mu}, \bm{\kappa} , \textbf{D})} \exp\left \{\sum_{i=1}^{n}\kappa_{j} \cos(\theta_i-\mu_i) +  \right. \\
		+  \sum_{i=1}^{n} \sum_{j=1}^{n} \cos(\theta_i)\cos(\theta_j)\left [d^{cc}_{ij}\cos(\mu_i)\cos(\mu_j) -  d^{cs}_{ij}\cos(\mu_i)\sin(\mu_j) - \right. \\
		\left. - d^{sc}_{ij}\sin(\mu_i)\cos(\mu_j) + d^{ss}_{ij}\sin(\mu_i)\sin(\mu_j) \right] + \\
		+  \sum_{i=1}^{n} \sum_{j=1}^{n} \cos(\theta_i)\sin(\theta_j)\left [d^{cc}_{ij}\sin(\mu_i)\sin(\mu_j) +  d^{cs}_{ij}\cos(\mu_i)\cos(\mu_j) -\right. \\
		\left. - d^{sc}_{ij}\sin(\mu_i)\sin(\mu_j) - d^{ss}_{ij}\sin(\mu_i)\cos(\mu_j) \right] + \\
		+  \sum_{i=1}^{n} \sum_{j=1}^{n} \sin(\theta_i)\cos(\theta_j)\left [d^{cc}_{ij}\sin(\mu_i)\cos(\mu_j) -  d^{cs}_{ij}\sin(\mu_i)\sin(\mu_j) + \right. \\
		\left. + d^{sc}_{ij}\cos(\mu_i)\cos(\mu_j) - d^{ss}_{ij}\cos(\mu_i)\sin(\mu_j) \right] + \\
		+  \sum_{i=1}^{n} \sum_{j=1}^{n} \sin(\theta_i)\sin(\theta_j)\left [d^{cc}_{ij}\sin(\mu_i)\sin(\mu_j) +  d^{cs}_{ij}\sin(\mu_i)\cos(\mu_j) +\right. \\
		\left. + d^{sc}_{ij}\cos(\mu_i)\sin(\mu_j) + d^{ss}_{ij}\cos(\mu_i)\cos(\mu_j) \right] \Biggl\} = \\
	\end{multline*}
	\begin{equation}\label{eq:ted}
		=\frac{1}{ Z(\bm{\mu}, \bm{\kappa} , \textbf{E})} \exp\left \{
		\bm{\kappa}^\intercal c(\bm{\theta}-\bm{\mu})+\left[\begin{matrix}
			c(\bm{\theta}) \\ 
			s(\bm{\theta})
		\end{matrix}\right]^\intercal \textbf{E} \left[\begin{matrix}
			c(\bm{\theta}) \\ 
			s(\bm{\theta})
		\end{matrix}\right] \right \}
	\end{equation}
	with	
	\begin{equation}
		\label{eq:Epar1} \textbf{E}=\left [\begin{matrix}
			\textbf{E}^{cc} & \textbf{E}^{cs} \\ 
			\textbf{E}^{sc} & \textbf{E}^{ss}
		\end{matrix} \right],
	\end{equation}
	and	
	\begin{align*}
		\textbf{E}^{cc}=\textbf{D}^{cc}\odot c(\bm{\mu})c(\bm{\mu})^{\intercal} -\textbf{D}^{cs}\odot c(\bm{\mu})s(\bm{\mu}&)^{\intercal} \\ -& \textbf{D}^{sc}\odot s(\bm{\mu})c(\bm{\mu})^{\intercal} + \textbf{D}^{ss}\odot s(\bm{\mu})s(\bm{\mu})^{\intercal}, \\
		\textbf{E}^{cs}=\textbf{D}^{cc}\odot c(\bm{\mu})s(\bm{\mu})^{\intercal} +\textbf{D}^{cs}\odot c(\bm{\mu})c(\bm{\mu}&)^{\intercal} \\ -& \textbf{D}^{sc}\odot s(\bm{\mu})s(\bm{\mu})^{\intercal} - \textbf{D}^{ss}\odot s(\bm{\mu})c(\bm{\mu})^{\intercal}, \\
		\textbf{E}^{sc}=\textbf{D}^{cc}\odot s(\bm{\mu})c(\bm{\mu})^{\intercal} -\textbf{D}^{cs}\odot s(\bm{\mu})s(\bm{\mu}&)^{\intercal} \\ +& \textbf{D}^{sc}\odot c(\bm{\mu})c(\bm{\mu})^{\intercal} - \textbf{D}^{ss}\odot c(\bm{\mu})s(\bm{\mu})^{\intercal}, 
	\end{align*}
	\begin{align*}
		\textbf{E}^{ss}=\textbf{D}^{cc}\odot s(\bm{\mu})s(\bm{\mu})^{\intercal} +\textbf{D}^{cs}\odot s(\bm{\mu})c(\bm{\mu}&)^{\intercal} \\ +& \textbf{D}^{sc}\odot c(\bm{\mu})s(\bm{\mu})^{\intercal} + \textbf{D}^{ss}\odot c(\bm{\mu})c(\bm{\mu})^{\intercal},
	\end{align*}
	where $\odot$ denotes Hadamard product and \textbf{D}
	is the interaction matrix in the extended multivariate von Mises model \cref{eq:fmvm1} with its submatrices 
	\begin{equation}\label{eq:originallD1}
		\textbf{D}=\left [\begin{matrix}
			\textbf{D}^{cc} & \textbf{D}^{cs} \\ 
			\textbf{D}^{sc} & \textbf{D}^{ss}
		\end{matrix} \right] .
	\end{equation}
	
	Further expanding the $\bm{\kappa}^\intercal c(\bm{\theta}-\bm{\mu})$ term and rewriting the $\left[\begin{matrix}
		c(\bm{\theta}) \\ 
		s(\bm{\theta})
	\end{matrix}\right]^\intercal \textbf{E} \left[\begin{matrix}
		c(\bm{\theta}) \\ 
		s(\bm{\theta})
	\end{matrix}\right]$ term in \cref{eq:ted} via
	\begin{align*}
		\left[\begin{matrix}
			c(\bm{\theta}) \\ 
			s(\bm{\theta})
		\end{matrix}\right]^\mathtt{T} \mathbf{E} \left[\begin{matrix}
			c(\bm{\theta}) \\ 
			s(\bm{\theta})
		\end{matrix}\right] &=\text{tr}\left( \left[\begin{matrix}
			c(\bm{\theta}) \\ 
			s(\bm{\theta})
		\end{matrix}\right]^\mathtt{T} \mathbf{E} \left[\begin{matrix}
			c(\bm{\theta}) \\ 
			s(\bm{\theta})
		\end{matrix}\right] \right)= 	\\
		&=\text{tr}\left( \left[\begin{matrix}
			c(\bm{\theta}) \\ 
			s(\bm{\theta})
		\end{matrix}\right] \left[\begin{matrix}
			c(\bm{\theta}) \\ 
			s(\bm{\theta})
		\end{matrix}\right]^\mathtt{T} \mathbf{E}  \right) = \text{vec}\left( \left[\begin{matrix}
			c(\bm{\theta}) \\ 
			s(\bm{\theta})
		\end{matrix}\right] \left[\begin{matrix}
			c(\bm{\theta}) \\ 
			s(\bm{\theta})
		\end{matrix}\right]^\mathtt{T}\right)^\mathtt{T} \text{vec} \left(  \mathbf{E}  \right)
	\end{align*}
	we express \cref{eq:fmvm1} in terms of the natural exponential parametrization
	\begin{equation}\label{eq:eMVM}
		f(\bm{\theta}|\bm{\eta}, \textbf{E}) = \exp\left \{ \left[\begin{matrix}
			c(\bm{\theta}) \\ 
			s(\bm{\theta})
		\end{matrix}\right]^\intercal
		\bm{\eta}
		+ \text{vec}\left( \left[\begin{matrix}
			c(\bm{\theta}) \\ 
			s(\bm{\theta})
		\end{matrix}\right] \left[\begin{matrix}
			c(\bm{\theta}) \\ 
			s(\bm{\theta})
		\end{matrix}\right]^\intercal\right)^\intercal \text{vec} \left(  \textbf{E}  \right) - \psi\left(\bm{\eta} , \textbf{E} \right) \right \}
	\end{equation}	
	where
	\begin{equation}
		\label{eq:etamap1} \bm{\eta}=\left[\begin{matrix}
			\bm{\kappa}\odot c(\bm{\mu}) \\
			\bm{\kappa}\odot s(\bm{\mu})
		\end{matrix}\right],
	\end{equation}
	$\text{tr}(\cdot)$ and $\text{vec}(\cdot)$ denote trace and vectorization of a matrix respectively, $\psi\left(\cdot,\cdot \right)$ is the logarithm of the normalizing constant or $\log$-partition function, and the canonical exponential family parameters $(\bm{\eta},\textbf{E})$ are given by \eqref{eq:etamap1} and \eqref{eq:Epar1}.
	
	Clearly, the exponential rewriting of the extended multivariate von Mises model is over-parametrized.  In order to guarantee that the Fisher information $\mathcal{I}$ \cref{eq:FIM} is positive definite, the canonical exponential family parameters of the transformed concentration, circular mean and interaction parameters need to be affinely independent since the Fisher information matrix is a  Gram matrix of log-likelihood differentials of \eqref{eq:etamap1} and \eqref{eq:Epar1} with respect to the inner product given by \eqref{eq:eMVM}. Thus, we can reparametrize the model \cref{eq:eMVM} to the minimal canonical form using the observation that
	\begin{gather}\label{eq:restrict}
		\\
		e_{ij}^{cc}\cos(\theta_{i})\cos(\theta_{j})+e_{ji}^{cc}\cos(\theta_{j})\cos(\theta_{i})=(e_{ij}^{cc}+e_{ji}^{cc})\cos(\theta_{i})\cos(\theta_{j}), \nonumber \\
		e_{ij}^{sc}\sin(\theta_{i})\cos(\theta_{j})+e_{ji}^{cs}\cos(\theta_{j})\sin(\theta_{i})=(e_{ij}^{sc}+d_{ji}^{cs})\sin(\theta_{i})\cos(\theta_{j}), \nonumber \\
		e_{ij}^{ss}\sin(\theta_{i})\sin(\theta_{j})+e_{ji}^{cc}\sin(\theta_{j})\sin(\theta_{i})=(e_{ij}^{ss}+e_{ji}^{ss})\sin(\theta_{i})\sin(\theta_{j}), \nonumber \\
		e_{ii}^{cc}(\cos(\theta_{i}))^2 + e_{ii}^{ss}(\sin(\theta_{i}))^2=e_{ii}^{ss} + (e_{ii}^{cc} -e_{ii}^{ss})(\cos(\theta_{i}))^2, \nonumber
	\end{gather}
	for $i,j=1,\ldots,n$, where $e_{ij}^{cc}, \ e_{ij}^{ss}, \ e_{ij}^{cs}, \ e_{ij}^{sc}$ are elements of $\textbf{E}^{cc}, \ \textbf{E}^{ss}, \ \textbf{E}^{cs}, \ \textbf{E}^{sc}$ submatrices of $\textbf{E}$ in \cref{eq:Epar1}. Based on the introduced reparametrization \cref{eq:restrict}, $\textbf{E}$ becomes symmetric and $e_{ii}^{cc} = - e_{ii}^{ss}$.
	
	As due to the symmetry of $\textbf{E}$ the non--diagonal parameters are counted twice, we  remove this redundancy by scaling $\textbf{D}$ in \cref{eq:fmvm1} by the factor of $2$. Now the transformation \cref{eq:Epar1} becomes
	\begin{multline}\label{eq:Epar3}
		\textbf{E}^{cc}=\tfrac{1}{2}\textbf{D}^{cc}\odot c(\bm{\mu})c(\bm{\mu})^\intercal -\tfrac{1}{2}\textbf{D}^{cs}\odot c(\bm{\mu})s(\bm{\mu})^\intercal + \tfrac{1}{2}\textbf{D}^{ss}\odot s(\bm{\mu})s(\bm{\mu})^\intercal - \tfrac{1}{2}\left.\textbf{D}^{cs}\right.^\intercal\odot s(\bm{\mu})c(\bm{\mu})^\intercal \\
		\textbf{E}^{cs}=\tfrac{1}{2}\textbf{D}^{cs}\odot c(\bm{\mu})c(\bm{\mu})^\intercal+\tfrac{1}{2}\textbf{D}^{cc}\odot c(\bm{\mu})s(\bm{\mu})^\intercal - \tfrac{1}{2}\textbf{D}^{ss}\odot s(\bm{\mu})c(\bm{\mu})^\intercal - \tfrac{1}{2}\left.\textbf{D}^{cs}\right.^\intercal\odot s(\bm{\mu})s(\bm{\mu})^\intercal \\
		\textbf{E}^{ss}=\tfrac{1}{2}\textbf{D}^{ss}\odot c(\bm{\mu})c(\bm{\mu})^\intercal+\tfrac{1}{2}\textbf{D}^{cc}\odot s(\bm{\mu})s(\bm{\mu})^\intercal +\tfrac{1}{2}\textbf{D}^{cs}\odot s(\bm{\mu})c(\bm{\mu})^\intercal + \tfrac{1}{2}\left.\textbf{D}^{cs}\right.^\intercal\odot c(\bm{\mu})s(\bm{\mu})^\intercal
	\end{multline}
	and as a consequence, the extended multivariate von Mises model \cref{eq:fmvm1} now becomes
	\begin{equation}\label{eq:fmvm2}
		f(\bm{\theta}|\bm{\mu}, \bm{\kappa} , \textbf{D}) = \frac{1}{ Z(\bm{\mu}, \bm{\kappa} , \textbf{D})} \exp \left \{ \bm{\kappa}^{\mathtt{T}}c(\bm{\theta}-\bm{\mu})+ \frac{1}{2} \left[\begin{matrix}
			c(\bm{\theta}-\bm{\mu}) \\ 
			s(\bm{\theta}-\bm{\mu})
		\end{matrix}\right]^\mathtt{T}\textbf{D} \left[\begin{matrix}
			c(\bm{\theta}-\bm{\mu}) \\ 
			s(\bm{\theta}-\bm{\mu})
		\end{matrix}\right] \right \},
	\end{equation}
	where $\textbf{D}$ is $2n \times 2n$ matrix with the same structure as $\textbf{E}$ in \cref{eq:Epar1}, that is $\textbf{D}$ is symmetric with the diagonal elements $d_{ii}^{cc}, \ d_{ii}^{ss}$ of the $\textbf{D}^{cc}, \ \textbf{D}^{ss}$ submatrices of $\textbf{D}$ in \cref{eq:originallD1} related by $d_{ii}^{cc} = - d_{ii}^{ss}$.
	
	The inverse transformations to the extended multivariate von Mises model concentration, circular mean and interaction matrix parametrization can be obtained as the solution to the system of equations \cref{eq:etamap1} and \cref{eq:Epar3} for $(\bm{\mu}, \bm{\kappa} , \textbf{D})$ in terms of $(\bm{\eta}, \textbf{E})$ in the following form
	\begin{gather*}
		\mu_{i} = 2 \arctan\left(\frac{\eta_{i}^{s}}{\eta_{i}^{c}+\sqrt{{\eta_{i}^{c}}^2+{\eta_{i}^{s}}^2}}\right) \mod 2\pi, \\
		\kappa_{i} = \sqrt{{\eta_{i}^{c}}^2+{\eta_{i}^{s}}^2}
	\end{gather*}
	where $\eta_{i}^{c}$ and $\eta_{i}^{s}$ are the exponential canonical parameters composing $\bm{\eta}$ in \cref{eq:etamap1} for  $i=1,\ldots,n$ associated with the cosine and sine respectively, and
	\begin{subequations}
		\begin{equation*}
			\textbf{D}^{cc}=2\textbf{E}^{cc}\odot c(\bm{\mu})c(\bm{\mu})^{\mathtt{T}} +2\textbf{E}^{cs}\odot c(\bm{\mu})s(\bm{\mu})^{\mathtt{T}} + 2\left.\textbf{E}^{cs}\right.^\intercal\odot s(\bm{\mu})c(\bm{\mu})^{\mathtt{T}} + 2\textbf{E}^{ss}\odot s(\bm{\mu})s(\bm{\mu})^{\mathtt{T}}
		\end{equation*}
		\begin{equation*}
			\textbf{D}^{cs}=2\textbf{E}^{cs}\odot c(\bm{\mu})c(\bm{\mu})^{\mathtt{T}} -2\textbf{E}^{cc}\odot c(\bm{\mu})s(\bm{\mu})^{\mathtt{T}} +2\textbf{E}^{ss}\odot s(\bm{\mu})c(\bm{\mu})^{\mathtt{T}} - 2\left.\textbf{E}^{cs}\right.^\intercal\odot s(\bm{\mu})s(\bm{\mu})^{\mathtt{T}}
		\end{equation*}
		\begin{equation*}
			\textbf{D}^{ss}=2\textbf{E}^{ss}\odot c(\bm{\mu})c(\bm{\mu})^{\mathtt{T}} -2\textbf{E}^{cs}\odot s(\bm{\mu})c(\bm{\mu})^{\mathtt{T}} -2\left.\textbf{E}^{cs}\right.^\intercal\odot c(\bm{\mu})s(\bm{\mu})^{\mathtt{T}} + 2\textbf{E}^{cc}\odot s(\bm{\mu})s(\bm{\mu})^{\mathtt{T}}.
		\end{equation*}
	\end{subequations}
	
	The multivariate von Mises model \cref{eq:mvm} and the extended multivariate von Mises model \cref{eq:fmvm1} coincide when the precision matrix $\textbf{D}$ is of the form
	\begin{equation*}
		\textbf{D}=\left [\begin{matrix}
			0 & 0 \\ 
			0 & \bm{\Lambda}
		\end{matrix} \right]. \\
	\end{equation*}
	In terms of differential geometry, the parametrizations constitute coordinate systems on the manifolds of probability measures. In this regard, the relationship between the dimensionality of the extended multivariate von Mises statistical model $S_{\text{eMvM}}$ given by \cref{eq:fmvm2} and the multivariate von Mises distribution $S_{\text{MvM}}$ given by \cref{eq:mvm} is
	\begin{equation*}
		\dim\left(S_{\text{fMvM}}\right)=2n(n+1)>\frac{n(n+3)}{2}=\dim\left(S_{\text{MvM}}\right),
	\end{equation*}
	and as a consequence, $S_{\text{MvM}} \subset S_{\text{eMvM}}$ for $\bm{\kappa}>0$, which implies that the multivariate von Mises statistical model is a submanifold of the extended multivariate von Mises model.
	
	\subsection{The extended multivariate von Mises submodel}\label{fmvMsub}
	To reduce the computational burden of sampling from a $2n\left(n+1\right)$ dimensional statistical model, where $n$ denotes the dimensionality of the supporting torus, and to improve the stability of the Fisher metric estimate, we further reduce the extended multivariate von Mises statistical model \cref{eq:fmvm2} dimensionality by removing interactions between cosines and sines for the same angle. 
	
	Considering the exponential rewriting of the extended multivariate von Mises model \cref{eq:eMVM} as a graphical interaction model \cite{cowell2007}, this is equivalent to removing direct feedback loops from the system. Mathematically, this means setting the elements $e_{ij}^{cc}, \ e_{ij}^{ss}, \ e_{ij}^{cs}$ of $\textbf{E}^{cc}, \ \textbf{E}^{ss}, \ \textbf{E}^{cs}$ submatrices of $\textbf{E}$ in \cref{eq:Epar1} for $i=j$ equal to zero. As a consequence $\textbf{D}^{cc}, \ \textbf{D}^{ss}, \ \textbf{D}^{cs}$ submatrices of the matrix $\textbf{D}$ representing interactions in \cref{eq:fmvm2} have zero diagonals. The dimension of this specific submodel is $2n^2$.
	
	\subsection{Extended multivariate von Mises Gibbs sampler} 
	\label{sec:GibbsSampler}
	To generate realizations from the \cref{eq:eMVM} distributed random vectors for the computation of the natural gradient estimates \cref{eq:natGradEst}, we use the multi-stage Gibbs sampler \cite{robert}.  The univariate conditionals of the exponential rewriting of the extended multivariate von Mises distribution can be obtained by expanding \cref{eq:eMVM}, moving all terms that do not depend on $\theta_k$ to the normalization constant, and using trigonometric identities, to the following form 
	\begin{multline}\label{eq:ted2}
		f(\theta_k | \theta_1,\ldots,\theta_{k-1},\theta_{k+1},\ldots, \theta_n) \propto  
		\\
		\exp \left \{ \left[\eta^c_{k} + \sum\limits_{\substack{i=1 \\ i\neq k}}^n (e_{ki}^{cc}+e_{ik}^{cc})\cos(\theta_i)+ \sum\limits_{\substack{i=1 \\ i\neq k}}^n (e_{ki}^{cs}+e_{ik}^{sc})\sin(\theta_i)\right] \cos(\theta_k) + \right.
		\\
		+  \left[\eta^s_{k} + \sum\limits_{\substack{i=1 \\ i\neq k}}^n (e_{ki}^{ss}+e_{ik}^{ss})\sin(\theta_i)+ \sum\limits_{\substack{i=1 \\ i\neq k}}^n (e_{ki}^{sc}+e_{ik}^{cs})\cos(\theta_i)\right] \sin(\theta_k) + 
		\\   
		+ e_{kk}^{cc}\cos(2\theta_k) + \left[e_{kk}^{sc}+e_{kk}^{cs}\right]\sin(2\theta_k) - e_{kk}^{cc}\cos(2\theta_k) \Biggl\}.
	\end{multline}
	Further rewriting \eqref{eq:ted2} yields	
	\begin{equation}\label{eq:gvm2}
		f(\theta_k | \theta_1,\ldots,\theta_{k-1},\theta_{k+1},\ldots, \theta_n) \propto  
		\exp \left\{ \gamma_k^1 \cos(\theta_k - \nu_k^1) + \gamma_k^2\cos(2(\theta_k-\nu_k^2)) \right\},
	\end{equation}
	where
	\begin{multline*}
		\gamma_k^1 =
		\left\{\left(\eta^c_{k} + \sum\limits_{\substack{i=1 \\ i\neq k}}^n (e_{ki}^{cc}+e_{ik}^{cc})\cos(\theta_i)+ \sum\limits_{\substack{i=1 \\ i\neq k}}^n (e_{ki}^{cs}+e_{ik}^{sc})\sin(\theta_i)\right)^2 \right. +  \\
		+ \left. \left(\eta^s_k + \sum\limits_{\substack{i=1 \\ i\neq k}}^n (e_{ki}^{ss}+e_{ik}^{ss})\sin(\theta_i)+ \sum\limits_{\substack{i=1 \\ i\neq k}}^n (e_{ki}^{sc}+e_{ik}^{cs})\cos(\theta_i)\right)^2\right\}^{\frac{1}{2}}, 
	\end{multline*}	
	\begin{equation*}
		\nu_k^1 = \text{atan2}\left(\frac{\eta^s_k + \sum\limits_{\substack{i=1 \\ i\neq k}}^n (e_{ki}^{ss}+e_{ik}^{ss})\sin(\theta_i)+ \sum\limits_{\substack{i=1 \\ i\neq k}}^n (e_{ki}^{sc}+e_{ik}^{cs})\cos(\theta_i)}{\eta^c_{k} + \sum\limits_{\substack{i=1 \\ i\neq k}}^n (e_{ki}^{cc}+e_{ik}^{cc})\cos(\theta_i)+ \sum\limits_{\substack{i=1 \\ i\neq k}}^n (e_{ki}^{cs}+e_{ik}^{sc})\sin(\theta_i)}\right),
	\end{equation*}
	\begin{equation*}
		\gamma_k^2 = \left\{\left(e_{kk}^{cc}-e_{kk}^{ss}\right)^2 + (e_{kk}^{sc}+e_{kk}^{cs})^2\right\}^\frac{1}{2}, 
	\end{equation*}	
	\begin{equation*}
		\nu_k^2 = \frac{1}{2}\text{atan2}\left(\frac{e_{kk}^{sc}+e_{kk}^{cs}}{e_{kk}^{cc}-e_{kk}^{ss}}\right), \\
	\end{equation*}
	where $\text{atan2}(\frac{y}{x})$ is the argument of the complex number $x+iy$, $\eta_i^c$ and $\eta_i^s$ are the canonical parameters in $\bm{\eta}$ associated with $\kappa_i\cos(\mu_i)$ and $\kappa_i\sin(\mu_i)$ in \cref{eq:etamap1} respectively, and $e_{ij}^{cs}$,  $e_{ij}^{sc}$,  $e_{ij}^{cc}$,  $e_{ij}^{ss}$ are the elements of the $\mathbf{E}^{cs}$, $\mathbf{E}^{sc}$, $\mathbf{E}^{cc}$, $\mathbf{E}^{cc}$ submatrices of $\mathbf{E}$ in \cref{eq:Epar1}.
	
	The probability density function \cref{eq:gvm2} is of the generalized von Mises distribution introduced in \cite{gatto2007}. We implement von Neumann's rejection sampling algorithm of \cite{gatto2008} to generate samples from the generalized multivariate von Mises distribution.
	
	Implementation of the extended multivariate von Mises Gibbs sampler is presented in Algorithm~\eqref{eMvMgibbs} in pseudocode form. The Gibbs sampler's number of iterations $m$ was determined experimentally and is set to 100.
	\begin{algorithm}
		\caption{Extended multivariate von Mises Gibbs sampler}
		\label{eMvMgibbs}
		\begin{algorithmic}
			\REQUIRE{$\bm{\eta}$, n, m}
			\ENSURE{$\mathbf{U}_{d\times n}\sim \text{EMvM}_d(\bm{\eta})$}
			\FOR{$k=1$ to $d-1$}
			\STATE {$\mathbf{U}_{k,1:} \leftarrow \text{Generate n i.i.d.} \ GvM_2\left(\sqrt{(\eta_k^c)^2+(\eta_k^s)^2}, \ \text{atan2}\left(\frac{\eta_k^s}{\eta_k^c}\right), \ \gamma_k^2, \ \nu_k^2\right)$}
			\ENDFOR
			\FOR{$j=1$ to $n$}
			\STATE{$\mathbf{U}_{d,j} \leftarrow GvM_2\left(\gamma_d^1, \ \nu_d^1, \ \gamma_d^2, \ \nu_d^2 | \mathbf{U}_{1,j},\ldots,\mathbf{U}_{d-1,j}\right)$}
			\ENDFOR
			\FOR{$i=1$ to $m$}
			\FOR{$k=1$ to $d$}
			\FOR{$j=1$ to $n$}
			\STATE{$\mathbf{U}_{k,j} \leftarrow GvM_2\left(\gamma_k^1, \ \nu_k^1, \ \gamma_k^2, \ \nu_k^2\right|\mathbf{U}_{1,j},\ldots,\mathbf{U}_{k-1,j},\mathbf{U}_{k+1,j},\ldots,\mathbf{U}_{d,j})$}
			\ENDFOR
			\ENDFOR
			\ENDFOR
		\end{algorithmic}
	\end{algorithm}
	
	\section{Adaptive quantile hill climbing}
	\label{sec:quntileHillClimb}
	
	For fixed $q_t$-quantile, iteratively solving \cref{eq:Jmax} with the expected fitness of the form \cref{eq:qESJtheta} is equivalent to solving
	\begin{gather}
		\bm{\theta}^{t+1}=\argmax_{\bm{\theta} \in \bm{\Theta}} p(\mathbf{x}|\bm{\theta}), \\
		\text{s.t.} \  \textbf{x}\sim P(\bm{\theta}^{t}) \left. \right| \textbf{F}(\textbf{x}) \geq \widehat{\textbf{F}_{1- \frac{1}{q_{t}}}^{\bm{\theta}^t}} . \label{eq:likeStep}
	\end{gather}
	The update step $\bm{\theta}^{t+1}$ is given by maximizing the likelihood $p(\mathbf{x}|\bm{\theta})$ over \cref{eq:likeStep}. In this context, the adaptive selection quantile can be formulated as a more general search method, summarized in the pseudocode in Algorithm \ref{alg:adaSQalg}.
	
	\begin{algorithm}
		\caption{Adaptive quantile hill climbing}
		\label{alg:adaSQalg}
		\begin{algorithmic}
			\REQUIRE{parametric family $P(\bm{\theta})$, \textbf{F(x)}, $N$, $q_{\max}$, $\beta$}
			\ENSURE{$\textbf{x}_{opt}$}
			\REPEAT
			\STATE{$\cos(\alpha^{t-1})\leftarrow\frac{<\Delta\bm{\theta}^{t-1}, \Delta\bm{\theta^}{t-2}>_{F(\bm{\theta}^{t-2})}}{||\Delta\bm{\theta}^{t-1}||_{ F(\bm{\theta}^{t-2})}||\Delta\bm{\theta}^{t-2}||_{ F(\bm{\theta}^{t-2})}}$}
			\STATE{$q_{t}\leftarrow\min\left[q_{t-1}\exp\left\{\beta \cos(\alpha^{t-1})\right\},q_{\max}\right]$}
			\STATE{Generate N i.i.d. samples $\textbf{x}_i\sim P(\bm{\theta}^{t-1})$}
			\STATE{Assign ranks $r_j$ to $\textbf{F}(\textbf{x}_i)$}
			\STATE{Update $\bm{\theta}^{t}$ by performing a maximum likelihood estimate of $P(\bm{\theta})$ from the ranked samples $\textbf{x}_{r_{N}},\textbf{x}_{r_{N-1}},\ldots,\textbf{x}_{r_{N-\frac{N}{q_{t}}}}$}
			\UNTIL{stopping criterion is met}
		\end{algorithmic}
	\end{algorithm}
	
	\section{Implementation details}
	\label{sec:implementation}
	We present a few valid technical details and algorithmic settings here. These include the map between the unit flat n-torus and the optimization configuration space \cref{subsec:periodicity}, the penalty function used to integrate nonlinear constraints into the optimization schedule \cref{subsec:constraint}, the non--overlapping constraint violation formulation \cref{sec:nonoverlap}, the adaptive learning rate \cref{subsec:learningrates}, the Fisher metric tensor scaling used to stabilize the unit natural gradients \cref{sec:spectralScaling}, the hyperparameter tuning method we developed \cref{sec:algparameters} and finally, the acceleration of computations through parallelization \cref{sec:parallel}.

	\subsection{Boundary mapping}
	\label{subsec:periodicity}
	
	The statistical model \cref{eq:fmvm1} we are working with consists of probability distributions with its support on the $n$-torus $T^n$ whose product space components are unit circles and is topologically equivalent to the identification space \cite{willard2012}
	\begin{equation*}
		\faktor{\prod_{i=1}^{n}{[0,2\pi]}} {\sim}
	\end{equation*}
	where the equivalence relation $\sim$ is defined by identifying all points such that
	\begin{gather*}
		(0,x_2,\ldots,x_n)\sim (2\pi,x_2,\ldots,x_n), \\
		(x_1,0,\ldots,x_n) \sim (x_1,2\pi,\ldots,x_n), \\
		\vdots \\
		(x_1,x_2,\ldots,0) \sim (x_1,x_2,\ldots,2\pi).
	\end{gather*}
	We need to map the $n$-cube $[0,2\pi]^n$ to an n-orthotope defined by the boundary constraints of the optimization problem in order to evaluate the fitness of each realization of the extended multivariate von Mises distributed random vector.    Specifically, we define a map $y_i: \ [0,2\pi)\rightarrow[l_i,u_i)$ via
	\begin{gather*}
		y_i=\frac{x_i}{2\pi}(u_i - l_i) + l_i \ , \ x_i \in \left[0,2\pi \right)
	\end{gather*}
	for $i=1,\ldots,n$ where $u_i$ and $l_i$ are $i$-th upper and lower bound, respectively.
	
	This kind of mapping is natural for variables with inherent periodicity but somewhat problematic for nonperiodic ones due to the discontinuity in the configuration space introduced by identifying lower and upper bounds, and more importantly, $\textbf{F}\left(x_1,\ldots,l_a,\ldots,x_n\right) \neq \textbf{F}\left(x_1,\ldots,u_a,\ldots,x_n\right)$ in the case when there is no period $p$ such that $l_a=u_a \mod p$. 
	
	To address this inconvenience, we define a different boundary map for aperiodic variable $x_a$, given by
	\begin{equation*}
		y_a=\left \{\begin{array}{cll}
			l_a+\frac{x_a}{\pi}(u_a - l_a)           & \text{if} & x_a \in \left[0,\pi \right), \\ 
			2u_a-l_a-\frac{x_i}{\pi}(u_a - l_a) & \text{if}  &x_a \in \left[\pi,2\pi \right),
		\end{array} \right.
	\end{equation*}
	at the expense of loosing injectivity of the boundary map $y_a:[0,2\pi)\rightarrow[l_a;u_a)$ and introducing additional extrema in the optimization landscape. In practice, this is not a problem since the algorithm is built to be robust in complex optimization landscapes.
	
	After combining the above, we have the following  boundary constraint mapping
	\begin{equation*}
		y_i=\left \{\begin{array}{ccll}
			l_i + \frac{x_i}{2\pi}(u_i - l_i) & \text{if} & x_i \ \text{is periodic} & \wedge\quad   x_i \in \left[0,2\pi \right), \\
			l_i+\frac{x_i}{\pi}(u_i - l_i) & \text{if} & x_i \ \text{is aperiodic} & \wedge\quad   x_i \in \left[0,\pi \right), \\ 
			2u_i-l_i-\frac{x_i}{\pi}(u_i - l_i) & \text{if}& x_i \ \text{is aperiodic}  & \wedge\quad  x_i \in \left[\pi,2\pi \right),
		\end{array} \right. 
	\end{equation*}
	for $i=1,\ldots,n$.

	\subsection{Constraint handling}
	\label{subsec:constraint}
	To address linear and nonlinear constraints, we implement a penalty function based on feasibility considerations \cite{deb}. The basic premise is to create an ordering on the set of candidate solutions, such that: i) any feasible solution is better than any infeasible one, ii) between two feasible solutions, the one with better objective function is preferred, and iii) between two infeasible solutions, the one with lower constraint violation is preferred. Then for the following minimization problem
	\begin{align}\label{eq:constraints}
		&\min_{\mathbf{x} \in \mathbf{x}} f\left(\mathbf{x}\right), \\ \nonumber
		&\text{s.t.} \ g_j\left(\mathbf{x}\right)\leq 0, \ j=1,\ldots,J,
	\end{align}
	where $g_j\left(\mathbf{x}\right)$ are inequality constraints, given a set of solutions $\left\{\mathbf{x}_i \ | \ i=1,\ldots, N \ \right\}$, the penalty function is expressed as
	\begin{equation}\label{eq:penalty}
		\textbf{F}\left(\mathbf{x}_i \right)=\left \{\begin{array}{cl}
			f\left(\mathbf{x}_i\right)& \text{if }  \ g_j\left(\mathbf{x}_i\right)\leq 0, \forall \ j=1,\ldots,J \\ 
			f^{\text{max}}+\sum_{k \in K} \frac{g_k\left(\mathbf{x}_i\right)}{g_{k}^{\max}}& \text{otherwise }   
		\end{array} \right.
	\end{equation}
	where $K=\left\{j:g_j\left(\mathbf{x}_i\right)> 0, \ j=1,\ldots,J \right\}$ and 
	\begin{gather}
		f^{\max}=\max \left\{f\left(\mathbf{x}_i\right) \ | \ \mathbf{x}_i:\ g_j\left(\mathbf{x}_i\right)\leq 0, \ i=1,\ldots,N \ , \ j=1,\ldots,J \right\}, \\
		g_{j}^{\max}=\max \left\{g_j\left(\mathbf{x}_i\right) \ | \ \mathbf{x}_i:\ g_j\left(\mathbf{x}_i\right)> 0, \ i=1,\ldots,N \right\}.
	\end{gather}
	Here, the penalty term in \cref{eq:penalty} for an infeasible solution $\mathbf{x}_i$ is the sum of the maximum of all feasible solutions sampled at a given iteration $f^{\max}$ and the sum of constraint violations normalized by the maximal constraint violation $g_{j}^{\max}$ for each constraint $g_{j}$.
	
	Note that this is particularly well suited for the adaptive selection quantile introduced in \cref{sec:adaptiveselectionquantile} since only the ordering is considered for the trust region updates.   Additionally, the penalty function \cref{eq:penalty} can be easily augmented for multiple objectives by using some suitable aggregation function $G\left(f_1\left(\mathbf{x}\right),\ldots,f_R\left(\mathbf{x}\right)\right)$ where $f_r\left(\mathbf{x}\right)$ for $r=1,\ldots,R$ are the multiple objective functions \cite{coello2007}.

	\subsection{$G$-packing overlap constraint evaluation for convex polygons and polyhedra}
	\label{sec:nonoverlap}
	
	To evaluate the intersection and the degree of constraint violation between convex polygons and polyhedra in candidate solutions, we use the method based on Phi-functions \cite{chernovSup}. Given the convex polytope 
	\begin{equation*}
		P_0=\text{conv}\left\{V_0 | V_0=(\textbf{v}_1,\textbf{v}_2,\dots, \textbf{v}_m),\textbf{v}_i \in \mathbb{R}^n\right\},
	\end{equation*}
	centred at the origin and defined by the convex hull of the vertices $\textbf{v}_i$, and given rotated and translated copies of the reference polytope $P_0$
	\begin{gather*}
		P_1= \textbf{R}_1 P_0 + \textbf{c}_1 , \\
		P_2= \textbf{R}_2 P_0 + \textbf{c}_2 ,
	\end{gather*}	
	for some rotation matrices $\textbf{R}_1, \ \textbf{R}_2$ and translation vectors $\textbf{c}_1, \ \textbf{c}_2$, separating hyperplane theorem~\cite{boyd2004convex} states that if $P_1$ and $P_2$ do not overlap, there exists a hyperplane that separates them.
	
	For convex polytopes of dimension $n=2$, the hyperplanes defined by the edges of $P_1$ and $P_2$ are all candidate separating hyperplanes, and it is adequate to check vertices of $P_1$ against $P_2$ hyperplanes and vice versa.
	For convex polytopes of dimension $n=3$, additional possible separating hyperplanes need to be defined by combining an edge from polytope $P_1$ and an edge from polytope $P_2$, apart from the hyperplanes defined by their respective faces.
	
	To implement this, vertices of $P_1$ are express in the coordinate system of $P_2$ denoted by $P_{12}$ and vertices of $P_2$ in the coordinate system of $P_1$ denoted by $P_{21}$ as
	\begin{gather*}
		P_{12} = \textbf{R}_2^{-1}\left[P_1-\textbf{c}_2\right], \\
		P_{21} = \textbf{R}_1^{-1}\left[P_2-\textbf{c}_1\right],
	\end{gather*}
	and a collection the hyperplanes $H_1$ characterizing $P_0$ is defined. In the $3$D case, additional collection of hyperplanes $H_2$ is defined by all combinations of an edge of $P_0$ and an edge of $P_{21}$, such that the hyperplane contains the edge of $P_0$,
	where the hyperplane normal vectors are set to the unit length with the direction outwards of $P_0$.
	
	Since $\lVert \textbf{h} \rVert=1$ by inserting vertices of $P_{12}$ and $P_{21}$ into the hyperplane equations we not only check for the existence of a separating hyperplane but in practice compute the euclidean distance between $P_1$ and $P_2$ which has the following closed form expression
	\begin{equation}\label{eq:dist}
		\text{dist} \left(P_1,P_2\right) =  \max\{\text{dist}_1 \left(P_1,P_2\right), \text{dist}_2 \left(P_1,P_2\right) \},
	\end{equation}
	where $\text{dist}_{1} \left(P_1,P_2\right) =  \max\{\max_{\textbf{h} \in H_1} \min_{\textbf{v} \in P_{12}} \textbf{h}^T [\textbf{v};1 ], \max_{\textbf{h} \in H_1} \min_{\textbf{v} \in P_{21}} \textbf{h}^T [\textbf{v};1] \}$, and $\text{dist}_{2} \left(P_1,P_2\right) =  \max\{\max_{\textbf{h} \in H_2} \min_{\textbf{v} \in P_{0}} -\textbf{h}^T [\textbf{v};1 ], \max_{\textbf{h} \in H_2} \min_{\textbf{v} \in P_{21}} \textbf{h}^T [\textbf{v};1] \}$.
	
	The distance function \cref{eq:dist} is a continuous and piecewise differentiable function, and $P_1$ and $P_2$ do not intersect if and only if $\text{dist} \left(P_1,P_2\right) \geq 0$.
	
	To evaluate whether a collection of convex polytopes $\mathcal{K}$ defined as the orbit of the convex polytope $P_0$ under the action of the CSG $G$ is a $\mathcal{G}$-packing, for $\mathcal{K}_{G}$ in \cref{eq:spacepacking} we define
	\begin{equation}\label{eq:distPack}
		\text{dist}\left(\mathcal{K}_{G}\right)=\min_{i,j,\{\alpha\}}\text{dist} \left(g_iP_0,g_{j\{\alpha\}}P_0\right),
	\end{equation} 
	where $g_i, \ g_{j\{\alpha\}} \in G$ such that
	\begin{gather}\label{eq:transSym}
		g_iP_0=R_iP_0 + \mathbf{ a}_i, \\ \label{transLat}
		g_{j\{\alpha\}}P_0 = R_jP_0 + \mathbf{ a}_j + \mathbf{ l}_{\{\alpha\}},
	\end{gather}
	for $i,j=1,\ldots,N$, $\{\alpha\}=\{(u_1,\ldots,u_n)|u_i\in\{-k,\ldots-2,-1,0,1,2,\ldots,k\}\}$ and $i\neq j$ if $\{\alpha\}=(0,\ldots,0)$ where $\mathbf{ l}_{\{\alpha\}} \in L$ \cref{eq:latticeG} and $\mathbf{ a}_i$ is of the form \cref{eq:alat}.	  In other words, we compute minimal Euclidean distances between orbits of $P_0$ \cref{eq:transSym} whose centroids lie inside the primitive cell and orbits of $P_0$ \eqref{transLat} whose centroids lie inside neighbouring primitive cells to up to twice the lattice basis vectors $\mathbf{ b}_1,\ldots,\mathbf{ b}_n$.
	
	During our experiments, if the upper bound on the size of the lattice vector generators was set to the corresponding lattice vector generators of the $\mathcal{G}$-packing of the circumscribed $(n-1)$-sphere of $P_0$, then to evaluate whether $\mathcal{K}_{G}$ is a $\mathcal{G}$-packing it was usually enough to assess the intersection between $g_iP_0$ and $g_{j\{\alpha\}}P_0$ for $\mathbf{ l}_{\{\alpha\}} \in L$ up to the first primitive cell ($k=1$) in every coordinate direction, although in some instances it was necessary to increase the value of $k$. For example, in the case of $p1$-packing of the pentacene representation, introduced in \Cref{sec:pentacene}, $k$ needed to by at least three. Generally, the value of $k$ depends on the shape of given polytope.
	
	Based on \cref{eq:distPack} $\mathcal{K}_{G}$ is a CSG packing if and only if $\text{dist}\left(\mathcal{K}_{G}\right)\geq 0$. In the case of $\text{dist}\left(\mathcal{K}_{G}\right)< 0$, we get a measure of constraint violation of a candidate solution defined by $\lvert \text{dist}\left(\mathcal{K}_{G}\right) \rvert$ used in penalty function \cref{eq:penalty} by setting $g\left(\mathbf{x}\right)=-\text{dist}\left(\mathcal{C}_{G}(\mathbf{x})\right)$ where $\mathcal{C}_{G}(\mathbf{x})$ is a $G\in\mathcal{G}$ configuration of the form \cref{eq:spacepackingConfig} defined by the candidate solution $\mathbf{x}$ that is not necessarily a packing.

	\subsection{Learning rates}
	\label{subsec:learningrates}
	Experiments show that having a single trust region radius for the exponential multivariate von Mises statistical manifold is insufficient and results in poor performance. Instead, we decide to transfer the unit gradient \cref{eq:unitGrad} back to the original circular mean, concentration and angle interaction parametrizations and perform the gradient ascent updates in those coordinates, allowing us to use an additional separate learning rate for each parameter group. In \cref{subsubsec:separation}, we introduce the aforementioned change of coordinates of the natural gradients. Additionally, we modify the adaptive learning rate scheme proposed in \cite{silva} described in \cref{subsubsec:adaptivelearningrates}.

	\subsubsection{Circular mean, concentration and precision update equations}
	\label{subsubsec:separation}
	
	By differentiating \cref{eq:etamap1} and \cref{eq:Epar3} with respect to time, we get a system of following linear equations
	\begin{equation*}
		\frac{d\eta^c_{i}}{dt} = -\kappa_{i}\sin(\mu_{i})\frac{d\mu_{i}}{dt} + \cos(\mu_{i})\frac{d\kappa_{i}}{dt},
	\end{equation*}
	\begin{equation*}
		\frac{d\eta^s_{i}}{dt} =\kappa_{i}\cos(\mu_{i})\frac{d\mu_{i}}{dt}  + \sin(\mu_{i})\frac{d\kappa_{i}}{dt},
	\end{equation*}
	\begin{multline*}
		\frac{de^{cc}_{ij}}{dt} =\cos(\mu_{i})\cos(\mu_{j})\frac{dd^{cc}_{ij}}{dt} - \cos(\mu_i)\sin(\mu_j)\frac{dd^{cs}_ {ij}}{dt} - \sin(\mu_i)\cos(\mu_j)\frac{dd^{sc}_{ij}}{dt} \\ + \sin(\mu_i)\sin(\mu_j)\frac{dd^{ss}_{ij}}{dt} -  d_{ij}^{cs}\cos(\mu_i)\cos(\mu_j)\frac{d\mu_{j}}{dt} 	
		-d_{ij}^{sc}\cos(\mu_i)\cos(\mu_i)\frac{d\mu_{i}}{dt} \\ - d_{ij}^{cc}\cos(\mu_j)\sin(\mu_i)\frac{d\mu_{i}}{dt}  - d_{ij}^{cc}\cos(\mu_i)\sin(\mu_j)\frac{d\mu_{j}}{dt}
		+ d_{ij}^{ss}\cos(\mu_i)\sin(\mu_j)\frac{d\mu_{i}}{dt} \\
		+ d_{ij}^{ss}\cos(\mu_j)\sin(\mu_i)\frac{d\mu_{j}}{dt} + d_{ij}^{cs}\sin(\mu_i)\sin(\mu_j)\frac{d\mu_{i}}{dt} + d_{ij}^{sc}\sin(\mu_i)\sin(\mu_j)\frac{d\mu_{j}}{dt},
	\end{multline*}
	\begin{multline*}
		\frac{de^{cs}_{ij}}{dt} =\cos(\mu_{i})\cos(\mu_{j})\frac{dd^{cs}_{ij}}{dt} + \cos(\mu_i)\sin(\mu_j)\frac{dd^{cc}_{ij}}{dt} - \sin(\mu_i)\cos(\mu_j)\frac{dd^{ss}_{ij}}{dt} \\ - \sin(\mu_i)\sin(\mu_j)\frac{dd^{sc}_{ij}}{dt} +  d_{ij}^{cc}\cos(\mu_i)\cos(\mu_j)\frac{d\mu_{j}}{dt} 	
		-d_{ij}^{ss}\cos(\mu_i)\cos(\mu_j)\frac{d\mu_{i}}{dt} \\ - d_{ij}^{cs}\cos(\mu_j)\sin(\mu_i)\frac{d\mu_{i}}{dt}  - d_{ij}^{cs}\cos(\mu_i)\sin(\mu_j)\frac{d\mu_{j}}{dt}
		- d_{ij}^{sc}\cos(\mu_i)\sin(\mu_j)\frac{d\mu_{i}}{dt} \\
		- d_{ij}^{sc}\cos(\mu_j)\sin(\mu_i)\frac{d\mu_{j}}{dt} - d_{ij}^{cc}\sin(\mu_i)\sin(\mu_j)\frac{d\mu_{i}}{dt} + d_{ij}^{ss}\sin(\mu_i)\sin(\mu_j)\frac{d\mu_{j}}{dt}, \\		 
	\end{multline*}
	\begin{multline*}
		\frac{de^{sc}_{ij}}{dt} =\cos(\mu_{i})\cos(\mu_{j})\frac{dd^{sc}_{ij}}{dt}  +\cos(\mu_j)\sin(\mu_i)\frac{dd^{cc}_{ij}}{dt}  -\cos(\mu_i)\sin(\mu_j)\frac{dd^{ss}_{ij}}{dt} \\ -\sin(\mu_i)\sin(\mu_j)\frac{dd^{cs}_{ij}}{dt}   +d_{ij}^{cc}\cos(\mu_i)\cos(\mu_j)\frac{d\mu_{i}}{dt} 	
		-d_{ij}^{ss}\cos(\mu_i)\cos(\mu_j)\frac{d\mu_{j}}{dt} \\ -d_{ij}^{cs}\cos(\mu_i)\sin(\mu_j)\frac{d\mu_{i}}{dt}  -d_{ij}^{cs}\cos(\mu_j)\sin(\mu_i)\frac{d\mu_{j}}{dt}
		-d_{ij}^{sc}\cos(\mu_j)\sin(\mu_i)\frac{d\mu_{i}}{dt} \\
		-d_{ij}^{sc}\cos(\mu_i)\sin(\mu_j)\frac{d\mu_{j}}{dt}
		-d_{ij}^{cc}\sin(\mu_i)\sin(\mu_j)\frac{d\mu_{j}}{dt} +d_{ij}^{ss}\sin(\mu_i)\sin(\mu_j)\frac{d\mu_{i}}{dt}, \\		 
	\end{multline*}
	\begin{multline*}
		\frac{de^{ss}_{ij}}{dt} =\cos(\mu_{i})\cos(\mu_{j})\frac{dd^{ss}_{ij}}{dt}  +\cos(\mu_j)\sin(\mu_i)\frac{dd^{cs}_{ij}}{dt}  +\cos(\mu_i)\sin(\mu_j)\frac{dd^{sc}_{ij}}{dt} \\ +\sin(\mu_i)\sin(\mu_j)\frac{dd^{cc}_{ij}}{dt}   +d_{ij}^{cs}\cos(\mu_i)\cos(\mu_j)\frac{d\mu_{i}}{dt} 	
		+d_{ij}^{sc}\cos(\mu_i)\cos(\mu_j)\frac{d\mu_{j}}{dt} \\ +d_{ij}^{cc}\cos(\mu_i)\sin(\mu_j)\frac{d\mu_{i}}{dt}  +d_{ij}^{cc}\cos(\mu_j)\sin(\mu_i)\frac{d\mu_{j}}{dt}
		-d_{ij}^{ss}\cos(\mu_j)\sin(\mu_i)\frac{d\mu_{i}}{dt} \\
		-d_{ij}^{ss}\cos(\mu_i)\sin(\mu_j)\frac{d\mu_{j}}{dt}
		-d_{ij}^{cs}\sin(\mu_i)\sin(\mu_j)\frac{d\mu_{j}}{dt} -d_{ij}^{sc}\sin(\mu_i)\sin(\mu_j)\frac{d\mu_{i}}{dt}.\\		 
	\end{multline*}
	Consequently, by solving this system for $\left(\frac{d\mu_{i}}{dt},\frac{d\kappa_{i}}{dt}, \frac{dd^{cc}_{ij}}{dt}, \frac{dd^{cs}_{ij}}{dt}, \frac{dd^{sc}_{ij}}{dt}, \frac{dd^{ss}_{ij}}{dt} \right)$ in terms of $\left(\frac{d\eta^c_{i}}{dt},\frac{d\eta^s_{i}}{dt},\frac{de^{cc}_{ij}}{dt}, \frac{de^{cs}_{ij}}{dt}, \frac{de^{sc}_{ij}}{dt}, \frac{de^{ss}_{ij}}{dt}\right)$, we express the circular mean, concentration and precision parameter time derivatives in terms of canonical exponential parameters time derivatives. Then, the coordinate change between these two tangent spaces is
	\begin{subequations}\label{eq:mudt}
		\begin{equation}
			\frac{d\mu_{i}}{dt}= -\frac{\sin(\mu_{i})}{\kappa_{i}}\frac{d\eta^c_{i}}{dt}  + \frac{\cos(\mu_{i})}{\kappa_{i}}\frac{d\eta^s_{i}}{dt},
		\end{equation}
		\begin{equation}
			\frac{d\kappa_{i}}{dt}= \cos(\mu_{i})\frac{d\eta^c_{i}}{dt}  + \sin(\mu_{i})\frac{d\eta^s_{i}}{dt},
		\end{equation}
		\begin{multline}
			\frac{dd^{cc}_{ij}}{dt} = 
			d_{ij}^{cs}\frac{d\mu_{j}}{dt} 
			+ d_{ij}^{sc}\frac{d\mu_{i}}{dt} 
			+ \cos(\mu_i)\cos(\mu_j)\frac{de^{cc}_{ij}}{dt} 
			+ \cos(\mu_i)\sin(\mu_j)\frac{de^{cs}_{ij}}{dt} \\ 
			+ \cos(\mu_j)\sin(\mu_i)\frac{de^{sc}_{ij}}{dt} 
			+ \sin(\mu_i)\sin(\mu_j)\frac{de^{ss}_{ij}}{dt}	,
		\end{multline}
		\begin{multline}
			\frac{dd^{cs}_{ij}}{dt} = 
			d_{ij}^{ss}\frac{d\mu_{i}}{dt} 
			- d_{ij}^{cc}\frac{d\mu_{j}}{dt} 
			+ \cos(\mu_i)\cos(\mu_j)\frac{de^{cs}_{ij}}{dt} 
			- \cos(\mu_i)\sin(\mu_j)\frac{de^{cc}_{ij}}{dt} \\ 
			+ \cos(\mu_j)\sin(\mu_i)\frac{de^{ss}_{ij}}{dt} 
			- \sin(\mu_i)\sin(\mu_j)\frac{de^{sc}_{ij}}{dt}	,	
		\end{multline}
		\begin{multline}
			\frac{dd^{sc}_{ij}}{dt} = 
			d_{ij}^{ss}\frac{d\mu_{j}}{dt} 
			- d_{ij}^{cc}\frac{d\mu_{i}}{dt} 
			+ \cos(\mu_i)\cos(\mu_j)\frac{de^{sc}_{ij}}{dt} 
			- \cos(\mu_j)\sin(\mu_i)\frac{de^{cc}_{ij}}{dt} \\ 
			+ \cos(\mu_i)\sin(\mu_j)\frac{de^{ss}_{ij}}{dt} 
			- \sin(\mu_i)\sin(\mu_j)\frac{de^{cs}_{ij}}{dt}		,
		\end{multline}
		\begin{multline}
			\frac{dd^{ss}_{ij}}{dt} = 
			-d_{ij}^{sc}\frac{d\mu_{j}}{dt} 
			- d_{ij}^{cs}\frac{d\mu_{i}}{dt} 
			+ \cos(\mu_i)\cos(\mu_j)\frac{de^{ss}_{ij}}{dt} 
			- \cos(\mu_j)\sin(\mu_i)\frac{de^{cs}_{ij}}{dt} \\ 
			- \cos(\mu_i)\sin(\mu_j)\frac{de^{sc}_{ij}}{dt} 
			+ \sin(\mu_i)\sin(\mu_j)\frac{de^{cc}_{ij}}{dt}		,
		\end{multline}
	\end{subequations}
	for $i,j=1,\dots ,n$ and $\bm{\kappa}>0$ where $\eta_i^c$ and $\eta_i^s$ are canonical parameters associated with $\kappa_i\cos(\mu_i)$ in \cref{eq:etamap1}  and $\kappa_i\sin(\mu_i)$ respectively, $e_{ij}^{cs}, \ e_{ij}^{sc}, \ e_{ij}^{cc}, \ e_{ij}^{ss}$ are elements of the $\mathbf{E}^{cs}, \ \mathbf{E}^{sc}, \ \mathbf{E}^{cc}, \ \mathbf{E}^{ss}$ submatrices of $\mathbf{E}$ in \cref{eq:Epar1} and $d_{ij}^{cs}, \ d_{ij}^{sc}, \ d_{ij}^{cc}, \ d_{ij}^{ss}$ are elements of the $\mathbf{D}^{cs}, \ \mathbf{D}^{sc}, \ \mathbf{D}^{cc}, \ \mathbf{D}^{cc}$ submatrices of of the precision matrix $\mathbf{D}$ in \cref{eq:originallD1}.
	
	Using the canonical parametrization time derivatives given by the flow associated with the natural gradient Monte-Carlo estimates \cref{eq:natGradEst} 
	\begin{gather*}
		\frac{d\bm{\eta}}{dt}= \widehat{\widetilde{\nabla}_{\bm{\eta}} J(\bm{\eta},\textbf{E})}, \\
		\frac{d\textbf{E}}{dt}= \widehat{\widetilde{\nabla}_{\textbf{E}} J(\bm{\eta},\textbf{E})},
	\end{gather*}
	and using the coordinate transformations \cref{eq:mudt}, provide a flow associated with the circular mean, concentration and interaction parameters
	\begin{gather*}
		\frac{d\bm{\mu}}{dt}=\widehat{\widetilde{\nabla}_{\bm{\mu}} J(\bm{\mu},\bm{\kappa}, \textbf{D})}, \\
		\frac{d\bm{\kappa}}{dt}=\widehat{\widetilde{\nabla}_{\bm{\kappa}} J(\bm{\mu},\bm{\kappa}, \textbf{D})},\\
		\frac{d\textbf{D}}{dt}=\widehat{\widetilde{\nabla}_{\textbf{D}} J(\bm{\mu},\bm{\kappa}, \textbf{D})}.
	\end{gather*} 
	Thus the updated equations are then given by
	\begin{subequations}\label{eq:muUpdate}
		\begin{gather}
			\bm{\mu}^{t+1}= \bm{\mu}^{t} +\gamma_{\bm{\mu}}\widehat{\widetilde{\nabla}_{\bm{\mu}} J(\bm{\mu}^{t},\bm{\kappa}^{t}, \textbf{D}^t)}, \\
			\bm{\kappa}^{t+1}=\bm{\kappa}^{t}+\gamma_{\bm{\kappa}}\widehat{\widetilde{\nabla}_{\bm{\kappa}} J(\bm{\mu}^{t},\bm{\kappa}^{t}, \textbf{D}^t)},\\
			\mathbf{D}^{t+1}=\mathbf{D}^{t} + \gamma_{\mathbf{D}}\widehat{\widetilde{\nabla}_{\mathbf{D}} J(\bm{\mu}^{t},\bm{\kappa}^{t}, \textbf{D}^t)},
		\end{gather}
	\end{subequations}
	where $\gamma_{\bm{\mu}}, \ \gamma_{\bm{\kappa}}, \gamma_{\mathbf{D}} $ are the respective learning rates for each parameter group.

	\subsubsection{Adaptive learning rates}
	\label{subsubsec:adaptivelearningrates}
	To further stabilize the dynamical system given \cref{eq:muUpdate}, we implement a method to adaptively adjust learning rates proposed in \cite{silva} by changing each variable's learning rate individually by comparing gradient directions at two consecutive steps. If the general trend is the same, the learning rate is increased to accelerate convergence. If the change in the path is significant, the learning rate is decreased to allow steps with smaller granularity.
	
	In our setting, instead of comparing gradients directly, we compare parameter update differences given by	
	\begin{gather*}
		\Delta \mu_{i}^{t} = \mu_{i}^{t} - \mu_{i}^{t-1}, \\
		\Delta \kappa_{i}^{t} = \kappa_{i}^{t} - \kappa_{i}^{t-1}, \\
		\Delta d_{ij}^{t} = d_{ij}^{t} - d_{ij}^{t-1},
	\end{gather*}
	at times $t$ and $t-1$ for $i,j=1,\ldots,n$, and we set the upper bound for the adaptive learning rates to the initial learning rate $\gamma^0$ since we want only more fine-tuned learning rates when the algorithm has already located an optimum basin.   The adaptive learning rate update equations then take the following form
	\begin{subequations}\label{eq:adaEq}
		\begin{gather}
			\gamma^{t}_{\mu_{i}}=\left \{\begin{matrix}
				\min \left \{c_{up}\gamma^{t-1}_{\mu_{i}} \ ,\ \gamma^{0}_{\mu_{i}} \right\} & \text{if} &  \text{sgn}\left(\Delta \mu_{i}^{t}\right)=\text{sgn}\left(\Delta \mu_{i}^{t-1}\right), \\ 
				c_{down}\gamma^{t-1}_{\mu_{i}} & \text{if} & \text{sgn}\left(\Delta \mu_{i}^{t}\right)\neq\text{sgn}\left(\Delta \mu_{i}^{t-1}\right),
			\end{matrix} \right.  \\
			\gamma^{t}_{\kappa_{i}}=\left \{\begin{matrix}
				\min \left \{c_{up}\gamma^{t-1}_{\kappa_{i}} \ ,\ \gamma^{0}_{\kappa_{i}} \right\} & \text{if} &  \text{sgn}\left( \Delta \kappa_{i}^{t} \right)=\text{sgn}\left(\Delta \kappa_{i}^{t-1}\right), \\ 
				c_{down}\gamma^{t-1}_{\kappa_{i}} & \text{if} & \text{sgn}\left(\Delta \kappa_{i}^{t}\right)\neq\text{sgn}\left(\Delta \kappa_{i}^{t-1}\right),\\
			\end{matrix} \right. \\
			\gamma^{t}_{d_{ij}}=\left \{\begin{matrix}
				\min \left \{c_{up}\gamma^{t-1}_{d_{ij}} \ ,\ \gamma^{0}_{d_{ij}} \right\} & \text{if} &  \text{sgn}\left(\Delta d_{ij}^{t}\right)=\text{sgn}\left(\Delta d_{ij}^{t-1}\right), \\ 
				c_{down}\gamma^{t-1}_{d_{ij}} & \text{if} & \text{sgn}\left(\Delta d_{ij}^{t}\right)\neq\text{sgn}\left(\Delta d_{ij}^{t-1}\right) ,\\
			\end{matrix} \right. 
		\end{gather}
	\end{subequations}
	for $i,j=1,\ldots ,n$ where $\text{sign}(\cdot)$ is the signum function and $c_{up}>1, \ c_{down}<1$ are real positive constants.
	
	Additionally, we use momentum constants $\alpha_{\bm{\mu}}>0$,  $\alpha_{\bm{\kappa}}>0$ and $\alpha_{\mathbf{D}}>0$ in the update equations by setting
	\begin{gather*}
		m_{\mu_{i}}^t=\widehat{\widetilde{\nabla}_{\mu_{i}} J(\bm{\mu}^{t},\bm{\kappa}^{t}, \textbf{D}^t)} + \alpha_{\bm{\mu}}m_{\mu_{i}}^{t-1} ,\\
		m_{\kappa_{i}}^t=\widehat{\widetilde{\nabla}_{\kappa_{i}} J(\bm{\mu}^{t},\bm{\kappa}^{t}, \textbf{D}^t)} + \alpha_{\bm{\kappa}}m_{\kappa_{i}}^{t-1} ,\\
		m_{d_{ij}}^t=\widehat{\widetilde{\nabla}_{d_{ij}} J(\bm{\mu}^{t},\bm{\kappa}^{t}, \textbf{D}^t)} + \alpha_{\mathbf{D}}m_{d_{ij}}^{t-1},
	\end{gather*}
	to further aid the trajectory stabilization with the final form of update equations being
	\begin{gather*}
		\mu_{i}^{t+1}=\mu_{i}^{t}+\gamma_{\mu_{i}}^t m_{\mu_{i}}^t ,\\
		\kappa_{i}^{t+1}=\kappa_{i}^{t}+\gamma_{\kappa_{i}}^t m_{\kappa_{i}}^t, \\
		d_{ij}^{t+1}=d_{ij}^{t}+\gamma_{d_{ij}}^t m_{d_{ij}}^t .
	\end{gather*}

	\subsection{Fisher metric spectral radius scaling}
	\label{sec:spectralScaling}
	Although the sample covariance matrix estimator used in the estimation of the Fisher metric expressed in the canonical exponential parameters in \cref{eq:natGradEst} is considered efficient and unbiased, this may not be the case when the intrinsic geometry of the space of positive definite symmetric matrices is taken into account \cite{trees2007}, possibly resulting in instabilities of the entropic trust region updates with subsequent large fluctuations in the parameter trajectories when the number of samples is not adequately large. 
	
	To improve the stability of the gradients in \cref{eq:updateIGO}, we set a new metric tensor by scaling the inverse of the Fisher matrix by its spectral radius according to
	\begin{equation*}
		\mathcal{I}^{\lambda_{\text{min}}}_{\bm{\theta}}=\frac{1}{\lambda_{\text{min}}}\mathcal{I}_{\bm{\theta}},
	\end{equation*}   
	where $\lambda_{\text{min}}$ is the smallest eigenvalue of $\mathcal{I}_{\bm{\theta}}$ or the inverse of the spectral radius of $\mathcal{I}_{\bm{\theta}}^{-1}$. Then the entropic trust region step size \cref{eq:truststep} takes the form
	\begin{equation}\label{eq:specTtrustStep}
		\delta\bm{\theta}^t=\Delta^{t}\frac{ \left(\mathcal{I}^{\lambda_{\text{min}}}_{\bm{\theta}^t} \right)^{-1} \nabla_{\bm{\theta}} J(\bm{\theta}^{t})}{\sqrt{{\nabla_{\bm{\theta}} J(\bm{\theta}^{t})}^{T}\left(\mathcal{I}^{\lambda_{\text{min}}}_{\bm{\theta}^t} \right)^{-1} \nabla_{\bm{\theta}} J(\bm{\theta}^{t})}}.
	\end{equation}	
	
	The spectral radius scaling can be regarded as an adaptive trust region radius $\Delta^{t}$ in \cref{eq:specTtrustStep} since
	\begin{equation*}
		\frac{\left(\mathcal{I}^{\lambda_{\text{min}}}_{\bm{\theta}^t} \right)^{-1} \nabla_{\bm{\theta}} J(\bm{\theta}^{t})}{\sqrt{{\nabla_{\bm{\theta}} J(\bm{\theta}^{t})}^{T}\left(\mathcal{I}^{\lambda_{\text{min}}}_{\bm{\theta}^t} \right)^{-1} \nabla_{\bm{\theta}} J(\bm{\theta}^{t})}}=\sqrt{\lambda_{\text{min}}}\frac{\left(\mathcal{I}_{\bm{\theta}^t} \right)^{-1} \nabla_{\bm{\theta}} J(\bm{\theta}^{t})}{\sqrt{{\nabla_{\bm{\theta}} J(\bm{\theta}^{t})}^{T}\left(\mathcal{I}_{\bm{\theta}^t} \right)^{-1} \nabla_{\bm{\theta}} J(\bm{\theta}^{t})}},
	\end{equation*}
	Finally, the entropic trust region update equations \cref{eq:updateIGO} change to
	\begin{equation*}
		\bm{\theta}^{t+1}=\bm{\theta}^{t}+\sqrt{\lambda_{\text{min}}}\frac{ \widetilde{\nabla} J(\bm{\theta}^{t})}{\parallel\widetilde{\nabla} J(\bm{\theta}^{t})\parallel}.
	\end{equation*}

	\subsection{Hyperparameters}
	\label{sec:algparameters}
	Multiple hyperparameters need to be set for the optimization algorithm to work accurately. In the following few paragraphs, we address the specific values we use for the problem of dense plane group packings of convex polygons. All the values were determined experimentally on the problem of densest $p2$-packing of a regular octagon, described in more detail in \cref{sec:proofofconcept}. 
	
	Initially, the number of samples used Monte--Carlo estimates of the Fisher metric tensor \cref{eq:expF} and the expected fitness gradient \cref{eq:adaSelEst} was set. As was already mentioned in section \ref{sec:spectralScaling}, the sample covariance matrix estimator is known to be ill-conditioned when the ratio $c=\frac{p}{N}$ (where $p$ is the rank of the covariance matrix and $N$ the number of samples) is not negligible. Since the computationally most demanding part of the optimization schedule is the extended multivariate von Mises Gibbs sampler \cref{sec:GibbsSampler}, the trade-off between the accuracy and stability of the algorithm's trajectories and speed of computations has to be balanced. We experimentally determined $c=0.12$ to be an acceptable balance between speed and stability.
	
	Concerning the adaptive selection quantile parameters in \cref{eq:adaSQ}, we set the starting $q$-quantile to $q_0=\frac{N}{100}$ and the schedule parameter $\beta$ in such a way that if after each iteration the cosine of the angle between consecutive steps is $0$, corresponding to no change in the direction of the updates, the adaptive selection $q$-quantile would reach the $N$th quantile after $2,000$ iterations; that is, $\beta=\ln\left(100\right)/2,000$. 
	
	Hyperparameter settings related to the adaptive learning rates described in section \cref{subsubsec:adaptivelearningrates} were obtained by a localized instance of grid search. First, the rate of change of the learning rates in \cref{eq:adaEq} was set to $c_{up}= 1.1$ and $c_{down}= 0.9$. Subsequently, we constructed the $6$D grid 
	\begin{equation*}
		\left\{\left(\gamma_{\bm{\mu}}^0,\gamma_{\bm{\kappa}}^0,\gamma_{\mathbf{D}}^0\right)|\gamma_{\mathbf{*}}^0 \in\left\{0.25;0.5;0.75;1\right\}\right\} \times \left\{\left(\alpha_{\bm{\mu}}^0,\alpha_{\bm{\kappa}}^0,\alpha_{\mathbf{D}}^0\right)|\alpha_{\mathbf{*}}^0 \in\left\{0;0.25;0.5;0.75\right\}\right\},
	\end{equation*}
	with $* \in \left\{\bm{\mu},\bm{\kappa},\mathbf{D} \right\}$, where $\gamma_{\mathbf{*}}^0$ denotes the learning rate and $\alpha_{\mathbf{*}}^0$ momentum parameters associated with the extended multivariate von Mises distribution parameters in \cref{eq:fmvm1} and performed $20$ optimization runs for each grid node. After the initial search, we increased the granularity of the grid by halving the size of the line segment between two neighbouring nodes and performed 20 optimization runs but only for the nearest neighbours of the node that attained the highest mean from the objective function of the best solutions found in each run in the initial grid search. We iterated this process of exploring only the nearest neighbour of the node that attained the best mean from all $20$ runs in the previous iteration. When no better combination of parameters was found, we increased the granularity of the grid by halving the size of the grid line segments and repeated the exploration process on a more fine-grained grid, starting with the node with the highest mean of the objective.
	
	The exponential decay $c_{\epsilon}$ of the $\bm{ \epsilon}^r$ neighbourhood in \cref{eq:epsilonC} used in the refining part of the algorithm has to be balanced between accuracy and speed. If $c_{\epsilon}$ is set too high, the neighbourhoods converge too fast, and it is possible that some $s$ the $\bm{ \epsilon}^r$ neighbourhoods of the best solution found at the $r$th run do not contain the optimal solution for all $r\geq s$. On the other hand, if $c_{\epsilon}$ is too low, it increases the number of refining runs and impacts the algorithm's efficiency. $c_{\epsilon}$ was set to $1.2$ by evaluating the difference from the theoretical packing density $\cref{eq:packDiff}$ with the objective of accuracy being $\Delta_{\mathcal{K}_{p2}}<10^{-7}$. 
	
	All the values of the hyperparameters used in computations presented in the experimental results \cref{sec:proofofconcept} and \cref{sec:additionalPackings} are listed in \cref{tab:1}.
	
	\begin{table}
		\centering
		\begin{tabular}{|c|c|}
			\hline
			$\gamma_{\bm{\mu}}^0$& 0.140625 \\
			\hline
			$\gamma_{\bm{\kappa}}^0$& 0.171875 \\
			\hline
			$\gamma_{\mathbf{D}}^0$& 0.21875 \\
			\hline
			$c_{up}$& 1.1 \\
			\hline
			$c_{down}$& 0.9 \\
			\hline
			$\alpha_{\bm{\mu}}$& 0.7109375 \\
			\hline
			$\alpha_{\bm{\kappa}}$& 0.1953125 \\
			\hline
			$\alpha_{\mathbf{D}}$& 0.578125 \\
			\hline
			$q_0$& $\frac{N}{100}$ \\
			\hline
			$\beta$& $\ln\left(100\right)/2,000$ \\
			\hline
			$c$ & 0.12 \\
			\hline
			Number of iterations& 8,000 \\
			\hline
			$c_{\epsilon}$ & 1.2 \\ 
			\hline
		\end{tabular}
		\caption{Used hyperparameter settings.}
		\label{tab:1}
	\end{table}

	\subsection{Parallelizing computations}
	\label{sec:parallel}
	All the computations were performed in MATLAB R2021b. To accelerate computations, we utilize the architecture of modern microprocessors using the Parallel Computing toolbox. In the following paragraphs, we provide details on our parallel algorithm implementation.
	
	Gibbs sampling for the exponential rewriting of the extended multivariate von Mises distribution described in \cref{sec:GibbsSampler} is suitable for parallel computations on the CPU since the Gibbs sampling simulation involves evaluating scalar operations, non--elementary function calls and branching. We divide the number of samples $\lambda$ needed for the trust region step size estimation \cref{eq:truststep} among $W$ number of workers available. Then every worker generates $\frac{\lambda}{T}$ independent extended multivariate von Mises with parameter $(\bm{\eta}^t,\textbf{E}^t)$ distributed samples.
	
	After workers complete the Gibbs sampling simulations, the results are combined and transferred to the GPU memory, and all subsequent computations are done on the GPU. GPU computations are particularly efficient for the objective function \cref{eq:density}, and constraint violation \cref{eq:distPack} evaluations since both involve just elementary function evaluations that are utilized using MATLAB's inherent vectorization design and are efficiently computed on the GPU for all $\lambda$ candidate solutions in parallel.
	
	\begin{figure}
		\centering
		\begin{subfigure}{0.48\textwidth}
			\includegraphics[width=1\linewidth]{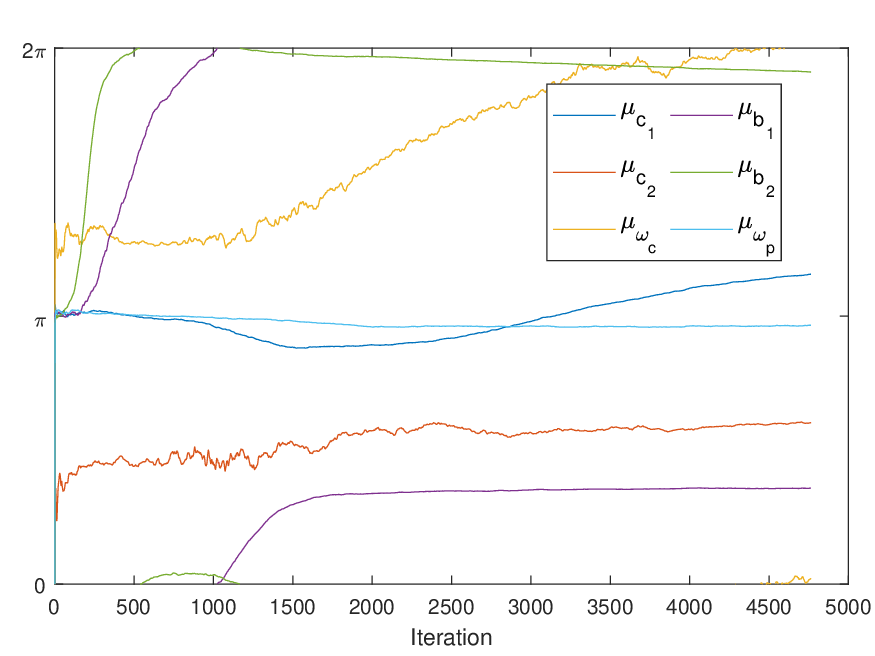}
		\end{subfigure}		
		\begin{subfigure}{0.48\textwidth}
			\includegraphics[width=1\linewidth]{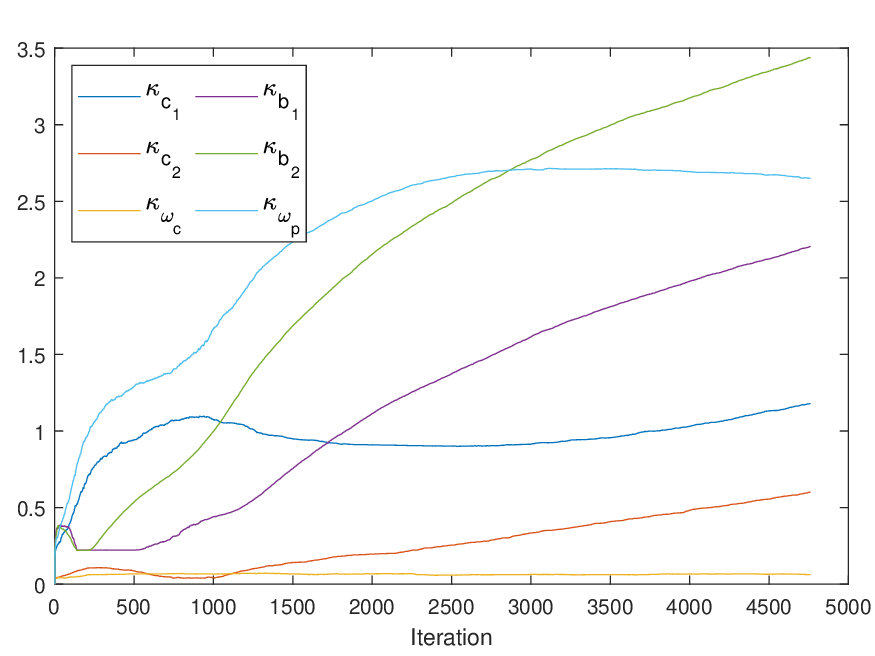}
		\end{subfigure}
		\caption{Trajectories of the extended multivariate von Mises distribution (left) circular mean $\bm{\mu}$ and (right) concentration $\bm{\kappa}$ parameters.}
		\label{fig:paperMuKappa}
	\end{figure}
	
	\section{Entropic trust region trajectories}
	\label{sec:trajectories}
\begin{figure}	
		\centering
		\begin{subfigure}{0.48\textwidth}
			\includegraphics[width=1\linewidth]{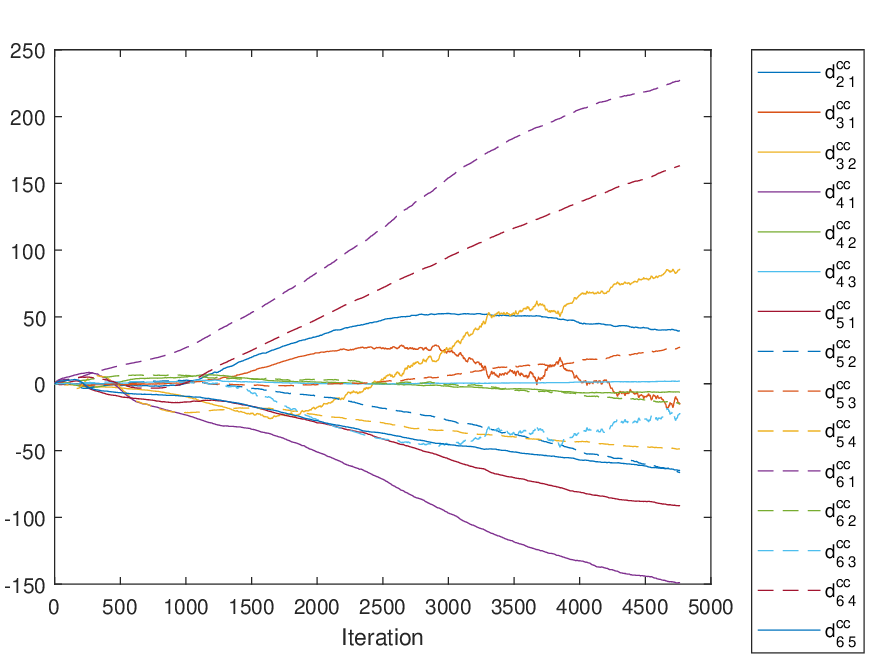}
		\end{subfigure}	
		\begin{subfigure}{0.48\textwidth}
			\includegraphics[width=1\linewidth]{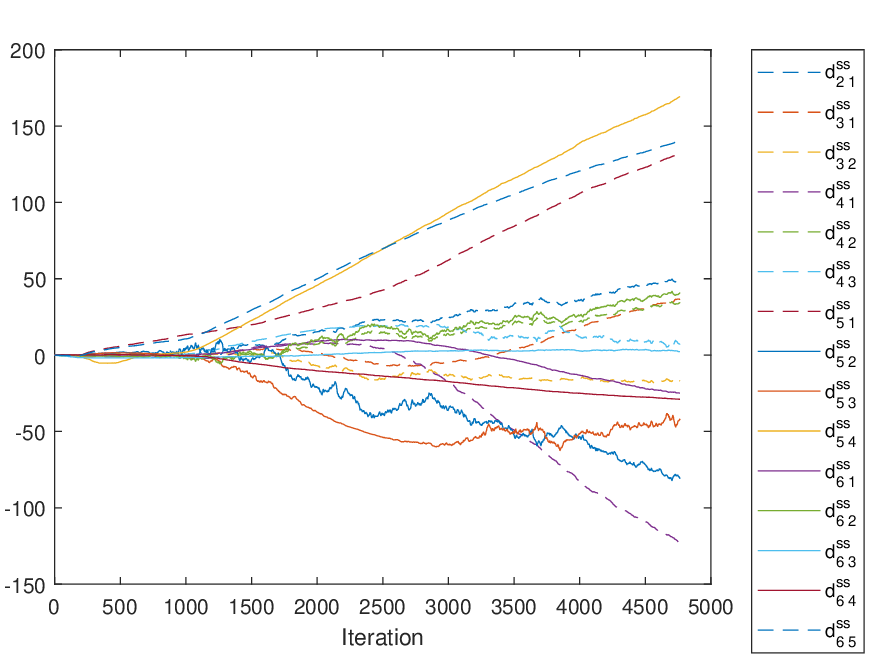}
		\end{subfigure}
		
		\begin{subfigure}{0.48\textwidth}
			\includegraphics[width=1\linewidth]{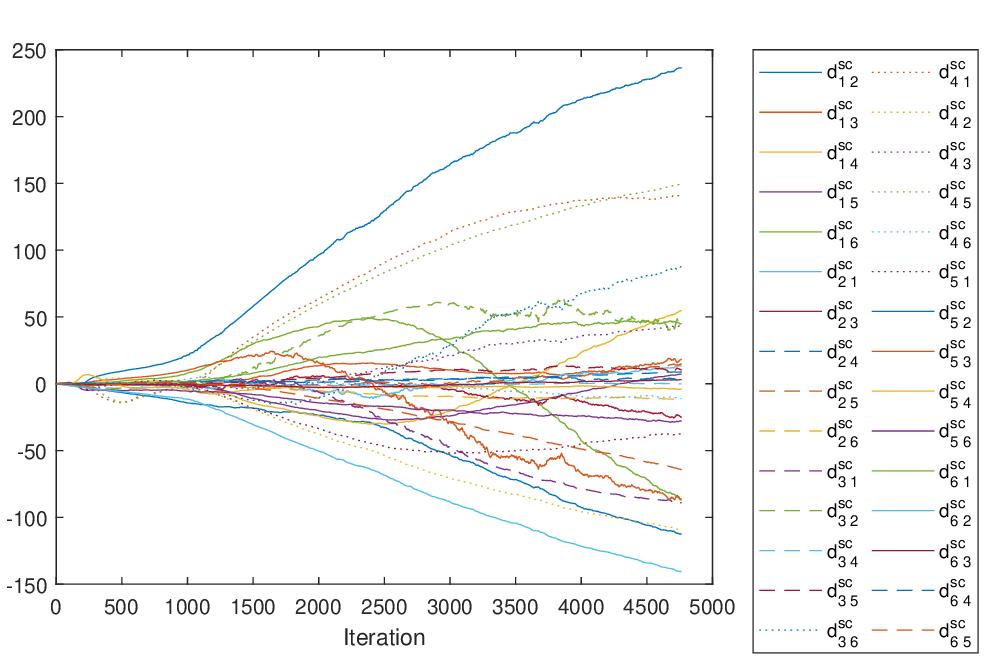}
		\end{subfigure}
		\caption{Trajectories of the extended multivariate von Mises distribution $\textbf{D}$ matrix elements: (top left) $\textbf{D}^{cc}$ representing $\cos\left(\theta_i\right)-\cos\left(\theta_j\right)$ interactions, (top right) elements of $\textbf{D}^{ss}$ representing $\sin\left(\theta_i\right)-\sin\left(\theta_j\right)$ interactions and (bottom) elements of $\textbf{D}^{sc}$ representing $\sin\left(\theta_i\right)-\cos\left(\theta_j\right)$ interactions.}
		\label{fig:paperD}
	\end{figure}

	A visualization of the evolution of the distribution parameters is presented in \cref{fig:paperMuKappa} and \cref{fig:paperD}. As we already mention in \cref{sec:toroidalDist}, it is generally difficult to precisely interpret the distribution parameters. Nevertheless, some intuition about the algorithmic behaviour can be extracted from the entropic trust region trajectories. 
	
	Most of the circular mean $\bm{\mu}$ parameters stabilize after $3,000$ iterations of the algorithm with the exception of the parameter $\mu_{\omega_c}$ which is related to the angle of rotation of the octagon in the asymmetric unit. This is due to the rotational symmetry of the octagon that induces multiple global optima in the $\omega_c$ subspace of the optimization landscape. To be precise, for $\omega_c=\frac{2\pi \left(k-1\right)}{8}, \ k=1,\dots,7 $ the packing density \cref{eq:density} is equal. 
	
	Similar behaviour can be observed in the trajectories of the concentration parameters $\bm{\kappa}$. These parameters are analogous to the inverse of dispersion measures in descriptive statistics. As the algorithm progresses, the concentration of the distribution gradually increases, except $\kappa_{\omega_c}$, representing the angle of rotation of the octagon, and stagnates close to zero. This behaviour can also be observed in the $1$D histogram and $2$D projections of $600$ realizations of the output distribution involving $\omega_c$ variable in \cref{fig:samplehist}.
	
	The evolution of the interaction matrix $\textbf{D}$ is shown in \cref{fig:paperD} for $\cos - \cos$ interactions, $\sin - \sin$ interactions and  $\cos - \sin$ interactions. From the geometric interpretation of the exponential family selection quantile-based entropic trust region in \cref{sec:geometryEntTrust} perspective, we can see that the extended multivariate von Mises distribution is indeed moving towards a higher degree of interaction between angles. In the sense of unsupervised learning, the algorithm is learning a representation of the optimization landscape, and \cref{fig:samplehist} can be regarded as a visualization of this representation in the setting of $p2$-packings of regular octagons.

	\begin{figure}
		\centering
		\begin{subfigure}[!b]{0.48\textwidth}
			\includegraphics[width=1\linewidth]{paperMean}
		\end{subfigure}
		\begin{subfigure}[!b]{0.504\textwidth}
			\includegraphics[width=1\linewidth]{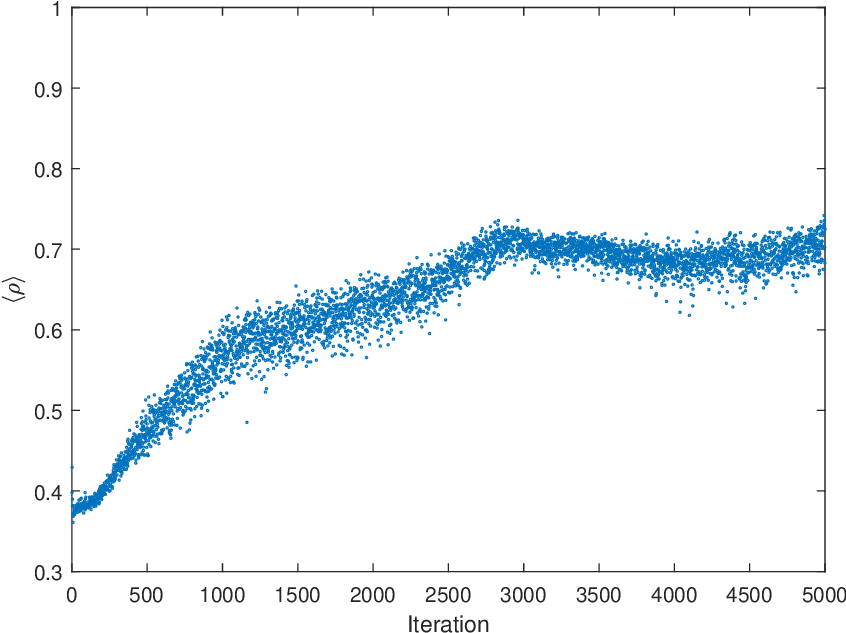}
		\end{subfigure}
		
		\begin{subfigure}[!b]{0.504\textwidth}
			\includegraphics[width=1\linewidth]{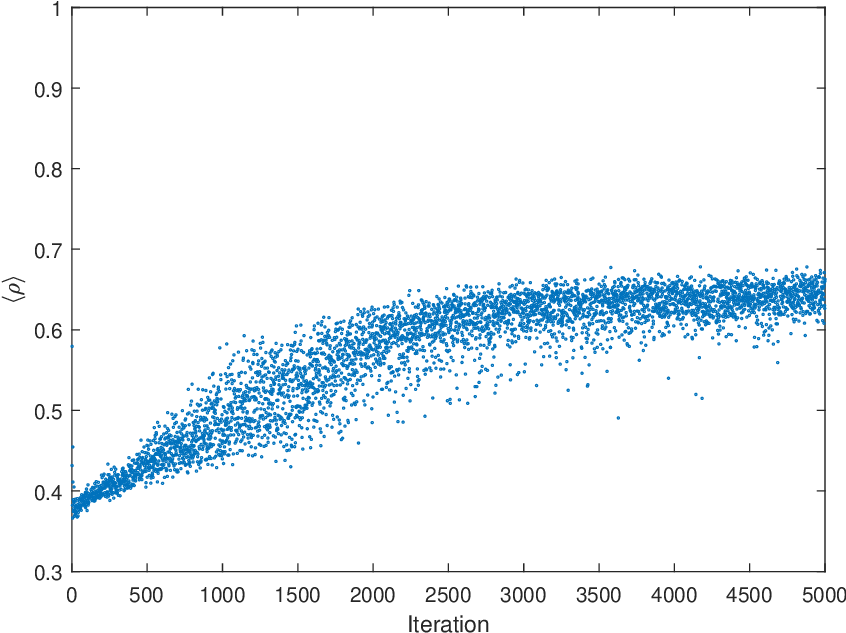}
		\end{subfigure}
		\caption{Evolution of average density in three runs initialized with the same seed. In each run, all the hyperparameters were fixed except the number of samples used in the Monte--Carlo estimates. (top left) $600$, (top right) $400$ and (bottom) $200$ realization from the extended multivariate von Mises distribution are used.}
		\label{fig:samplesComparison}
	\end{figure}

	\section{Influence of sampling pool on stability}
	\label{sec:samplingPool}
	The Fisher metric tensor $\mathcal{I}_{\bm{\theta}}$ used in the trust region step \cref{eq:truststep} has in the case of $p2$ group rank $r=72$. Based on the ratio $c=0.12$ of the Fisher matrix rank to the number of samples mentioned in \cref{sec:algparameters}, the current implementation of the algorithm uses $600$ samples generated at each iteration for the step size estimates in the densest $p2$-packing problem. Though less than unity, $c$ cannot be considered a negligible quantity and results in problems with the variance of the sample covariance matrix estimator of the Fisher matrix \cref{eq:expF}, even though we use adaptive learning rates and the Fisher metric scaling discussed in \cref{subsubsec:adaptivelearningrates} and \cref{sec:spectralScaling} to mitigate the effects of large estimator deviations. 
	
	To illustrate this issue, we compare three optimization runs with the random number generator initialized with the same seed and hyperparameter settings but a different number of samples generated at each iteration. \Cref{fig:samplesComparison} shows evolution of the average density $\left<\rho\right>$ where $600$, $400$ and $200$ are used in the trust region step size Monte--Carlo estimates \cref{eq:natGradEst}.   Visually, the run where $200$ samples are used is more scattered than in the cases with $400$ and $600$,  which implies larger fluctuations in the trust region updates between iterations. Numerically, the variance of average density first differences $\Delta^t\left<\rho\right> = \left<\rho\right>^t-\left<\rho\right>^{t-1}$ is $3.87 \ 10^{-4}$ in the case of $600$ samples, $4.5\ 10^{-4}$ for the $400$ samples, and $1.2\ 10^{-3}$ for the $200$ samples, and this supports our observations. The density of the best solution found in the case with the $600$ samples is $\rho \left(\mathcal{K}_{p2} \right)=0.897526117081202$, in the case of $400$ samples $\rho \left(\mathcal{K}_{p2} \right)=0.893840811317134$, and $\rho \left(\mathcal{K}_{p2} \right)=0.890457479411185$ for the $200$ samples. Although the variance of the average mean differences between the $600$ and $400$ sample instances is rather small, the highest density solution found during $5,000$ iterations was in the case of the former.

	\begin{figure}
		\centering
		\includegraphics[width=0.3\linewidth]{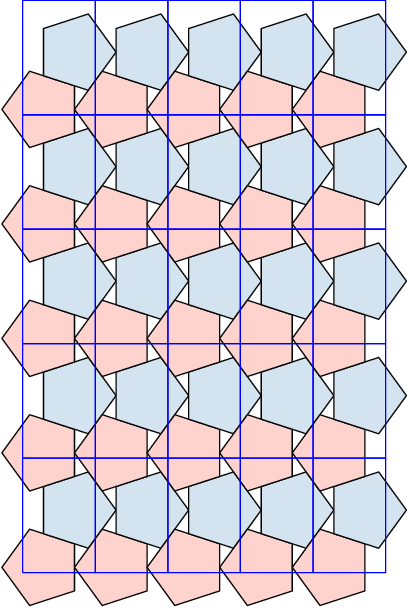}
		\includegraphics[width=0.55\linewidth]{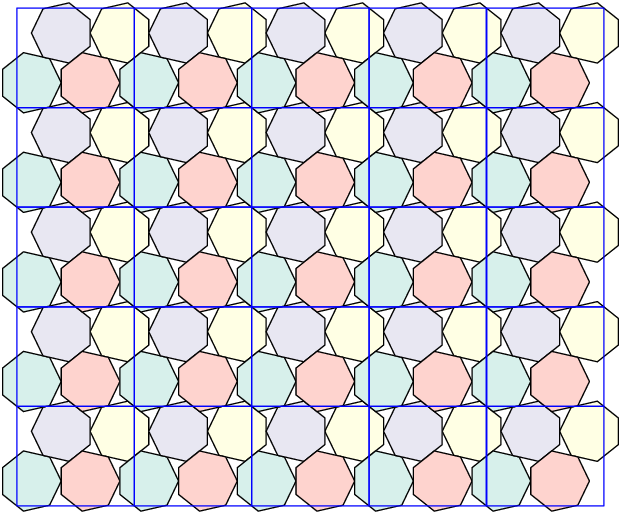}
		\caption{Visualization of $25$ cells of the output configurations of the densest (left) regular pentagon $pg$-packing and (right) regular heptagon $p2gg$-packing.}
		\label{fig:p2gon5final}
	\end{figure}
	
	\section{Additional plane group packings}
	\label{sec:additionalPackings}
	To demonstrate the robustness of the entropic trust region search, we present additional densest plane group experiments on a few selected convex polygons for which their densest packings are known. In all test cases, the hyperparameters are set according to \cref{tab:1}. Numerical comparisons of the results are presented in \cref{tab:tabel}.
	
	\subsection{$pg$-packing of regular pentagon}
	General optimal packing of congruent copies of the regular pentagon is believed to be a double lattice configuration with the density
	\begin{equation*}
		\rho_{\text{opt}}=\frac{5-\sqrt{5}}{3} \approx 0.92131067
	\end{equation*}
	\cite{halesSup}. In terms of plane groups, a double lattice is the group $p2$ already introduced in the $p2$-octagon packing case study in \cref{sec:proofofconcept}. The $p2$ group can be viewed as a collection of $2$ lattices related to each other by a $2$-fold rotational symmetry that permutes the two lattices. 
	
	Interestingly, our packing algorithm converged to the general optimal configuration when applied to search for the densest $pg$ packings of regular pentagons. The $pg$ plane group is not a semi-direct product of a lattice group and a point group since it contains a glide reflection symmetry operation that is neither an element of the point group nor the lattice group associated with $pg$. The crystal system is rectangular, meaning the angle $\omega_p$ between the primitive cell's basic vectors is fixed to $90^{\circ}$, which reduces the number of optimization variables to $5$ compared to the $p2$ plane group.
	
	The output configuration of the densest $pg$-packing of the regular pentagon optimization schedule is visualized in \cref{fig:p2gon5final} (left) with packing density $\rho \left(\mathcal{K}_{pg} \right)=0.9213-1060131385$ and the theoretical optimum difference $\Delta_{\mathcal{K}_{pg}}=7.2852 \ 10^{-8}$, implying that the densest packing of regular pentagon can be realized using a glide reflection symmetry instead of the $2$-fold rotational symmetry in $p2$ plane group.
	
	\subsection{$p2gg$-packing of regular heptagon}
	The optimal double lattice packing that is a $p2$-packing of the regular heptagon is
	\begin{equation*}
		\rho_{\text{opt}}=\frac{2}{97}\left(-111+492\cos\left(\frac{\pi}{7}\right)-356\cos^2\left(\frac{\pi}{7}\right)\right)\approx 0.89269068
	\end{equation*} 
	\cite{kuperberg} and is conjectured to be the general optimal packing of regular heptagons. 
	
	The output configuration of our densest $p2gg$-packing optimization schedule is shown in \cref{fig:p2gon5final} (right) with packing density $\rho \left(\mathcal{K}_{p2gg} \right)=0.89269066997639$ and theoretical optimum difference $\Delta_{\mathcal{K}_{p2gg}}=1.6150 \ 10^{-08}$, showing that the densest configuration of regular heptagons can be obtained as packing in a higher symmetry group than it was previously known.
	
	$p2gg$ plane group symmetry operations are given by a $2$-fold rotational symmetry and $2$ glide reflection along perpendicular mirror planes. The crystal system is the same as in the $pg$ group, that is rectangular, restricting the $pg$-packing configuration space dimension to $5$ and the fractional coordinates of the heptagon to $0\leq c_1 \leq \frac{1}{2}$ and $0\leq c_2 \leq \frac{1}{2}$.
	
	\subsection{$p4$-packing of non-regular pentagon}
\begin{figure}
		\centering
		\includegraphics[width=0.345\linewidth]{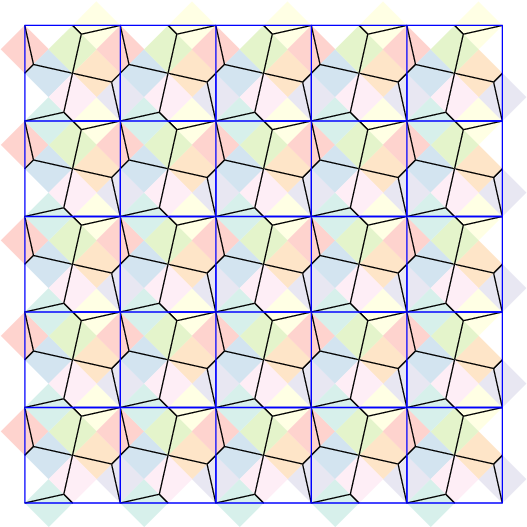}
		\includegraphics[width=0.45\linewidth]{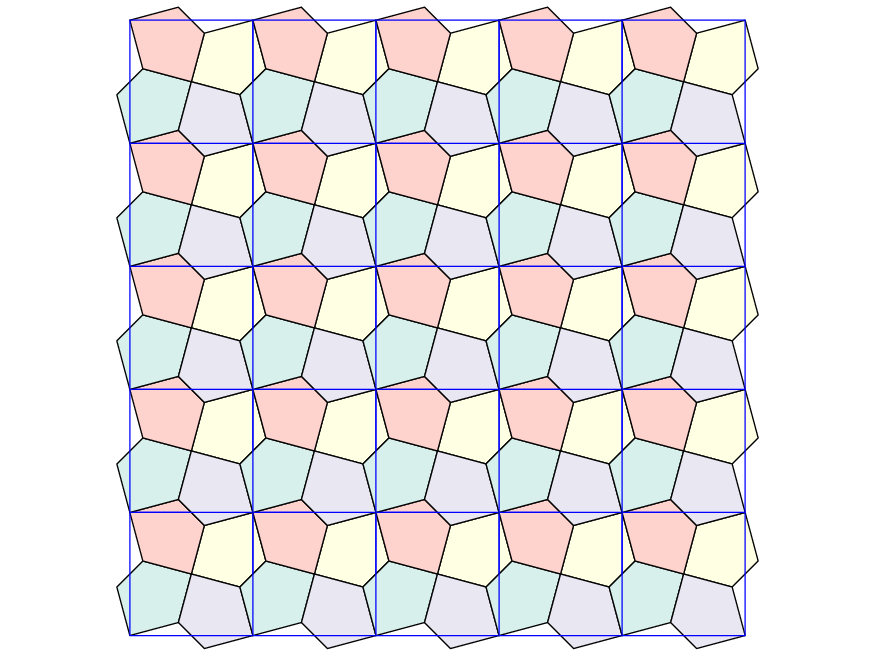}
		\caption{(Left) The Cairo tiling as a $p4gm$-packing of squares with a motif. (Right) a visualization of $25$ cells of the output configurations of the densest non-regular pentagon $p4$-packing.}
		\label{fig:p4CairoFinal}
	\end{figure}

	An non-regular pentagon with the sequence of internal angles $120^\circ, \ 120^\circ, \ 90^\circ, \ 120^\circ, \ 90^\circ$ tiles the Euclidean plane, and its tiling is called the Cairo pentagonal tiling.
	Although the full symmetry of the Cairo tiling is $p4gm$, it is not a $p4gm$-packing but rather a $p4$-packing. The $p4gm$ symmetry of the Cairo tiling is in fact a $p4gm$-packing of squares with a motif as shown in \cref{fig:p4CairoFinal} (left).
	
	Similarly to the $p2$ group, $p4$ is a semi-direct product of a point group with a $4$-fold rotational symmetry and a lattice group belonging to the square crystal system, imposing restrictions on the shape of the primitive cell. That is, the sizes of the lattice group generators $\textbf{b}_1$ and $\textbf{b}_2$ are equal, and the angle between $\textbf{b}_1$ and $\textbf{b}_2$ fixed to $\omega_p = \frac{\pi}{2} $. \cref{fig:p4CairoFinal} (right) shows the output configuration of the ETRPA procedure with density $\rho \left(\mathcal{K}_{p4} \right)=0.99999999503997$ and theoretical optimum difference $\Delta_{\mathcal{K}_{p4}} = 4.9600 \ 10^{-9}$.
	
	\subsection{$p3$-packing of regular hexagon}
	\label{sec:p3Pack}
	
	The regular hexagon tiles the $2$D Euclidean space, which corresponds to an optimal packing density of $\rho_{\text{opt}}=1$.
	
\begin{figure}[t]
		\centering
		\includegraphics[trim={0 -30 0 0},width=0.53\linewidth]{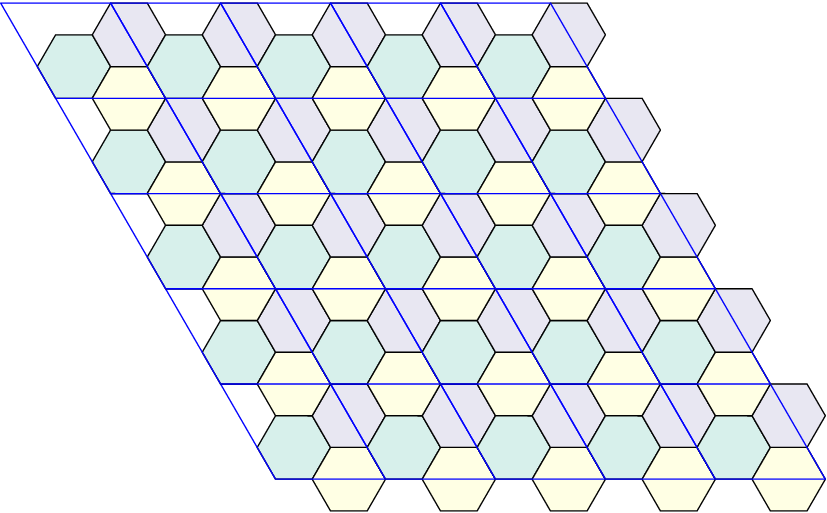}
		
		\includegraphics[width=0.47\linewidth]{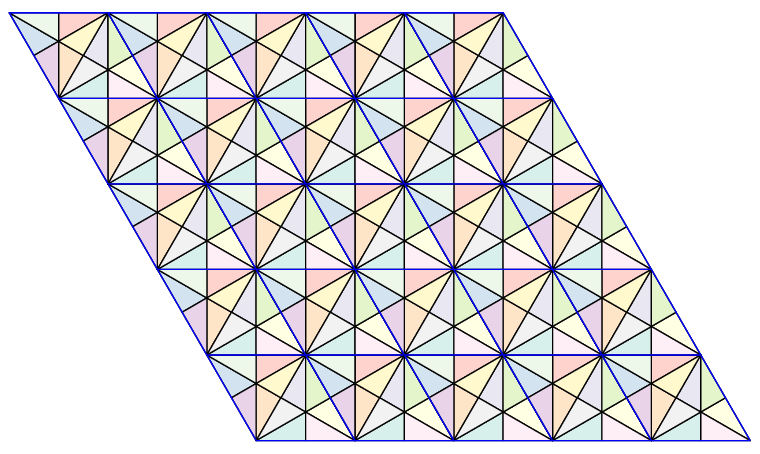}
		\includegraphics[width=0.45\linewidth]{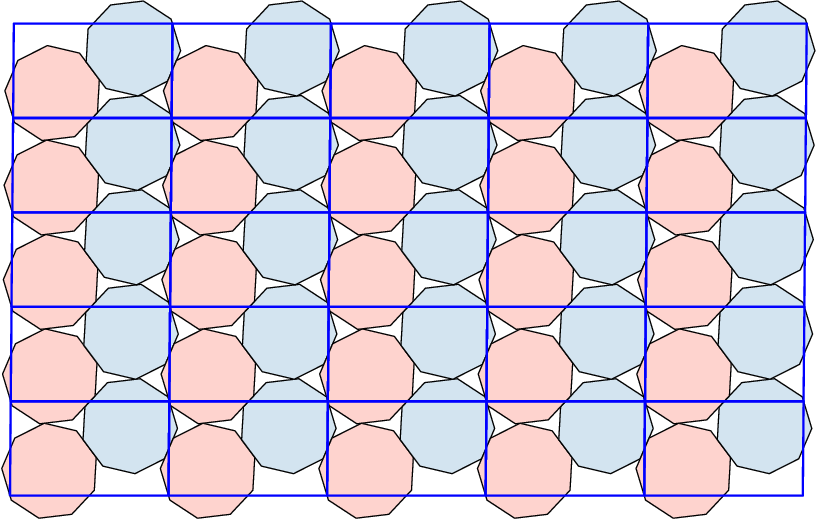}
		\caption{Visualization of $25$ cells of the output configurations of the densest (top) regular hexagon $p3$-packing, (bottom left) $30-60-90$ triangle $p6mm$-packing and (bottom right) regular enneagon $p2$-packing.}
		\label{fig:p6mtrianglefinal}
	\end{figure}

	\cref{fig:p6mtrianglefinal} (top) shows the output configuration when we applied the entropic trust region packing algorithm to the regular hexagon in the plane group $p3$. $p3$ is a semi-direct product of a lattice group and a point group consisting of a $3$-fold rotational symmetry. The crystal system is hexagonal, which means that some restrictions on the shape of the primitive cell are imposed. Specifically, $\parallel\textbf{b}_1\parallel = \parallel\textbf{b}_2\parallel$, that is, the lattice group generators $\textbf{b}_1$ and $\textbf{b}_2$ are of equal size, and the angle between $\textbf{b}_1$ and $\textbf{b}_2$ is fixed to $120^{\circ}$, reducing the degrees of freedom of the packing problem to $4$. That is the fractional coordinates $0\leq c_1 \leq \frac{2}{3}$ and $0\leq c_2 \leq \frac{2}{3}$ of the triangle's centroid in the asymmetric unit, the angle $\omega_c$ of the rotation of the triangle, and the size of one of the edges of the primitive cell $\parallel\textbf{b}_1\parallel$. 
	
	The $p3$ group also forces additional linear constraints on the position of the triangle's centroid coordinates $c_1$ and $c_2$, which are
	\begin{align*}
		c_2-\min\{1-c_1,\frac{c_1}{2}+\frac{1}{2}\} 		   \leq 0, \\
		c_1-\frac{c_2}{2}-\frac{1}{2} \leq 0.
	\end{align*}  
	These restrictions are treated as additional inequality constraints in \cref{eq:constraints} and are incorporated within the penalty function \cref{eq:penalty} computations.
	
	The output packing density of our optimization schedule is $\rho \left(\mathcal{K}_{p1} \right)=0.99999993-380570$ with the optimal packing difference of $\Delta_{\mathcal{K}_{p1}}=6.6194 \ 10^{-08}$. In practice, this configuration is a tiling showcasing three-fold rotational symmetries of regular hexagonal tiling.
	
	\subsection{$p6mm$-packing of $30-60-90$ triangle}
	The $30-60-90$ triangle is an irregular polygon with internal angles of $30^\circ,\ 60^\circ$ and $90^\circ$. Since it tiles the Euclidean plane, its exact packing density is$\rho_{\text{opt}}=1$. The $p6mm$ group consist of a $6$-fold rotational symmetry and reflections through $6$ mirror planes. The crystal system is the same as in the $p3$ plane group mentioned earlier, which is  hexagonal. The $p6mm$ group induces additional linear constraints on the position of the triangle's centroid fractional coordinates $c_1$ and $c_2$ by
	\begin{align*}
		2c_1-c_2-1 		   \leq 0, \\
		-\frac{c_1}{2}+c_2 \leq 0.
	\end{align*}  
	and are treated as additional inequality constraints in \cref{eq:constraints} as in the case of $p3$-packing of hexagons \cref{sec:p3Pack} (bottom left).
	
	The output configuration of our optimization schedule applied to the densest $p6mm$-packing for the $30-60-90$ triangle search is shown in \cref{fig:p6mtrianglefinal} (bottom left) with packing density of $\rho \left(\mathcal{K}_{p6mm} \right)=0.99999999871467$ and optimal packing difference $\Delta_{\mathcal{K}_{p6mm}}=1.2853e-09$, confirming that the $30-60-90$ triangle constitutes the primitive cell of the $p6mm$ plane group.
	
	\subsection{$p2$-packing of regular enneagon}
	
	Additionally, we applied the entropic trust region to the densest $p2$-packing of a regular enneagon. Although theoretical optimal packing is unknown, the densest packing of regular enneagons with the $p2$ plane group symmetry, introduced in section \cref{sec:proofofconcept}, was reported $0.90103007842093$ in \cite{de2011}. The density of the output configuration shown \cref{fig:p6mtrianglefinal} (bottom right) is $\rho \left(\mathcal{K}_{p2} \right)=0.901030017272363$.
		
	\bibliographystyle{siamplain}
	\bibliography{references}	
	
\end{document}